\documentclass[12pt]{article} 
\usepackage[sectionbib]{natbib}
\usepackage{array,epsfig,fancyheadings,rotating}
\usepackage[]{hyperref}  
\usepackage{sectsty, secdot}
\usepackage{color,xcolor}
\usepackage{bm}
\usepackage{algorithm,algpseudocode,float,algorithmicx}
\sectionfont{\fontsize{12}{14pt plus.8pt minus .6pt}\selectfont}
\renewcommand{\theequation}{\thesection\arabic{equation}}
\subsectionfont{\fontsize{12}{14pt plus.8pt minus .6pt}\selectfont}
\newcommand{\overbar}[1]{\mkern 2.5mu\overline{\mkern-2.5mu#1\mkern-2.5mu}\mkern 2.5mu}
\newcommand{\kl}{\textrm{KL}} 
\usepackage{soul}
\newcommand{\bbE}{\mathbb{E}}
\def\mc{\mathcal}
\def\bs{\boldsymbol}

\newcommand{\ttE}{\mathtt{E}}
\newcommand{\ttF}{\mathtt{F}}
\def\vm{\boldsymbol{m}}

\def\R{\mathbb R}

\def\P{\mathbb P}
\def\H{\mathcal H}

\def\S{\mathcal S}
\newcommand{\bbP}{\mathbb{P}}
\def\vB{\boldsymbol{B}}
\def\E{\mathbb E}

\def\u{\boldsymbol{u}}
\def\X{\boldsymbol{X}}

\def\vZ{\boldsymbol{Z}}
\def\vu{\boldsymbol{u}}

\def\vX{\boldsymbol{X}}
\def\Cov{\mathrm{Cov}}
\def\m{\mathbf m}

\def\S{\mathcal{S}}
\def\mb{\mathbb}
\def\bbeta{\bm \beta}

\def\var{\mbox{var}}
\def\wh{\widehat}
\def\wt{\widetilde}
\def\mr{\mathrm}

\def\btheta{\boldsymbol{\theta}}

\def\l{\left(}
\def\r{\right)}

\def\lno{\left\|}
\def\rno{\right\|}

\usepackage{amsmath}
\usepackage{amssymb}
\usepackage{amsfonts}
\usepackage{multirow}
\usepackage{amsthm}
\usepackage{setspace}
\numberwithin{equation}{section}
\newtheorem{theorem}{Theorem}
\newtheorem{lemma}{Lemma}

\newtheorem{proposition}{Proposition}
\theoremstyle{definition}
\newtheorem{definition}{Definition}

\newtheorem{remark}{Remark}
\newtheorem{assumption}{Assumption}
\newcommand{\norm}[1]{\left\Vert #1 \right\Vert}

\usepackage{geometry}
\usepackage{authblk}
\begin{document}


\title{On the Optimality of Functional Sliced Inverse Regression}
\author[a]{Rui Chen}
\author[b]{Songtao Tian\thanks{Co-first author.}}
\author[c]{Dongming Huang}
\author[d]{Qian Lin}
\author[e]{Jun S. Liu\thanks{Corresponding author.}}
\affil[a]{Department of Statistics and Data Science, Tsinghua University
}
\affil[b]{Department of Mathematical Sciences, Tsinghua University
}
\affil[c]{Department of Statistics and Data Science, National University of Singapore
}
\affil[d]{Department of Statistics and Data Science, Tsinghua University;~Beijing Academy of Artificial Intelligence, Beijing, 100084, China
}
\affil[e]{Department of Statistics, Harvard University, 02138
}
\renewcommand*{\Affilfont}{\small\it} 
\renewcommand\Authands{ and } 
\date{} 
\maketitle
\renewcommand{\baselinestretch}{2}
\begin{quotation}
\begin{spacing}{1.5}
\noindent
\small Abstract:~
In this paper, 
we prove   that functional  sliced inverse regression (FSIR) achieves the optimal (minimax) rate for estimating the central space in functional sufficient dimension reduction problems.
First, we  provide a concentration inequality for the FSIR estimator  of the covariance of the conditional mean. Based on this inequality, we establish
the root-$n$  consistency of  the FSIR estimator of the image of covariance of the conditional mean.
Second,  we  apply the most widely used truncated scheme to estimate   the  inverse of the covariance operator and identify    the truncation parameter that ensures 
 that FSIR can achieve the  optimal minimax convergence rate for  estimating the central space. 
Finally, we conduct simulations to demonstrate the optimal choice of truncation parameter and the estimation efficiency of FSIR.
To the best of our knowledge, this is the first paper to  rigorously  prove the minimax optimality of FSIR in estimating the central space for multiple-index models and general $Y$ (not necessarily discrete).  
\end{spacing}
\vspace{9pt}

\par
\end{quotation}\par

	\def\thefigure{\arabic{figure}}
	\def\thetable{\arabic{table}}
	
	\renewcommand{\theequation}{\thesection.\arabic{equation}}

	\fontsize{12}{14pt plus.8pt minus .6pt}\selectfont
\newpage
	
	
	\renewcommand{\baselinestretch}{2}
	
	\markright{ \hbox{\footnotesize\rm Statistica Sinica
		}\hfill\\
		\hbox{\footnotesize\rm
		}\hfill }
	
	\markboth{\hfill{\footnotesize\rm FIRSTNAME1 LASTNAME1 AND FIRSTNAME2 LASTNAME2} \hfill}
	{\hfill {\footnotesize\rm FILL IN A SHORT RUNNING TITLE} \hfill}
	
	\renewcommand{\thefootnote}{}
	$\ $\par

\section{Introduction}
Sufficient Dimension Reduction (SDR) aims to identify a low-dimensional subspace that captures the most important features of the data that are relevant to the response variable. It is a useful tool for researchers to perform both exploratory  data analysis and detailed model developments when the dimension of the covariates is high.
Concretely, for a pair of random variables $(\boldsymbol{X},Y)\in\R^p\times\R$, 
 an \textit{effective dimension reduction} (EDR) subspace is a  subspace $\mathcal{S}\subset \R^p$ such that $Y$ is independent of $\bs X$ given $P_{\S}\boldsymbol{X}$ (where $P_{\S}$ denotes the projection operator from $\R^p$ to $\S$), which can be represented as: 
  \begin{equation}\label{eqn:sdr}
     Y\perp\hspace{-2mm}\perp \boldsymbol{X} \mid P_{\S}\boldsymbol{X}.
 \end{equation}
 SDR targets at estimating the intersection of all EDR subspaces, which is shown to be again an EDR subspace  under mild conditions \citep{cook1996graphics}. This intersection is  often referred to as the  \textit{central space} and denoted by $\mathcal{S}_{Y\mid \boldsymbol{X}}$.
 To find the central space $\mathcal{S}_{Y\mid \boldsymbol{X}}$, researchers have developed a variety of methods: \textit{sliced inverse regression} (SIR, \citealt{li1991sliced}), \textit{sliced average variance estimation} (SAVE, \citealt{Cook1991}), \textit{principal hessian directions} (PHD, \citealt{KerChauLi1992}), \textit{minimum average variance estimation} (MAVE,  \citealt{Xia2009}), \textit{directional regression} (DR, \citealt{LiBing2007}), etc. 
SIR is  one of the most popular SDR methods and its  asymptotic properties  are of particular interest.
For more details, readers can refer to \cite{Hsing1992,Zhu1995,zhu2006sliced,wu2011asymptotic,linzl2018,lin2021optimality,tan2020sparse,huang2023sliced}.

There has been a growing interest in statistical modeling of functional data (i.e., the predictors are random functions in some function space $\mc H$), and researchers have extended existing multivariate SDR algorithms to accommodate this type of data. 
In a functional SDR algorithm,
the space $\R^p$ is replaced by the function  space $\mathcal{H}$,  
 a functional EDR subspace can be defined as a subspace $\S\subset\mathcal{H}$ such that \eqref{eqn:sdr} holds, and the intersection of all functional EDR subspaces can be referred to as the functional central space.

\textit{Functional sliced inverse regression} (FSIR) proposed by \cite{Ferre2003fsir} is one of the earliest  functional SDR algorithms. 
They  estimate the central space $\mathcal{S}_{Y\mid \boldsymbol{X}}$ under a  \textit{multiple-index model}. 
The model is mathematically represented as 
 \begin{equation}\label{eq:lmifsir}
     Y=f(\langle \bbeta_1,\boldsymbol{X}\rangle ,\dots, \langle \bbeta_d,\boldsymbol{X}\rangle ,\varepsilon)
 \end{equation}
 where $f:\R^{d+1}\to \R$ is an unknown (non-parametric) link function,  the predictor $\boldsymbol{X}$ and indices $\bbeta_{i}$'s are  functions in $L^{2}[0,1]$, the separable Hilbert space of square-integrable curves on $[0,1]$, and ${\varepsilon}$ is a random noise independent of $\boldsymbol{X}$. 
Although the individual $\bbeta_{j}$'s are unidentifiable because of 
the flexibility of the link  function $f$, 
 space $\mathcal{S}_{Y\mid \boldsymbol{X}}:=\mbox{span}\{\bbeta_{1},\cdots,\bbeta_{d}\}$ is estimable. 
Several key findings concerning FSIR have been established since its introduction. 
\cite{Ferre2003fsir}  showed the consistency of the FSIR estimator  under some technical assumptions. 
However, they did not establish a similar convergence rate for FSIR as those for the multivariate SIR obtained by \cite{Hsing1992} and \cite{Zhu1995}. 
\cite{Forzanicook2007note} found that the $\sqrt{n}$-consistency of  the central space  can not be achieved by the FSIR estimator  unless some very restrictive conditions on the covariance operator of the predictor are imposed. 
\cite{yao2015effective} introduced  the technique of functional cumulative slicing estimation (FCSE) for SDR, focusing on scenarios with sparse designs, and also obtained its  convergence rate.
  \cite{lianfsdr2015} showed that the convergence rate of the FSIR estimator for a discrete response $Y$ is the same as that for  \textit{functional linear regression} in \citep{Hall2007mcflr}, but they did not provide 
  a rigorous proof of the optimality of FSIR, which is deemed quite challenging. 
Recent developments in this field can be found in \cite{lian2014Sefsdr,Wanglain2020,tian2023functional}.   

Recently, there have been significant advances in understanding the behavior of SIR
for high dimensional data (i.e., $\rho:=\lim p/n$ is a constant or $\infty$). \cite{linzl2018} 
  established the  phase transition phenomenon of SIR in high dimensions, i.e., SIR can estimate the central space consistently   if and only if $\rho=0$.
  \cite{lin2021optimality} further obtained the minimax convergence rate of estimating the central space in high dimensions and showed that SIR can achieve the minimax rate. In a different setting, \cite{tan2020sparse}  studied  the minimax rate in high dimensions   under various loss functions and proposed a computationally tractable adaptive estimation scheme for sparse SIR.
\cite{huang2023sliced} generalized the  minimax rate results to cases with a large structural dimension $d$ (i.e., there is no constant upper bound on  the dimension of the central space $d$).
These work in high dimensional SIR  inspires our research to address the theoretical challenges of FSIR and bridge the aforementioned theoretical gap between FSIR and multivariate SIR.

In the present paper, we focus on examining the  error bound of FSIR under very mild conditions. We show that FSIR can estimate the central space optimally for general $Y$ (not necessarily discrete) over a large class of distributions.

\subsection{Major contributions}
Our main results are summarized as follows: 
\begin{itemize}
    \item[\textbf{(i)}] To study the asymptotic properties of FSIR under  very general settings, we introduce a fairly mild condition called the  \textit{weak sliced stable condition } (WSSC) for functional data (see Definition~\ref{def:weak slice stable}). 
    \item[\textbf{(ii)} ]
    Under the above WSSC, we prove a concentration inequality for the FSIR estimator $\widehat{\Gamma}_{e}$ around its population counterpart  $\Gamma_e:=\var(\E[\boldsymbol{X}\mid Y])$ (see Lemma \ref{lem:cvgem}). 
    \item[\textbf{(iii)}]  Based on the  concentration inequality, we  show that the space spanned by the top $d$ eigenfunctions of $\widehat{\Gamma}_{e}$    is a root-$n$ consistent estimator of the image of $\Gamma_{e}$ (see Theorem \ref{thm:rootnconsistency}). This part is a crucial step to  our key results, the minimax rate optimality of FSIR estimator for the central space. 
    \item[\textbf{(iv)}] Having \textbf{(i)}-\textbf{(iii)} established, we apply the most widely used truncated scheme to estimate   the  inverse of the covariance operator of the predictor and then  establish the consistency of the FSIR estimator for the central space. Furthermore, we identify the optimal truncation parameter to achieve the minimax optimal convergence rate for FISR in Theorem \ref{thm:cfsdrs}. It turns out that the converging rate of FSIR is the same as the minimax rate for estimating the slope in functional linear regression \citep{Hall2007mcflr}. 
    \item[\textbf{(v)} ] Finally, we show that the convergence rate we obtained in \textbf{(iv)} is minimax  rate-optimal for multiple-index models over a large class of distributions (see Theorem \ref{thm:opmmultiple}). 
    \item[\textbf{(vi)} ] Simulation studies  show that the optimal choice of $m$ matches the theoretical ones and   illustrate the efficiency of FSIR on both synthetic and real data.

\end{itemize}

To the best of our knowledge,   this is the first work that rigorously establishes the optimality of FSIR  for a general response $Y$. Our results  provide a precise characterization of the difficulty associated with the estimation of the functional central space in terms of the minimax rates over  a wide range of distributions.  
It not only enriches the existing theoretical results of FSIR, but also opens up new possibilities for extending other well-understood results derived from high-dimensional data to those related to functional data.

\bigskip

\subsection{Notations and organization}\label{sec:pan} 
Throughout the paper, we take $\H=L^{2}[0,1]$ to be the separable Hilbert space of square-integrable curves on $[0,1]$ with the inner product $\langle f,g\rangle=\displaystyle{\int_{0}^1}f(u)g(u)\,\mathrm{d}u$ and norm $\|f\|:=\sqrt{\langle f,f\rangle}$ for $f,g\in\mathcal H$.

For an operator $T$ on $\H$, $\|T\|$ denotes its operator norm with respect to $\langle \cdot,\cdot\rangle $, i.e., 
\begin{align*}
\|T\|:={\sup_{\bs{u}\in\mathbb{S}_\mathcal{H}}}\|T(\bs{u})\|
\end{align*}
where $\mathbb{S}_{\mathcal{H}}=\{\bs{u}\in\mathcal{H}:\|\bs{u}\|=1\}$. $\mbox{Im}(T)$    denotes the closure of the image of $T$, $P_{T}$   the projection operator from $\H$ to $\mbox{Im}(T)$,     and $T^*$ the adjoint operator of $T$ (a bounded linear operator).  If $T$ is self-adjoint,  $\lambda^{+}_{\min}(T)$ denotes  the  infimum of the positive spectrum of $T$ and $T^{\dagger}$   the Moore--Penrose pseudo-inverse  of $T$.
Abusing notations, we also denote by $P_S$ the projection operator onto a closed space $S\subseteq \mathcal{H}$.  For any $x,y\in\H$, $x\otimes y$ is the operator of $\H$ to itself, defined by $x\otimes y(z)=\langle x,z\rangle y, \forall z\in\H$. 
For any random element $\boldsymbol{X}=\bs X_t\in \H$, its mean function is defined as $(\mb E\bs X)_t=\mb E[\bs X_t]$.  For any random operator $T$ on $\H$, the mean $\E[T]$ is defined as the unique operator on $\H$ such that for all $z\in \H$, 
$ (\E[T])(z)=\E[T(z)]$. Specifically,
the covariance operator of $\boldsymbol{X}$, $\var(\boldsymbol{X})$, is defined as $\var(\boldsymbol{X})(z)=\E\left( \langle \boldsymbol{X}, z \rangle \boldsymbol{X} \right) - \langle\E\boldsymbol{X}, z \rangle \E\boldsymbol{X}$. 
For a pair of random variables $(\boldsymbol{X},Y)\in \H\times \R$, $\Gamma$ and $\Gamma_{e}$ denote the covariance  operator of $\boldsymbol{X}$ and $\E[\boldsymbol{X}\mid Y]$ respectively,  i.e.,  
\begin{align}\label{eq:def_Gamma}
    \Gamma:=\var(\boldsymbol{X})  \qquad\mbox{and } \quad \Gamma_{e}:=\var(\E[\boldsymbol{X}\mid Y]).
\end{align}

For two sequences $a_n$ and $b_n$, we denote $a_n\lesssim  b_n$ (resp. $a_n\gtrsim b_n$ ) if there
exists  a  positive constant $C$ such that $a_n\leqslant Cb_n$ (resp. $a_n\geqslant C b_n$), respectively. We denote
$a_n\asymp b_n$ if both $a_n\lesssim  b_n$ and $a_n\gtrsim b_n$ hold. 
For a random sequence $X_n$, we denote by $X_n=O_{p}(a_n)$ that $\forall\varepsilon>0$, there exists a constant $C_\varepsilon>0$, such that
$\sup_n\mathbb{P}(|X_n|\geqslant C_\varepsilon a_n)\leqslant\varepsilon$. 
Let $[k]$ denote $\{1,2,\dots,k\}$ for some positive integer $k\geqslant1$.

 The rest of this paper is organized as follows. We first  provide a brief review of FSIR in Section \ref{sec:prefsir}.  After introducing the weak sliced stable condition for functional data in Section \ref{sec:Weak sliced stable condition for functional data}, we establish the  root-$n$ consistency of the estimated inverse regression subspace in Section \ref{sec:Consistency inverse regression subspace}. Lastly, the minimax rate optimality of FSIR are shown in  Section \ref{sec:Optimality of  FSIR} and the numerical experiments are  reported in Section \ref{sec:simulationfsir}. All proofs are deferred to the supplementary files.

\section{Optimal Truncated FSIR}\label{sec:prefsir}
Without loss of generality, we assume that $\vX\in\mc H$ satisfies $\mathbb{E}[\boldsymbol{X}]=0$ throughout the paper. As is usually done in functional data analysis \citep{Ferre2003fsir,lian2014Sefsdr,lianfsdr2015}, we assume  that $\mathbb{E}[\|\boldsymbol{X}\|^4]<\infty$, 
which implies
that $\Gamma$ is a trace class \citep{hsing2015theoretical} and $\boldsymbol{X}$ possesses the following Karhunen--Lo\'{e}ve expansion: 
 \begin{align}\label{eq:KL}
\boldsymbol{X}=\sum^{\infty}_{i=1}\xi_{i}\phi_{i}
 \end{align}
where $\xi_{i}$'s  are random variables satisfying   $\E[\xi^{2}_{i}]=\lambda_{i}$  and $\E[\xi_{i}\xi_{j}]=0$ for $i\neq j$ and $\{\phi_{i}\}^{\infty}_{i=1}$ are the eigenfunctions of  $\Gamma$ in \eqref{eq:def_Gamma} associated with  the decreasing eigenvalues sequence  $\{\lambda_{i}\}^{\infty}_{i=1}$. 
In addition, we assume that $\Gamma$ is non-singular (i.e., $\lambda_i>0, \forall i$) as the literature on functional data analysis usually does. Since $\Gamma$ is compact ($\Gamma$ is a trace class), by spectral decomposition theorem of compact operators, we know that $\{\phi_{i}\}^{\infty}_{i=1}$ forms a complete basis of $\mc H$.

In order for FSIR to produce a consistent  estimator of the functional central space $\mathcal{S}_{Y\mid \boldsymbol{X}}$ for $(\vX,Y)$ from the multiple index model \eqref{eq:lmifsir}, people often assume that the joint distribution of $(\boldsymbol{X},Y)$ satisfies the following conditions (see e.g., \cite{Ferre2003fsir,lian2014Sefsdr,lianfsdr2015}).
\begin{assumption}\label{as:Linearity condition and Coverage condition}
The joint distribution of $(\bs X,Y)$ satisfies
\begin{itemize}
 \item[\textbf{i)}] \textit{Linearity condition:}
For any $\bs b\in \mathcal{H}$, $\E[\langle \bs b,\boldsymbol{X}\rangle \mid (\langle \bbeta_1,\boldsymbol{X}\rangle ,\dots, \langle \bbeta_d,\boldsymbol{X}\rangle )]$ is linear in $\langle \bbeta_1,\boldsymbol{X}\rangle ,\dots, \langle \bbeta_d,\boldsymbol{X}\rangle $.
 \item[\textbf{ii)}]\textit{Coverage condition:} $\mathrm{Rank}\left(\mathrm{var}(\E[\bs X|Y])\right)=d$.
\end{itemize}
\end{assumption}
Both of these conditions are  natural generalizations of the multivariate ones that appear  in the multivariate  SIR literature \citep{li1991sliced,hall1993almost,LiHsing2010dec}.
They are necessary for \cite{Ferre2003fsir} to  establish that  the \textit{inverse regression subspace}  
$\S_{e}:=\mr{span}\{\E[\boldsymbol{X}\mid Y=y]\mid y\in\R \}$ equals to the space $ \Gamma \mathcal{S}_{Y\mid \boldsymbol{X}}:=\mbox{span}\{\Gamma\bbeta_{1}, \dots,\Gamma\bbeta_d\}$.
Since $\S_{e}=\mr{Im}(\Gamma_e)$, FSIR  estimates $\S_{Y\mid \boldsymbol{X}}$ by estimating $\Gamma^{-1}\mr{Im}(\Gamma_e)$.

The FSIR  procedure for estimating  $\Gamma_e$ can be briefly summarized as follows. 
Given $n$ i.i.d. samples $\{(\boldsymbol{X}_{i},Y_{i})\}_{i=1}^n$ from the multiple index model \eqref{eq:lmifsir}, FSIR  sorts the samples  according to the order statistics $Y_{(i)}$ and then divide the samples into $H(\geqslant d)$ equal-size slices (for the simplicity of notation, we assume that $n=Hc$ for some positive integer $c$). We re-index the data as 
\begin{equation*}
	Y_{h,j}=Y_{(c(h-1)+j)}  \qquad\mbox{and} \qquad  \boldsymbol{X}_{h,j}=\boldsymbol{X}_{(c(h-1)+j)}
\end{equation*} 
where $\boldsymbol{X}_{(k)}$ is the concomitant of $Y_{(k)}$ \citep{yang1977}. Let $\S_{h}$ be the $h$-th interval $(Y_{(h-1,c)},Y_{(h,c)}]$ for $h=2,\dots,H-1$,  $\S_1=\{y\mid y\leqslant Y_{(1,c)}\}$ and  $\S_{H}=\{y\mid y>Y_{(H-1,c)} \}$. 
 Consequently, $\mathfrak{S}_{H}(n):=\{\mathcal{S}_{h}, h=1,..,H\}$ is a partition of $\R$ and is   referred to as the \textit{sliced partition}.
FSIR  estimates the conditional covariance  $\Gamma_{e}$ via
\begin{equation}\label{eq:def of hat Gamma e}
    \widehat{\Gamma}_{e}:=\frac{1}{H}\sum_{h=1}^{H}\overbar{\boldsymbol{X}}_{h,\cdot}\otimes \overbar{\boldsymbol{X}}_{h,\cdot}
\end{equation}
 where $\overbar{\boldsymbol{X}}_{h,\cdot}:=\frac{1}{c}\sum^{c}_{j=1}\boldsymbol{X}_{h,j}$ is the sample mean in the $h$-th slice.

To estimate $\Gamma^{-1}$, one need to resort to some truncation scheme. Given $n$ i.i.d. samples  $\{(\boldsymbol{X}_{i},Y_{i})\}_{i=1}^n$, a straightforward estimator of 
 $\Gamma^{-1}$ is $\widehat{\Gamma}^{\dagger}$, the pseudo-inverse of the sample covariance operator $\widehat{\Gamma}:=\frac{1}{n}\sum^{n}_{i=1}\boldsymbol{X}_{i}\otimes \boldsymbol{X}_{i}$. However, it is not practical since the operator $\Gamma$ is  compact ($\Gamma$ is a trace class) and then $\Gamma^{-1}$ is  unbounded. To circumvent this technical difficulty, one may apply some truncation strategies such as the operations in 
\cite{Ferre2003fsir}, which we briefly review as follows. We choose an integer $m$
and define the truncated covariance operator $\Gamma_{m}:=\Pi_{m}\Gamma\Pi_{m}$ where $\Pi_{m}:=\sum^{m}_{i=1}\phi_{i}\otimes\phi_{i}$ is the truncation  projection operator. Since   each  $\Gamma_{m}$ is of  finite rank,  we are able  to estimate $\Gamma_{m}^{\dagger}$. Specifically,  let the sample truncation operator $\widehat{\Pi}_{m}:=\sum_{i=1}^{m}\widehat{\phi}_{i}\otimes \widehat{\phi}_{i}$ and the sample  truncated covariance operator $\widehat{\Gamma}_{m}:=\widehat{\Pi}_{m}\widehat{\Gamma}\widehat{\Pi}_{m}$, where $\{\widehat{\phi}_{m}\}_{i=1}^m$ are the  top $m$   eigenfunctions  of $\widehat{\Gamma}$.  Then the estimator of $\Gamma_{m}^{\dagger}$ can be defined as $\widehat{\Gamma}^{\dagger}_{m}$. It is clear that $\|\Gamma-\Gamma_{m}\|\xrightarrow{m\rightarrow \infty} 0$ and the space $\Gamma_{m}\S_{Y\mid \boldsymbol{X}}$ would be close to the space $\Gamma\S_{Y\mid \boldsymbol{X}}$  when $m$ is sufficiently large. Thus we can accurately estimate $\Gamma^{-1}$ by $\wh\Gamma_m^\dag$ for sufficiently large $m$.


  

\begin{algorithm}[H]
\setstretch{1.5}
\begin{algorithmic}
\caption {FSIR \citep{Ferre2003fsir}.}\label{alg:FSIRm}
\State 
\begin{enumerate}
\item Standardize $\{\boldsymbol{X}_{i}\}_{i=1}^n$, i.e.,    
   $\bs Z_i:=\boldsymbol{X}_{i}-n^{- 1}\sum_{i=1}^n \bs X_{i}$; 
	\item Divide the $n$ samples $\{(\bs Z_i,Y_i)\}_{i=1}^n$ into $H$ equally sized slices according to the order statistics $Y_{(i)},1\leqslant i \leqslant n$; 
\item Calculate  $\widehat{\Gamma} =\frac{1}{n}\sum_{i=1}^{n}\bs Z_{i}\otimes \bs Z_{i}$ with its top $m$ eigenvalues $\{\widehat{\lambda}_i,1\leqslant i\leqslant m\}$ and the corresponding eigenfunctions $\{\widehat{\phi_{i}}:1\leqslant i \leqslant m\}$, where 
$m$ is the tuning parameter. Let $\widehat{\Gamma}_{m} =\sum^{m}_{i=1}\widehat{\lambda}_{i}\widehat{\phi}_{i}\otimes\widehat{\phi}_{i}$;

\item Calculate 
$\overbar{\bs Z}_{h,\cdot}=\frac{1}{c}\sum_{j=1}^{c}\bs Z_{h,j}$, $h=1,2,\cdots,H$ and 
  $\widehat{\Gamma}_{e}=\frac{1}{H}\sum_{h=1}^H\overbar{\bs Z}_{h,\cdot}\otimes \overbar{\bs Z}_{h,\cdot}$ similarly as \eqref{eq:def of hat Gamma e}
  ;
  
  \item  Find  the top $d$ eigenfunctions of $ \widehat{\Gamma}_{e}$, denoted by $\wh\bbeta_{k}' 
 ~(k=1,\dots,d)$
 and calculate   $\widehat{\bbeta}_{k}=\widehat{\Gamma}^{\dagger}_{m}\wh\bbeta_{k}'$.
\end{enumerate}
\State  Return $\mathrm{span}\left\{\widehat{\bs\beta}_1,...,\widehat{\bs\beta}_d\right\}$.
\end{algorithmic}
\end{algorithm}

In introducing our {\it optimal truncated FSIR algorithm} (FSIR-OT), we commence by revisiting the classical FSIR algorithm proposed in \cite{Ferre2003fsir} as shown in Algorithm\ref{alg:FSIRm}. It is worth noting that \cite{Ferre2003fsir} did not provide any specific guidance on the choice of the tuning parameter $m$. Our FSIR-OT algorithm provides an optimal selection criterion for $m$, namely, $m\propto n^{1/(\alpha+2\beta)}$, where $\alpha$ and $\beta$ are defined in Assumption \ref{assumption: rate-type condition}.
Under this optimal choice, we  prove that FSIR-OT can achieve
the minimax rate for estimating the central space in the next section.

\section{Minimax rate optimality of FSIR-OT}\label{sec:ssabfsir}
Throughout the paper,  the number of indexes $d$ is assumed known and fixed.
By analyzing  the asymptotic behaviors of FSIR-OT, we derive  the minimax rate optimality of FSIR-OT.
We begin by proposing  a fairly mild condition called  the weak sliced stable condition (WSSC) for functional data. Then, we show that the top $d$ eigenfunctions of the estimated conditional covariance  $\widehat{\Gamma}_{e}$ span a consistent estimator $\widehat{\S}_{e}$ of the inverse regression subspace $\S_{e}$, with convergence rate of $n^{-1/2}$ based on WSSC. 
Lastly, we establish the consistency of the FSIR-OT estimator of the central space  and show that the convergence rate is minimax optimal over a large class of distributions.

\subsection{Weak sliced stable condition for functional data}\label{sec:Weak sliced stable condition for functional data}
The sliced stable condition (SSC) was first introduced in \cite{linzl2018} to analyze the asymptotic behavior of SIR in high dimensions  such as the phase transition phenomenon. \cite{lin2021optimality} showed  the   optimality of SIR in high dimensions based on SSC. \cite{huang2023sliced} weakened SSC to weak sliced stable condition (WSSC) to establish the optimality of SIR in more general settings. This inspires us to extend WSSC  to functional data. Throughout the paper, $\gamma$ is a fixed small positive constant.

\begin{definition}[\textbf{Weak Sliced Stable Condition}]\label{def:weak slice stable}
Let $Y\in\R$ be a random variable, $K$ a positive integer and $\tau >1$ a constant. 
A partition $\mathcal{B}_{H}:=\{-\infty=a_0<a_1<\dots<a_{H-1}<a_H=\infty\}$ of $\R$ is called a $\gamma$-partition if 
	\begin{align}\label{eq:gamma partition}
	\frac{1-\gamma}{H} \leqslant \mathbb{P}(a_h\leqslant Y\leqslant a_{h+1})\leqslant \frac{1+\gamma}{H}, \qquad \forall h=0,1,\ldots, H-1.
	\end{align}
A continuous curve $\boldsymbol{\kappa}(y): \R\to\mc H$ is said to be  \textit{ weak $( K, \tau)$-sliced stable} w.r.t. $Y$, if 
	for any $H\geqslant K$ and any $\gamma$-partition $\mathcal{B}_{H}$, 
	it holds that
    	\begin{equation}\label{eq:essfk}
		\frac{1}{H}\sum^{H-1}_{h=0}\var\left(\langle\bs{u},\bs{\kappa}(Y)\rangle\mid a_h\leqslant Y\leqslant a_{h+1}\right)
		\leqslant \frac1\tau\var\left(\langle\bs{u},\bs{\kappa}(Y)\rangle\right)\quad(\forall\bs{u}\in\mathbb{S}_{\mathcal{H}}).
	\end{equation}
\end{definition}
Compared with the original SSC 
(e.g.,  the equation $(4)$ in \cite{linzl2018}) for the central curve $\bs m(y):=\bbE[\vX|Y=y]$, WSSC condition is less restrictive. 
The  average of the variances (the left hand side of \eqref{eq:essfk}) is only required to be sufficiently small by WSSC condition, in contrast to that, it needs to vanish as $H\to\infty$ by the original SSC. 
In fact, as we will in Theorem~\ref{thm:rootnconsistency}, the constant  $\tau$ in \eqref{eq:essfk} only needs to be greater than $\frac{6\norm{\Gamma_{e}}}{\lambda^{+}_{\min}(\Gamma_{e})}$ to guarantee the consistency of the FSIR-OT estimator.
Furthermore, the following lemma shows that WSSC of $\vm(y)$ is readily fulfilled under certain mild prerequisites. 
\begin{lemma}\label{lem:monment condition to sliced stable}
	Suppose that the  joint distribution of $(\vX,Y)\in \mc H\times\mb{R}$ satisfies the following  conditions:
\begin{itemize}
    \item[${\bf i)}$] for any  $\bs{u}\in \mb S_{\mc H}$,  $\mathbb{E}\left[|\langle\bs{u},\vX\rangle|^{\ell}\right]\leqslant c_{1}$ holds for absolute constants $\ell>2$ and $c_{1}>0$;
    \item[${\bf ii)}$]$Y$ is a continuous random variable;
    \item[${\bf iii)}$]the central curve $\vm(y):=\bbE[\vX|Y=y]$ is continuous.
\end{itemize}
Then for any $\tau>1$, there exists  an integer  $K=K(\tau,d)\geqslant d$ such that $\vm(y)$ is  weak  $(K,\tau)$-sliced stable w.r.t. $Y$.
\end{lemma}

\begin{assumption}\label{as: slice stable}
 The central curve $\bs{m}(y)=\E[\bs X|Y=y]$ is weak $(K,\tau)$-sliced stable with respect to $Y$ for two positive  constants $K$ and $\tau$ (i.e., WSSC).
\end{assumption}

 We note that the requirement of $K$ being a constant is mild since  $d$ is bounded.
With the help of WSSC, we can now bound the distance between $\widetilde{\Gamma}_{e}$ and $\Gamma_e$, where 
$$\widetilde{\Gamma}_{e}:=\frac{1}{H}\sum_{h:\mc S_h  \in \mathfrak{S}_{H}(n)}\overline{\bs{m}}_{h}\otimes \overline{\bs{m}}_{h} \ \text{ and }  \ \overline{\bs{m}}_{h} := \mathbb{E}[\bs{m}(Y) \mid Y\in \S_{h}]=\E[\bs X|Y\in\S_h].$$  Here $\mathfrak{S}_{H}(n)$ is the sliced partition defined in Section \ref{sec:prefsir}. This bound is key to obtaining a concentration inequality for the FSIR-OT estimator $\wh\Gamma_e$ of the conditional covariance  $\Gamma_e$. 
\begin{proposition}\label{prop:lgvgm}
Under Assumption \ref{as: slice stable}, there exist positive constants $C$ and $H_0\geqslant K$,  such that for all $H>H_{0}$,  if $n>1+4H/\gamma$ is sufficiently large,  we have 
\begin{equation}\label{eq:eslicedsta}
\mb{P}\left(\left|\left\langle\left( \widetilde{\Gamma}_{e}-\Gamma_{e}\right)(\bs{u}),\bs{u}  \right\rangle\right|
\leqslant \frac3\tau \left\langle\Gamma_{e}(\bs{u}),\bs{u}\right\rangle,\forall\bs{u}\in\mb{\mb{S}}_{\H}\right)\geqslant 1-CH^2\sqrt{n+1}\exp\left(\frac{-\gamma^2(n+1)}{32H^2}\right).
\end{equation}
\end{proposition}

When   $\tau>\frac{6\norm{\Gamma_{e}}}{\lambda^{+}_{\min}(\Gamma_{e})}$, we know that 
$\mbox{Im}(\widetilde{\Gamma}_{e})=\mbox{Im}(\Gamma_{e})$ holds with high probability (see Lemma 5 in Appendix for details).
 In other words,  Proposition \ref{prop:lgvgm} implies that $\mr{Im}(\widetilde{\Gamma}_{e})$ is a consistent estimator of $\mbox{Im}(\Gamma_{e})$ even if $\widetilde{\Gamma}_{e}$   is not    a consistent estimator of $\Gamma_{e}$.

\subsection{Root-$n$ consistency of the FSIR-OT estimator for inverse regression subspace}\label{sec:Consistency inverse regression subspace}
We study asymptotic behaviors of the  FSIR-OT estimator $\widehat{\S}_{e}$ of the inverse regression subspace $\S_{e}$. As in most studies in functional data analysis \citep{Hall2007mcflr, lei2014adaptive,lianfsdr2015,Wanglain2020},
we introduce the following assumption:
\begin{assumption}\label{as:moment}
There exists a constant $c_{1}>0$  such that  $\E[\xi_{i}^{4}]/\lambda^{2}_{i}\leqslant c_{1}$ uniformly  for all $i\in\mb Z_{+}$ where $\xi_i$ and $\lambda_i$ are defined in \eqref{eq:KL}.
\end{assumption}
Now we are ready to state our  first main result,  
which is similar to the `key lemma' in \cite{linzl2018},  a crucial tool for developing the phase transition phenomenon and establishing the minimax optimality of the high dimensional SIR.
\begin{lemma}\label{lem:cvgem}
Suppose that  Assumptions \ref{as: slice stable} and \ref{as:moment} hold.
For any fixed integer $H>H_{0}$ ($H_0$ is defined in Proposition \ref{prop:lgvgm}) and any sufficiently large  $n>1+4H/\gamma$,  we have 
\begin{align*}
\norm{ \widehat{\Gamma}_{e}-\Gamma_{e}}=O_p\left(\frac1\tau+ \sqrt{\frac{1}{n}}\right) \quad\mbox{ and }\quad
\norm{\widehat{\Gamma}_{e}-\widetilde{\Gamma}_{e}}=O_p\left(\sqrt{\frac{1}{n}}\right).
\end{align*}  
\end{lemma}
The  $\tau$ term  in the first equation of   Lemma \ref{lem:cvgem} suggests  that $\widehat{\Gamma}_{e}$ may not be a consistent estimator of $\Gamma_{e}$. However,    we are  interested in estimating the  space  $\mc S_e=\mbox{Im}(\Gamma_e)$ rather than $\Gamma_{e}$ itself and the $\tau$ term would not affect the convergence rate of $\norm{ P_{\widehat{\S}_{e}}-P_{\S_{e}}}$ as long as $\tau$ is sufficiently large.
This will be  elaborated  in  the following theorem, our second main result.
\begin{theorem}\label{thm:rootnconsistency}
 Consider the same conditions and constants as in  Lemma \ref{lem:cvgem} and suppose that $\tau>\frac{6\norm{\Gamma_{e}}}{\lambda^{+}_{\min}(\Gamma_{e})}$. It holds that 
\begin{equation}\label{eq:consisse}
   \mb E\left[\norm{ P_{\widehat{\S}_{e}}-P_{\S_{e}}}^2\right]\lesssim \frac{1}{n}
\end{equation}
where the expectation $\bbE$ is taken with respect to the randomness of the sample. 
\end{theorem}
Equation \eqref{eq:consisse}  implies that $\widehat{\S}_{e}$ is a root-$n$ consistent estimator  of the inverse regression subspace  $\S_e$. This  is a crucial step to  establish the minimax rate optimality of FSIR-OT estimator for the central space.

\subsection{Optimality of  FSIR-OT}\label{sec:Optimality of  FSIR}
 In order to obtain  the convergence rate  of the FSIR-OT estimator of the central space, we need a further assumption, which  is commonly imposed in functional data analysis (see e.g., \cite{Hall2007mcflr,lei2014adaptive,lianfsdr2015}). 
\begin{assumption}[Rate-type condition]\label{assumption: rate-type condition}
 There exist  positive constants $\alpha$, $\beta$, $c_{2}$ and $c_2'$ satisfying 
\begin{align*}
\alpha>1, \quad \frac{1}{2}\alpha+1<\beta, \quad  \lambda_{j}-\lambda_{j+1}\geqslant c_{2}j^{-\alpha-1} \mbox{ and }  |b_{ij}|\leqslant c_{2}'j^{-\beta}\quad(\forall i\in[d],j\in\mb Z_{+})
\end{align*}
where $b_{ij}:=\langle\eta_i,\phi_j\rangle$ for $\{\eta_i\}_{i=1}^d$ the generalized eigenfunctions of $\Gamma_e$ associated with top $d$ eigenvalues $\{\mu_i\}_{i=1}^d$ (i.e.,  $\Gamma_e\eta_i=\mu_i\Gamma\eta_i$).
\end{assumption}
The assumption on the eigenvalues $\lambda_{j}$ of $\Gamma$ requires a gap between  adjacent eigenvalues and ensures the accuracy of the 
 estimation of  eigenfunctions of $\Gamma$. It also implies a lower bound on the decay rate of $\lambda_j$: $\lambda_{j}\gtrsim  j^{-\alpha}$.
The assumption on the coefficients  $b_{ij}$ implies  that they do not decrease too slowly  with respect to $j$ 
uniformly for all $i$. It also implies that any basis $\{\wt\bbeta_i\}_{i=1}^d$ of $\S_{Y\mid \boldsymbol{X}}$ such that $\wt\bbeta_i=\sum_{j=1}^\infty\wt{b}_{ij}\phi_j$ satisfies $|\wt b_{ij}|\lesssim j^{-\beta}$.  
 The inequality $ \frac{1}{2}\alpha+1<\beta$   requires  that  the generalized eigenfunction $\eta_i$ is smoother than the covariate function $\vX$.

The conditions in Assumption~\ref{assumption: rate-type condition} have been imposed in \cite{Hall2007mcflr} for showing that the minimax rate of functional linear regression models is $n^{-(2\beta-1)/(\alpha+2\beta)}$. 
\cite{lianfsdr2015} also made use of some similar conditions 
to show  that the FSIR estimator of the  central space $\S_{Y\mid \boldsymbol{X}}$ for discrete $Y$ (i.e., $Y$ only takes finite values) can achieve  the same  convergence rate  as the one for estimating the slope in functional linear regression.

Now  we state  our  third main result, an upper bound on the  convergence rate of the FSIR-OT estimator of the central space. 
  \begin{theorem}\label{thm:cfsdrs}
  Suppose Assumptions \ref{as:Linearity condition and Coverage condition} to \ref{assumption: rate-type condition} hold with constants $\alpha$, $\beta$ and  $\tau>\frac{6\norm{\Gamma_{e}}}{\lambda^{+}_{\min}(\Gamma_{e})}$.   
By choosing   $m\asymp n^{\frac{1}{\alpha+2\beta}}$, we can get that for any fixed integer  $H>H_{0}$ ($H_0$ is defined in Proposition \ref{prop:lgvgm}) and any sufficiently large  $n>1+4H/\gamma$, we have  
\begin{align*}
\norm{P_{\widehat{\S}_{Y\mid \boldsymbol{X}}}-P_{\S_{Y\mid \boldsymbol{X}}}}^2=O_p\left(n^{\frac{-(2\beta-1)}{\alpha+2\beta}}\right)
\end{align*}
where $\widehat{\S}_{Y\mid \boldsymbol{X}}=\widehat{\Gamma}^{\dagger}_{m}\widehat{\S}_{e}$ is the estimated central space given by FSIR-OT.
\end{theorem}

The convergence rate we have derived for FSIR-OT   is  the same as the minimax rate for estimating the slope
in functional linear regression  \citep{Hall2007mcflr}. 
While the convergence rate appears to be the same as  that  in \cite{lianfsdr2015},  their study only considered  the case where the response $Y$ is discrete. Moreover, their work lacked a  proof for the optimality of FSIR-OT in estimating the central space.
\cite{yao2015effective} also introduced the FCSE method, focusing on scenarios with sparse designs, wherein only limited, noisy, and irregular observations are available for some or all subjects. However, they did not provide any analysis regarding the minimax optimality.
In the following, we will provide a rigorous proof that our convergence rate is indeed minimax rate-optimal over a large class of distributions, which is highly nontrivial.
To do this, we first introduce a class of distributions:
\begin{equation*}
\mathfrak{M}\left(\alpha,\beta,\tau,c_0,C_0\right):= \left\{
\left(\vX,Y\right)\quad\vline
\begin{aligned}
&Y=f(\langle \bbeta_1,\boldsymbol{X}\rangle ,\dots, \langle \bbeta_d,\boldsymbol{X}\rangle ,\varepsilon);\\
&\vX,\bbeta_i\in\mc H:=L^{2}[0,1]\quad(i=1,\dots,d);\\
&\varepsilon\text{ is a random noise independent of $\vX$};\\
& (\vX,Y) \text{ satisfies Assumption \ref{as:Linearity condition and Coverage condition}-\ref{assumption: rate-type condition} };\\
&c_0\leqslant\lambda_d(\Gamma_e)\leqslant\dots\leqslant\lambda_1(\Gamma_e)\leqslant C_0; \\
& \|\Gamma\|\leqslant C_0,\quad \lambda_{\min}(\Gamma |_{\mathcal{S}_{e}})\geqslant c_0 
\end{aligned} \right\}
\end{equation*}
where $c_0$ and $C_0$ are two positive universal constants.

 Then we have the following minimax lower bound for estimating the central space over $\mathfrak{M}\left(\alpha,\beta,\tau,c_0,C_0\right)$.
 \begin{theorem}\label{thm:opmmultiple}
For any given positive constants $\alpha$, $\beta$ and $\tau$ satisfying $\alpha>1,\frac{1}{2}\alpha+1<\beta$
 and $\tau>\frac{6C_0}{c_0}$, there exists an absolute constant $\vartheta>0$ that only depends on $\alpha$  and $\beta$,
  such that for any sufficiently large $n$, it holds that
    \begin{align*}
       \inf_{\wh{\mathcal{S}}_{Y\mid \boldsymbol{X}}}\sup_{\mc M\in \mathfrak{M}\left(\alpha,\beta,\tau,c_0,C_0\right)}\bbP_{\mc M}\left(\norm{P_{\wh{\mathcal{S}}_{Y\mid \boldsymbol{X}}}-P_{\mathcal{S}_{Y\mid \boldsymbol{X}}}}^{2}\geqslant \vartheta n^{-\frac{2\beta-1}{\alpha+2\beta}}\right)\geqslant0.9
    \end{align*}
where $\wh{\mathcal{S}}_{Y\mid \boldsymbol{X}}$ is taken over all possible estimators of $\mathcal{S}_{Y\mid \boldsymbol{X}}$ based on the training data $\{(\boldsymbol{X}_{i},Y_{i})\}^{n}_{i=1}$.
 \end{theorem}

 The main tool we used in proving this minimax lower bound is Fano's Lemma (see e.g., \citep{yu1997assouad}). 
The major challenge is to construct a specific family of distributions that are far apart from each other in the parameter space, yet close to each other in terms of Kullback--Leibler divergence. An important contribution in this paper is the construction of such distributions.

Theorems \ref{thm:cfsdrs} and \ref{thm:opmmultiple} together show that  the FSIR-OT estimator $P_{\widehat{\S}_{Y\mid \boldsymbol{X}}}$  is minimax rate-optimal for estimating the central space.

\section{Numerical Studies}\label{sec:simulationfsir}
In this section, we 
present several numerical experiments to illustrate the  behavior of the FSIR-OT algorithm.  The first experiment  demonstrates the optimal choice of the truncation parameter $m$ for
estimating the central space. 
The results corroborate the conclusion of Theorem \ref{thm:cfsdrs} that the choice in of FSIR-OT  (namely $m\asymp n^{\frac{1}{\alpha+2\beta}}$) is optimal. 
The second experiment focuses on the   estimation  performance  of FSIR-OT on synthetic data.  Lastly, we analyze a real data set on bike rentals using  FSIR algorithms.  
By comparing our algorithm with the FCSE algorithm of \cite{yao2015effective} and the regularized FSIR (\citealt{lianfsdr2015}, RFSIR), we demonstrate  advantages of FSIR-OT on both synthetic and real datasets. Similar to FSIR-OT, FCSE performs a truncation operation on the covariance operator, controlled by parameter \( m \), whereas RFSIR employs ridge-type regularization characterized by a regularization parameter \( \rho \).

\subsection{Generalized signal noise ratio (gSNR) of multiple index models}
Recall that the  signal-to-noise ratio (SNR) for the linear model $Y=\langle\bbeta,\bs{X}\rangle+\varepsilon$, where $\varepsilon\sim N(0,\sigma^2)$,
is defined as
\begin{align*}
\text{SNR}=\frac{\bbE[\langle\bbeta,\bs{X}\rangle^{2}]}{\bbE[Y^{2}]}=\frac{\langle\Gamma\bbeta,\bbeta\rangle}{\sigma^2+\langle\Gamma\bbeta,\bbeta\rangle}.
\end{align*}
A simple calculation shows that
\begin{equation*}
\Gamma_e=\frac{\Gamma\bbeta\otimes\Gamma\bbeta}{\langle\Gamma\bbeta,\bbeta\rangle+\sigma^2}, \mbox{ \ \  and  \ \ } \lambda(\Gamma_e)=\frac{\|\Gamma\bbeta\|^2}{\langle\Gamma\bbeta,\bbeta\rangle+\sigma^2},
\end{equation*}
where $\lambda(\Gamma_e)$ is the unique non-zero eigenvalue of $\Gamma_e$. This leads to the following identity for the linear model:
$$
\lambda (\Gamma_e)=\frac{\|\Gamma\bbeta\|^2}{\langle\Gamma\bbeta,\bbeta\rangle} \text{SNR}.
$$
Thus, in a  multiple index model we call $\lambda$, the smallest non-zero eigenvalue of $\Gamma_e$, the model's generalized SNR (gSNR) .

\subsection{Optimal choice of truncation parameter $m$}\label{sec:experiment optimal m}
Throughout this section, we set  $H=15$ and $\varepsilon\sim N(0,2)$. 
We note that  the results are not sensitive to the choice of $H$. The  guidelines for the choice of $H$ in practice are presented Section H.2 of  Supplementary Material.
The  experimental results  for other noise levels (with variances of $1$ and $0.25$, respectively) and other $H$ are shown in  Section H.3 of  Supplementary Material.

The following model is first considered:  

\begin{itemize}
    \item[(I)] $Y=\langle \bs{\beta}_1,\boldsymbol{X}\rangle +\varepsilon$,
    where $\boldsymbol{X}=\sum_{j=1}^{100}j^{-3/4}X_j\phi_{j}$   and $\bs{\beta}_1=\sum_{j\geqslant 1}(-1)^{j}j^{-2}\phi_{j}$.  Here  $X_j\overset{\mathrm{iid}}{\sim}N(0,1)$,  $\phi_{1}=1, \phi_{j+1}=\sqrt{2}\cos(j\pi t), j\geqslant 1$.
\end{itemize}

Note that the construction of $\bs X$ here is equivalent to a construction that satisfies the assumption that $\Gamma$ is non-singular (i.e., $\lambda_i>0, \forall i$). A detailed explanation is deferred to Section H.1 of  Supplementary Material.

For this model $\alpha=3/2$ and $\beta=2$, so the optimal choice of $m$ used by FSIR-OT satisfies $m\propto n^{2/11}$. The gSNRs of Model I 
are 0.791,  0.498, and 0.333, respectively, when the noise variances are 0.25, 1, and 2. 

 To evaluate the performance of  FSIR-OT, we consider the  subspace estimation error defined as $\mc D(\bs{\widehat B};\bs B):=\left\|P_{\bs{\widehat B}}-P_{\bs B}\right\|$ where $\bs{\widehat B}:=(\widehat\bbeta_1,\dots,\widehat\bbeta_d):\R^d\to L^2[0,1]$ and $\bs B:=(\bbeta_1,\dots,\bbeta_d):\R^d\to L^2[0,1]$.
This metric takes value in $[0,1]$ and, the smaller it is, the better the performance.
Each trial is repeated  $100$ times for reliability. 

\begin{figure}[H]
	\centering
	\begin{minipage}{0.50\textwidth}
		\includegraphics[width=\textwidth]{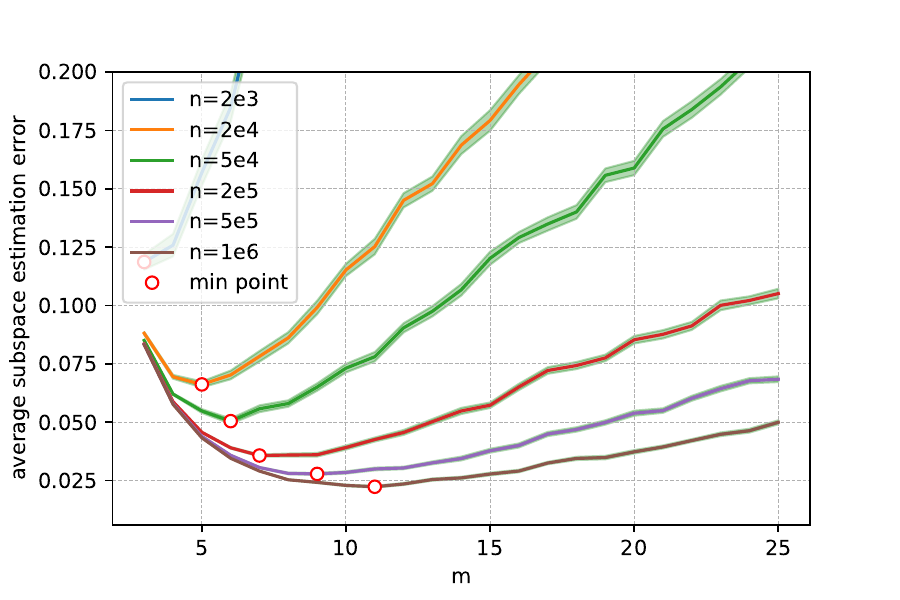}
	\end{minipage}\hfill
	\begin{minipage}{0.50\textwidth}
		\includegraphics[width=\textwidth]{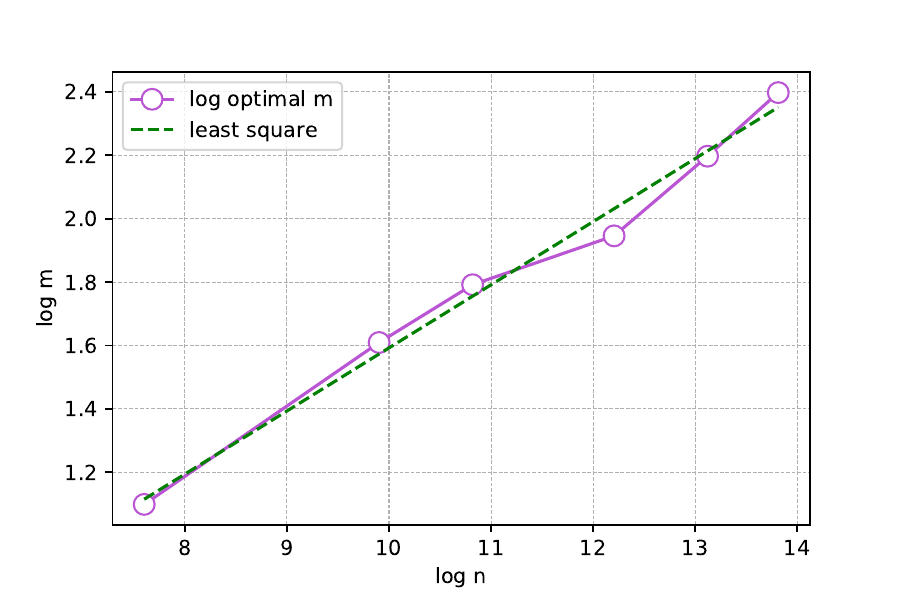}
	\end{minipage}
\caption{
Experiments for the optimal choice of truncation parameter $m$ with $\varepsilon\sim N(0,2)$ and $H=15$. 
Left: average subspace estimation error with increasing $m$ for different $n$. 
Right: linear trend of $\log(m^*)$ against $\log(n)$, with a slope of $0.2$  and $R^2>0.98$.}
\label{fig:optimal m noise2H15}
\end{figure}

The left panel of Figure \ref{fig:optimal m noise2H15} is the average subspace estimation error under Model (I) where $n$ ranges in $\{2\times 10^3,2\times 10^4,5\times10^4,2\times 10^5,5\times 10^5,10^6\}$, $m$ ranges in $\{3,4,\dots,25\}$. The optimal value of $m$ (denoted by $m^*$) for each $n$ is marked with a red circle. 
Among the 100 replicates for every $n$, the number of times that the minimal estimation error occurs at $m^*$ is $48,32,29,29,26,26$, respectively. 
The shaded areas represent the standard error bands associated with these estimates (all smaller than  $0.009$). 
The right panel of Figure~\ref{fig:optimal m noise2H15} illustrates the linear dependence of  $\log(m^*)$ on $\log(n)$. The solid line characterizes the linear trend of $\log(m^*)$ against $\log(n)$. The dotted line is their least-squares fitting, with its slope estimated as $0.2$, which is close to the theoretical value of $2/11$. These results are consistent with the theoretically optimal choice of $m$ in FSIR-OT.

\subsection{Subspace estimation error performance in  synthetic data}\label{sec: synthetic data}
In this section, we compare  FSIR-OT with RFSIR and FCSE  for model (I) from Section~\ref{sec:experiment optimal m}  and the following two models:
\begin{itemize}
\item[(II)] $Y=\langle \bs{\beta}_{1},\boldsymbol{X}\rangle +100\langle \bs{\beta}_{2},\boldsymbol{X}\rangle ^{3}+\varepsilon$, 
     where $\bs{\beta}_{1}(t)=\sqrt2\sin(\frac{3\pi t}{2})$,  $\bs{\beta}_{2}(t)=\sqrt2\sin(\frac{5\pi t}{2})$ for $t\in [0,1]$, and  $\boldsymbol{X}$ is the standard Brownian motion on $[0,1]$ (The Brownian motion is approximated by the top 100 eigenfunctions of the Karhunen--Lo\`{e}ve decomposition in practical implementation).
    \item[(III)] $Y=\exp(\langle \bs{\beta},\boldsymbol{X}\rangle )+\varepsilon$,
    where $\boldsymbol{X}$ is the standard Brownian motion on $[0,1]$, and $\bs{\beta}=\sqrt2\sin(\frac{3\pi t}{2})$.
\end{itemize} 
For model II and model III, we compute the estimated gSNR by $\lambda_{d}(\wh\Gamma_e)$, the $d$-th eigenvalue of the SIR estimate of $\Gamma_e$ based on  $2000$ replicates, where $n = 10000$. The mean gSNRs (standard deviation) of model II are 0.020 (0.001), 0.009 (0.001), and 0.003 (0.001), respectively, when the noise variances are 0.25, 1, and 2. The mean gSNRs (standard deviation) of model III are 0.729 (0.01), 0.536 (0.01), and 0.305 (0.01), respectively, when the noise variances are 0.25, 1, and 2.

\begin{figure}[H]
	\centering
	\begin{minipage}{0.33\textwidth}
		\includegraphics[width=\textwidth]{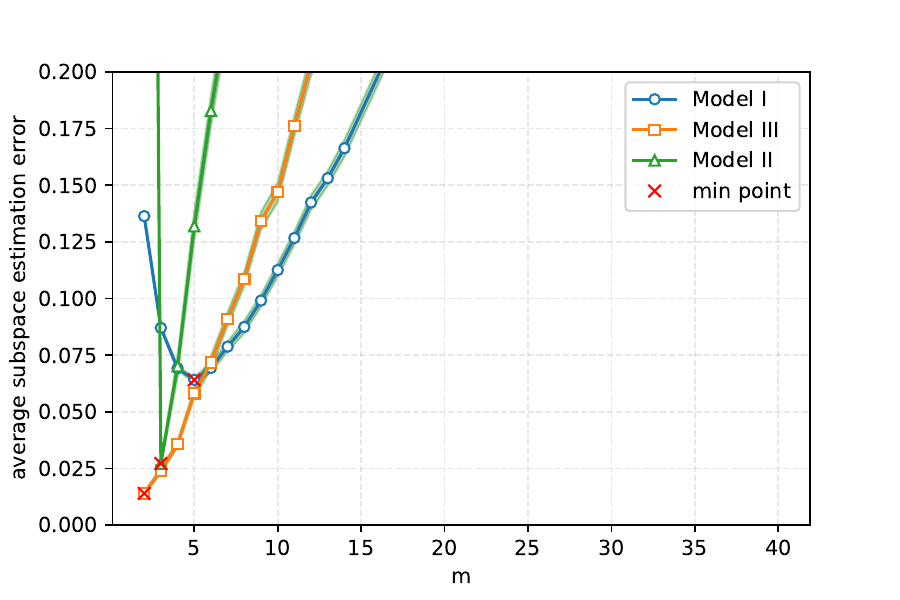}
	\end{minipage}\hfill
	\begin{minipage}{0.33\textwidth}
		\includegraphics[width=\textwidth]{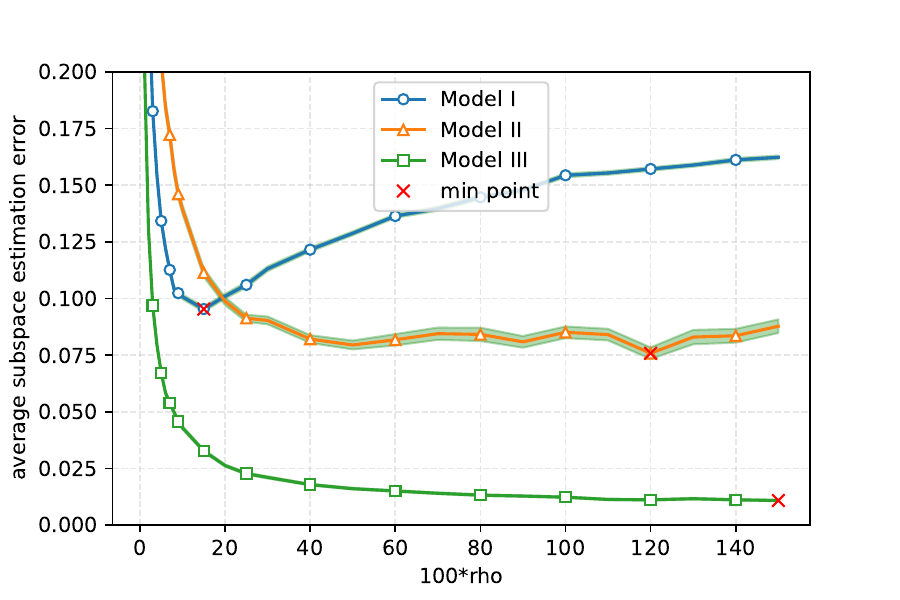}
	\end{minipage}
	\begin{minipage}{0.33\textwidth}
		\includegraphics[width=\textwidth]{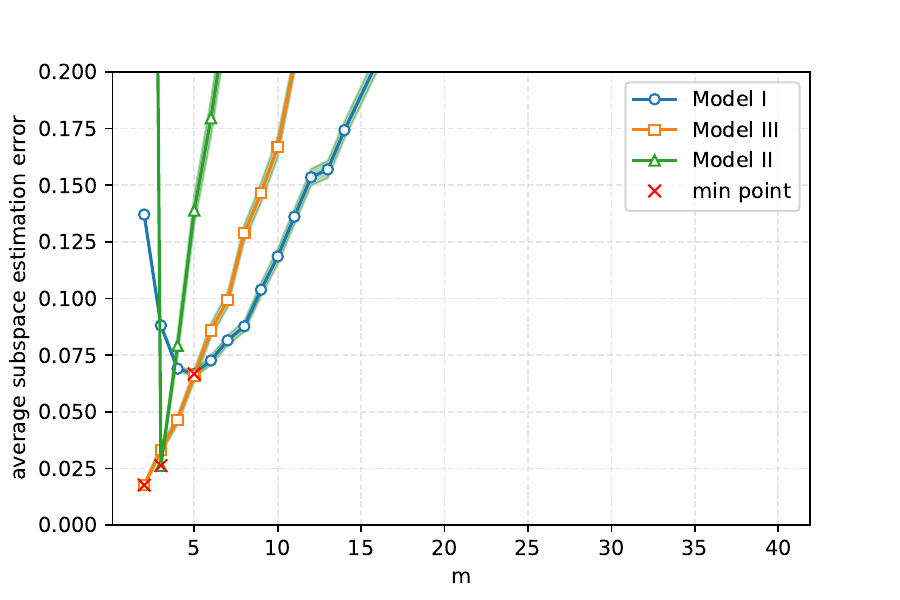}
	\end{minipage}
\caption{
Average subspace estimation error of FSIR-OT, RFSIR and FCSE for various models in the case of $\varepsilon\sim N(0,2)$ and $H=15$. 
The standard errors are all below $0.01$. 
Left: FSIR-OT with different truncation parameter $m$; 
Middle: RFSIR with different values of the regularization parameter $\rho$; Right: FCSE with different truncation parameter $m$.
}
\label{fig:error 3models noise2H15}
\end{figure}

For each model, we calculate the average subspace estimation error of FSIR-OT, RFSIR and FCSE 
based on $100$ replicates, where $n=20000$, the truncation parameter  $m$ of FSIR-OT and FCSE  ranges in $\{2,3,\dots,13,14,20,30,40\}$, and the regularization parameter $\rho$ in RFSIR  ranges in $0.01\times \{1,2,\cdots,9,10,15,20,25,30,40,\cdots,140,150\}$. Detailed results are presented in Figure \ref{fig:error 3models noise2H15}, 
where we mark the minimal error in each model with red `$\times$'
 and denote the corresponding value of truncation (or regularization) parameter by $m^*$ (or $\rho^*$). 
The shaded areas represent the corresponding standard errors, all of which are less than $0.012$. 
For FSIR-OT,  
the  minimal errors for $\mc M_1$, $\mc M_2$, and $\mc M_3$ are  $0.06,0.03$, and $0.01$ 
respectively. 
Among the 100 replicates for every model, the number of times that the minimal estimation error occurs at $m^*$ is $33,80$ and $70$, respectively. 
For RFSIR,  the corresponding minimal errors are 0.10, 0.08, and $0.01$, respectively. 
Among the 100 replicates for every model, the number of times that the minimal estimation error occurs at $\rho
^*$ is $37,13$ and $18$, respectively. For FCSE,  the corresponding minimal errors are $0.07,0.03$, and $0.02$. Among the 100 replicates for every model, the number of times that the minimal estimation error occurs at $m^*$ is $24,23$ and $28$, respectively.

The results here suggest that the performance of FSIR-OT is generally superior to, or at the very least equivalent to, that of  RFSIR and FCSE.

\subsection{Application to real data}\label{sec: real data}
\begin{figure}[H] 
	\centering
	\begin{minipage}{0.50\textwidth}
		\includegraphics[width=\textwidth]{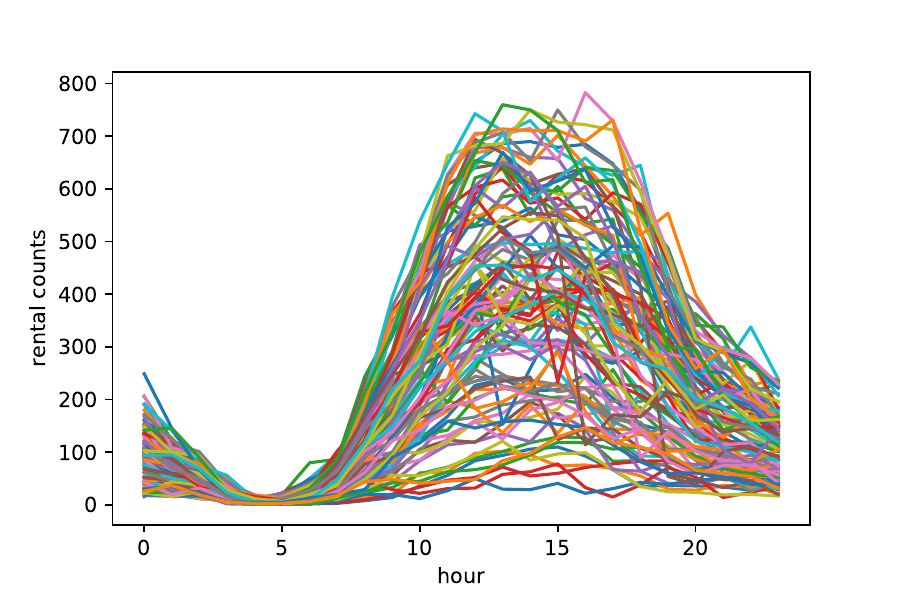}
	\end{minipage}\hfill
	\begin{minipage}{0.50\textwidth}
		\includegraphics[width=\textwidth]{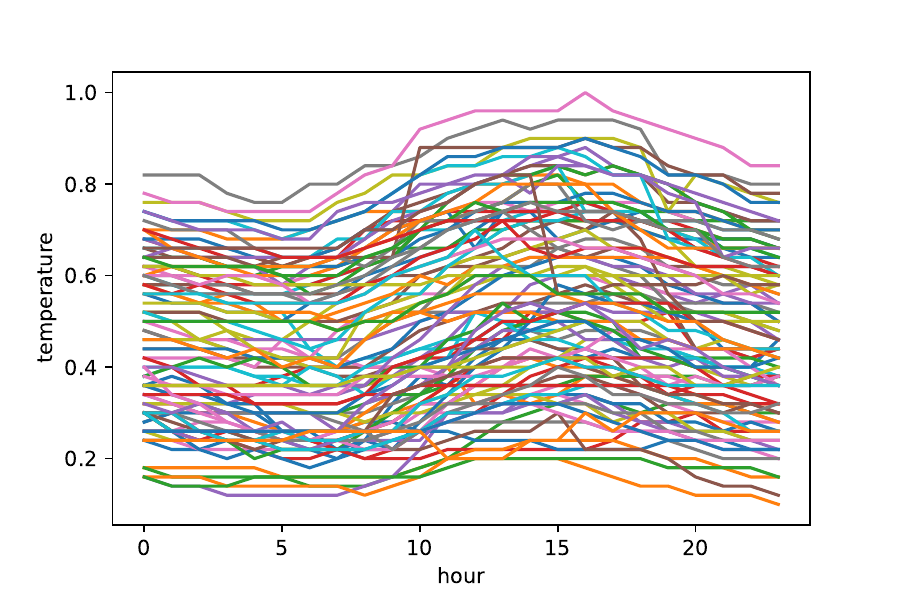}
	\end{minipage}
 \caption{Bike sharing data}
\label{The bike sharing data}
\end{figure}

In the following, we apply FSIR-OT to 
 a  business data analysis problem regarding bike sharing. The data are   available from
 \url{https://archive.ics.uci.edu/ml/datasets/Bike+Sharing+Dataset}.  The main purpose is to analyze  how the bike rental counts are affected by the temperature on Saturdays. After removing data from $3$ Saturdays  with missing information, we plot hourly bike rental  counts and hourly normalized temperature  (values  divided by the maximum 41\textdegree C)   on $102$  Saturdays in Figure \ref{The bike sharing data}. In the following experiments, we treat hourly normalized  temperature and the logarithm of daily average bike rental counts  as predictor function and scalar response respectively.

In order to compare the estimation error performance of FSIR-OT with 
 RFSIR and FCSE for estimating the central space, we employ  
dimension reduction using these  algorithms  with  $H=15$ as
an intermediate step in modelling the relation between the predictor and response. 
Specifically, given any training samples $\{(\vX_i,Y_i)\}_{i=1}^n$, we utilize each dimension reduction algorithm to obtain a set of low-dimensional predictors $\bs x_i$ for $i\in [n]$.
 Then, we employ Gaussian process regression to fit a nonparametric regression model  based on samples $\{(\bs x_i,Y_i)\}_{i=1}^n$. We randomly select $90$ samples as the training data and  then  calculate the out-of-sample mean squared error (MSE) using the remaining samples.

Since $d$ is unknown in most real applications, including this one, we follow the PCA approach by calculating the sum of the first $k$  $(k\leq 24)$  eigenvalues of $\wh\Gamma_e$. We found for this dataset that the first five eigenvalues account for $99.8\%$ of the summation of all eigenvalues. Therefore, we narrowed the selection range of $d$ to $\{1,2,3,4,5\}$. 
For each chosen  $d$, we then selected the corresponding value of $m$ satisfying $m\geq d$.
The experiment is repeated $100$ times and the mean and standard error are presented in Table~ \ref{tab:real data error}.  
Among the 100 replicates for every method (FSIR-OT, FCSE, and RFSIR), the number of times that the minimal estimation error occurs at the optimal pair of tuning parameters (i.e., $(m^*,d^*)$, $(m^*,d^*)$, and $(\rho^*,d^*)$, respectively)  is $6$, $4$, and $8$, respectively.

 From Table \ref{tab:real data error}, 
it can be concluded that FSIR-OT performs better than FCSE and RFSIR  if both methods are fine-tuned. 
Furthermore, the best result of FSIR-OT is observed at $d=2$, while those for FCSE and RFSIR are both observed at $d=4$.
This means that, if all methods were further fine-tuned, FSIR-OT would have provided a more accurate and simpler (lower dimensional) model for the relationship between the response variable and the predictor than the other two methods.

\begin{table}[H]
\begin{center}
\renewcommand\arraystretch{0.5}
\begin{tabular}{|c|c|c|c|c|c|c|}
\hline
$H=15$ & $m$   & $2$ & $4$ & $6$ & $8$ & $10$   \\
\hline
\hline
\multirow{5}{*}{FSIR-OT}
& $d=1$  & \textbf{0.209} & \textbf{0.212} & \textbf{0.207} & \textbf{0.211} & \textbf{0.206} \\
&  & \footnotesize{(0.010)} & \footnotesize{(0.010)} & \footnotesize{(0.010)} & \footnotesize{(0.009)} & \footnotesize{(0.009)}
\\
& $d=2$  & \textbf{0.229} & \textbf{0.196} & \textbf{\underline{0.188}} & \textbf{0.200} & \textbf{0.215} \\
&  & \footnotesize{(0.013)} & \footnotesize{(0.009)} & \footnotesize{(0.008)} & \footnotesize{(0.010)} & \footnotesize{(0.010)}
\\
& $d=3$  &  & \textbf{0.209} & \textbf{0.193} & \textbf{0.208} & \textbf{0.207} \\
&  && \footnotesize{(0.012)} & \footnotesize{(0.009)} & \footnotesize{(0.010)} & \footnotesize{(0.010)}
\\
& $d=4$  &  & \textbf{0.229} & \textbf{0.216} & \textbf{0.213} & \textbf{0.224} \\
&  &  & \footnotesize{(0.011)} & \footnotesize{(0.011)} & \footnotesize{(0.009)} & \footnotesize{(0.010)}
\\
& $d=5$  &  &  & \textbf{0.245} & \textbf{0.284} & \textbf{0.316} \\
&  & &  & \footnotesize{(0.014)} & \footnotesize{(0.021)} & \footnotesize{(0.025)}
\\[1pt]
\hline
\multirow{5}{*}{FCSE}
& $d=1$  & \textbf{0.207} & \textbf{0.206} & \textbf{0.190} & \textbf{0.214} & \textbf{0.230} \\
&  & \footnotesize{(0.010)} & \footnotesize{(0.010)} & \footnotesize{(0.008)} & \footnotesize{(0.010)} & \footnotesize{(0.012)}
\\
& $d=2$  & \textbf{0.215} & \textbf{0.222} & \textbf{0.202} & \textbf{0.197} & \textbf{0.195} \\
&  & \footnotesize{(0.010)} & \footnotesize{(0.009)} & \footnotesize{(0.009)} & \footnotesize{(0.010)} & \footnotesize{(0.010)}
\\
& $d=3$  &  & \textbf{0.216} & \textbf{0.209} & \textbf{0.214} & \textbf{0.207} \\
&  &  & \footnotesize{(0.011)} & \footnotesize{(0.010)} & \footnotesize{(0.010)} & \footnotesize{(0.011)}
\\
& $d=4$  & & \textbf{\underline{0.190}} & \textbf{0.223} & \textbf{0.220} & \textbf{0.207} \\
&  &  & \footnotesize{(0.007)} & \footnotesize{(0.010)} & \footnotesize{(0.012)} & \footnotesize{(0.010)}
\\
& $d=5$  &  &  & \textbf{0.254} & \textbf{0.255} & \textbf{0.302} \\
&  &  &  & \footnotesize{(0.012)} & \footnotesize{(0.015)} & \footnotesize{(0.039)}

\\[1pt]
\hline
 & $\rho$    & 0.044 & 0.101 & 0.159 & 0.216  & 0.274   \\
\hline
\multirow{5}{*}{RFSIR} 
& $d=1$  & \textbf{0.236} & \textbf{0.222} & \textbf{0.244} & \textbf{0.219} & \textbf{0.221} \\
&  & \footnotesize{(0.011)} & \footnotesize{(0.012)} & \footnotesize{(0.012)} & \footnotesize{(0.010)} & \footnotesize{(0.011)}
\\
& $d=2$  & \textbf{0.206} & \textbf{0.224} & \textbf{0.230} & \textbf{0.235} & \textbf{0.236} \\
&  & \footnotesize{(0.011)} & \footnotesize{(0.011)} & \footnotesize{(0.011)} & \footnotesize{(0.011)} & \footnotesize{(0.012)}
\\
& $d=3$  & \textbf{0.219} & \textbf{0.218} & \textbf{0.212} & \textbf{0.216} & \textbf{0.232} \\
&  & \footnotesize{(0.009)} & \footnotesize{(0.011)} & \footnotesize{(0.010)} & \footnotesize{(0.009)} & \footnotesize{(0.011)}
\\
& $d=4$  & \textbf{0.198} & \textbf{0.215} & \textbf{0.207} & \textbf{0.197} & \textbf{\underline{0.189}} \\
&  & \footnotesize{(0.010)} & \footnotesize{(0.010)} & \footnotesize{(0.010)} & \footnotesize{(0.008)} & \footnotesize{(0.008)}
\\
& $d=5$  & \textbf{0.208} & \textbf{0.193} & \textbf{0.211} & \textbf{0.211} & \textbf{0.234} \\
&  & \footnotesize{(0.010)} & \footnotesize{(0.009)} & \footnotesize{(0.012)} & \footnotesize{(0.011)} & \footnotesize{(0.011)}

\\[1pt]
\hline
\end{tabular}
\caption{
The mean (standard error) of the out-of-sample MSE for predicting logarithm of daily average bike rental counts using projected predictors after different dimension reduction methods.}\label{tab:real data error}
\end{center}
\end{table}

\begin{remark}
In dealing with real data, a crucial question is how to select the optimal 
$m$. The selection method provided in Theorem \ref{thm:cfsdrs} is based on asymptotic theory, which aims to provide minimax optimality results of FSIR under general conditions and is not directly applicable to real data. To date, the problem of selecting the optimal 
$m$ for a particular data set remains  unresolved, as shown in  \cite{Hall2007mcflr} and \cite{lianfsdr2015}.

To utilize  the asymptotic results of Theorem \ref{thm:cfsdrs} for selecting  $m$ in practice, we first estimate parameters $\alpha$ and $\beta$ according to Assumption \ref{assumption: rate-type condition}. 
Specifically, 
we first obtain  the $d$ eigenfunctions, $\widehat{\boldsymbol{\beta}}'_k$ ($k=1,\dots,d$), of $\widehat{\Gamma}_e$ and set $\wh\eta_k = \widehat{\Gamma}^{-1} \widehat{\boldsymbol{\beta}}'_k$. Then we estimate $\alpha$ and $\beta$ according to  $\wh\lambda_{j}-\wh\lambda_{j+1}\geqslant c_{2}j^{-\wh\alpha-1} \mbox{ and }  |\wh {b}_{ij}|:=\langle\wh\eta_i,\wh\phi_j\rangle\leqslant c_{2}'j^{-\wh\beta}$ where $\widehat{\Gamma}_e, \widehat{\Gamma}, \wh\lambda_j $ and $\wh\phi_j$ are defined in Algorithm \ref{alg:FSIRm}.
For example, we calculate $\wh\alpha=-(\frac{\ln(\wh\lambda_{j}-\wh\lambda_{j+1})}{\ln j}+1)$ and $\wh\beta=-\frac{\ln|\wh{b}_{ij}|}{\ln j}$ for sufficiently large $j$ respectively.
After we  get $\wh\alpha$ and $\wh\beta$, we choose $m$ in the interval
 $[n^{\frac{1}{\widehat{\alpha} + 2\widehat{\beta}}}/\log(n), n^{\frac{1}{\widehat{\alpha} + 2\widehat{\beta}}} \cdot \log(n)]$ and choose $\rho$ in $[n^{-\frac{\widehat{\alpha}}{\widehat{\alpha} + 2\widehat{\beta}}}/\log(n), n^{-\frac{\widehat{\alpha}}{\widehat{\alpha} + 2\widehat{\beta}}} \cdot \log(n)]$ (see \cite{lianfsdr2015}). 
 This approach significantly narrows down the choice range for $m$ and $\rho$ and is also consistent with our asymptotic results.
In our experiments, feasible values for $m$ were within 
$\{1,2,\dots,11\}$ and that for $\rho$ were within $[0.015,0.302]$. For the ease of presentation, we  selected $5$ representative values each for $m$ and $\rho$. Detailed results are presented in  Table \ref{tab:real data error}.
\end{remark}
\section{Discussion}

In this paper, we established the minimax rate-optimality of FSIR-OT for estimating the functional central space. 
Specifically,  we first prove an upper bound on the convergence rate of the FSIR-OT estimator of the functional central space under very mild assumptions. 
Then we establish a minimax lower bound on the estimation of the functional central space over a large class of distributions. These two results together show optimality of FSIR-OT. 
Our results not only enrich the theoretical understanding of FSIR-OT but also indicate the possibility of extending the findings of multivariate SDR methods to functional data.

There are some open questions related to the findings in this paper. 
First, the structural dimension $d$ is assumed to be bounded in the current paper. 
It is still unclear whether this restriction can be relaxed so that the minimax convergence rate of the functional central space estimation can be determined even when $d$ is large (i.e., there is no constant upper bound on $d$). 
Second, recent studies have revealed the dependence of the estimation error on the \textit{gSNR} defined as $\lambda_d(\Cov\left(\bbE[\vX \mid Y]\right))$ for multivariate SIR \citep{lin2021optimality,huang2023sliced}). 
Exploring the role of gSNR in the estimation of the functional central space will be an interesting next step.

\section*{Supplementary Materials}
Supplement to ``On the Optimality of Functional Sliced Inverse
Regression''. The supplementary material includes
the proofs for all the theoretical results in the paper.

\section*{Acknowledgements}
Lin's research was supported in part by the National Natural Science Foundation of China (Grant 92370122, 11971257). 
Huang's research was supported in part by NUS Start-up Grant A-0004824-00-0 and Singapore Ministry of Education AcRF Tier 1 Grant A-8000466-00-00. 
Liu's research was  supported in part by the National Science Foundation of the United States (DMS-2015411).

	\par

	\par
	

\bibliographystyle{chicago}      
\bibliography{SDR}   

\begin{thebibliography}{}

\bibitem[\protect\citeauthoryear{Angelova}{Angelova}{2012}]{angelova2012moments}
Angelova, J.~A. (2012).
\newblock On moments of sample mean and variance.
\newblock {\em Int. J. Pure Appl. Math\/}~{\em 79\/}(1), 67--85.

\bibitem[\protect\citeauthoryear{Cook}{Cook}{1996}]{cook1996graphics}
Cook, R.~D. (1996).
\newblock Graphics for regressions with a binary response.
\newblock {\em Journal of the American Statistical Association\/}~{\em 91\/}(435), 983--992.

\bibitem[\protect\citeauthoryear{Cook and Weisberg}{Cook and Weisberg}{1991}]{Cook1991}
Cook, R.~D. and S.~Weisberg (1991).
\newblock Sliced inverse regression for dimension reduction: Comment.
\newblock {\em Journal of the American Statistical Association\/}~{\em 86\/}(414), 328--332.

\bibitem[\protect\citeauthoryear{Ferr{\'e} and Yao}{Ferr{\'e} and Yao}{2003}]{Ferre2003fsir}
Ferr{\'e}, L. and A.-F. Yao (2003).
\newblock Functional sliced inverse regression analysis.
\newblock {\em Statistics\/}~{\em 37\/}(6), 475--488.

\bibitem[\protect\citeauthoryear{Forzani and Cook}{Forzani and Cook}{2007}]{Forzanicook2007note}
Forzani, L. and R.~D. Cook (2007).
\newblock A note on smoothed functional inverse regression.
\newblock {\em Statistica Sinica\/}, 1677--1681.

\bibitem[\protect\citeauthoryear{Hall and Horowitz}{Hall and Horowitz}{2007}]{Hall2007mcflr}
Hall, P. and J.~L. Horowitz (2007).
\newblock Methodology and convergence rates for functional linear regression.
\newblock {\em The Annals of Statistics\/}~{\em 35\/}(1), 70--91.

\bibitem[\protect\citeauthoryear{Hall and Li}{Hall and Li}{1993}]{hall1993almost}
Hall, P. and K.-C. Li (1993).
\newblock On almost linearity of low dimensional projections from high dimensional data.
\newblock {\em The annals of Statistics\/}, 867--889.

\bibitem[\protect\citeauthoryear{Hsing and Carroll}{Hsing and Carroll}{1992}]{Hsing1992}
Hsing, T. and R.~J. Carroll (1992).
\newblock An asymptotic theory for sliced inverse regression.
\newblock {\em The Annals of Statistics\/}, 1040--1061.

\bibitem[\protect\citeauthoryear{Hsing and Eubank}{Hsing and Eubank}{2015}]{hsing2015theoretical}
Hsing, T. and R.~Eubank (2015).
\newblock {\em Theoretical foundations of functional data analysis, with an introduction to linear operators}, Volume 997.
\newblock John Wiley \& Sons.

\bibitem[\protect\citeauthoryear{Huang, Tian, and Lin}{Huang et~al.}{2023}]{huang2023sliced}
Huang, D., S.~Tian, and Q.~Lin (2023).
\newblock Sliced inverse regression with large structural dimensions.
\newblock {\em arXiv preprint arXiv:2305.04340\/}.

\bibitem[\protect\citeauthoryear{Lei}{Lei}{2014}]{lei2014adaptive}
Lei, J. (2014).
\newblock Adaptive global testing for functional linear models.
\newblock {\em Journal of the American Statistical Association\/}~{\em 109\/}(506), 624--634.

\bibitem[\protect\citeauthoryear{Li and Wang}{Li and Wang}{2007}]{LiBing2007}
Li, B. and S.~Wang (2007).
\newblock On directional regression for dimension reduction.
\newblock {\em Journal of the American Statistical Association\/}~{\em 102\/}(479), 997--1008.

\bibitem[\protect\citeauthoryear{Li}{Li}{1991}]{li1991sliced}
Li, K.-C. (1991).
\newblock Sliced inverse regression for dimension reduction.
\newblock {\em Journal of the American Statistical Association\/}~{\em 86\/}(414), 316--327.

\bibitem[\protect\citeauthoryear{Li}{Li}{1992}]{KerChauLi1992}
Li, K.-C. (1992).
\newblock On principal hessian directions for data visualization and dimension reduction: Another application of stein's lemma.
\newblock {\em Journal of the American Statistical Association\/}~{\em 87\/}(420), 1025--1039.

\bibitem[\protect\citeauthoryear{Li and Hsing}{Li and Hsing}{2010}]{LiHsing2010dec}
Li, Y. and T.~Hsing (2010).
\newblock Deciding the dimension of effective dimension reduction space for functional and high-dimensional data.
\newblock {\em The Annals of Statistics\/}~{\em 38\/}(5), 3028--3062.

\bibitem[\protect\citeauthoryear{Lian}{Lian}{2015}]{lianfsdr2015}
Lian, H. (2015).
\newblock Functional sufficient dimension reduction: Convergence rates and multiple functional case.
\newblock {\em Journal of Statistical Planning and Inference\/}~{\em 167}, 58--68.

\bibitem[\protect\citeauthoryear{Lian and Li}{Lian and Li}{2014}]{lian2014Sefsdr}
Lian, H. and G.~Li (2014).
\newblock Series expansion for functional sufficient dimension reduction.
\newblock {\em Journal of Multivariate Analysis\/}~{\em 124}, 150--165.

\bibitem[\protect\citeauthoryear{Lin, Li, Huang, and Liu}{Lin et~al.}{2021a}]{lin2021optimality}
Lin, Q., X.~Li, D.~Huang, and J.~S. Liu (2021a).
\newblock On the optimality of sliced inverse regression in high dimensions.
\newblock {\em The Annals of Statistics\/}~{\em 49\/}(1), 1--20.

\bibitem[\protect\citeauthoryear{Lin, Li, Huang, and Liu}{Lin et~al.}{2021b}]{lin2021optimalitysupplement}
Lin, Q., X.~Li, D.~Huang, and J.~S. Liu (2021b).
\newblock Supplementary to “on optimality of sliced inverse regression in high dimensions”.
\newblock {\em The Annals of Statistics\/}~{\em 49\/}(1).

\bibitem[\protect\citeauthoryear{Lin, Zhao, and Liu}{Lin et~al.}{2018a}]{linzl2018}
Lin, Q., Z.~Zhao, and J.~S. Liu (2018a).
\newblock On consistency and sparsity for sliced inverse regression in high dimensions.
\newblock {\em The Annals of Statistics\/}~{\em 46\/}(2), 580--610.

\bibitem[\protect\citeauthoryear{Lin, Zhao, and Liu}{Lin et~al.}{2018b}]{lin2018supplement}
Lin, Q., Z.~Zhao, and J.~S. Liu (2018b).
\newblock Supplement to “on consistency and sparsity for sliced inverse regression in high dimensions.”.

\bibitem[\protect\citeauthoryear{Seelmann}{Seelmann}{2014}]{Seelmann2014NotesOT}
Seelmann, A. (2014).
\newblock Notes on the sin 2$\theta$ theorem.

\bibitem[\protect\citeauthoryear{Tan, Shi, and Yu}{Tan et~al.}{2020}]{tan2020sparse}
Tan, K., L.~Shi, and Z.~Yu (2020).
\newblock Sparse sir: Optimal rates and adaptive estimation.
\newblock {\em The Annals of Statistics\/}~{\em 48\/}(1), 64--85.

\bibitem[\protect\citeauthoryear{Tian, Yu, and Chen}{Tian et~al.}{2023}]{tian2023functional}
Tian, S., Z.~Yu, and R.~Chen (2023).
\newblock Functional slicing-free inverse regression via martingale difference divergence operator.
\newblock {\em arXiv preprint arXiv:2307.12537\/}.

\bibitem[\protect\citeauthoryear{Tsybakov}{Tsybakov}{2009}]{Tsybakov:1315296}
Tsybakov, A.~B. (2009).
\newblock {\em {Introduction to Nonparametric Estimation}}.
\newblock Springer series in statistics. Dordrecht: Springer.

\bibitem[\protect\citeauthoryear{Wang and Lian}{Wang and Lian}{2020}]{Wanglain2020}
Wang, G. and H.~Lian (2020).
\newblock Functional sliced inverse regression in a reproducing kernel hilbert space.
\newblock {\em Statistica Sinica\/}~{\em 30\/}(1), 17--33.

\bibitem[\protect\citeauthoryear{Wu and Li}{Wu and Li}{2011}]{wu2011asymptotic}
Wu, Y. and L.~Li (2011).
\newblock Asymptotic properties of sufficient dimension reduction with a diverging number of predictors.
\newblock {\em Statistica Sinica\/}~{\em 2011\/}(21), 707.

\bibitem[\protect\citeauthoryear{Xia, Tong, Li, and Zhu}{Xia et~al.}{2009}]{Xia2009}
Xia, Y., H.~Tong, W.~K. Li, and L.-X. Zhu (2009).
\newblock An adaptive estimation of dimension reduction space.
\newblock pp.\  299--346.

\bibitem[\protect\citeauthoryear{Yang}{Yang}{1977}]{yang1977}
Yang, S.~S. (1977).
\newblock General distribution theory of the concomitants of order statistics.
\newblock {\em The Annals of Statistics\/}~{\em 5\/}(5), 996--1002.

\bibitem[\protect\citeauthoryear{Yao, Lei, and Wu}{Yao et~al.}{2015}]{yao2015effective}
Yao, F., E.~Lei, and Y.~Wu (2015).
\newblock Effective dimension reduction for sparse functional data.
\newblock {\em Biometrika\/}~{\em 102\/}(2), 421--437.

\bibitem[\protect\citeauthoryear{Yu}{Yu}{1997}]{yu1997assouad}
Yu, B. (1997).
\newblock {Assouad, Fano, and Le Cam}.
\newblock In {\em Festschrift for Lucien Le Cam}, pp.\  423--435. Springer.

\bibitem[\protect\citeauthoryear{Zhu, Miao, and Peng}{Zhu et~al.}{2006}]{zhu2006sliced}
Zhu, L., B.~Miao, and H.~Peng (2006).
\newblock On sliced inverse regression with high-dimensional covariates.
\newblock {\em Journal of the American Statistical Association\/}~{\em 101\/}(474).

\bibitem[\protect\citeauthoryear{Zhu and Ng}{Zhu and Ng}{1995}]{Zhu1995}
Zhu, L.-X. and K.~W. Ng (1995).
\newblock Asymptotics of sliced inverse regression.
\newblock {\em Statistica Sinica\/}, 727--736.

\end{thebibliography}

\bigskip
\vskip .65cm
\noindent
Rui Chen, Center for Statistical Science, Department of Industrial Engineering, Tsinghua University
\vskip 2pt
\noindent
E-mail: chenrui\_fzu@163.com
\vskip 2pt

\noindent
Songtao Tian, Department of Mathematical Sciences, Tsinghua University
\vskip 2pt
\noindent
E-mail: tst20@mails.tsinghua.edu.cn

\noindent
Dongming Huang, Department of Statistics and Data Science, National University of Singapore
\vskip 2pt
\noindent
E-mail: stahd@nus.edu.sg

\noindent
Qian Lin, Department of Statistics and Data Science, Tsinghua University;~Beijing Academy of Artificial Intelligence, Beijing, 100084, China
\vskip 2pt
\noindent
E-mail: qianlin@tsinghua.edu.cn

\noindent
Jun S. Liu, Department of Statistics, Harvard University
\vskip 2pt
\noindent
E-mail: jliu@stat.harvard.edu

\newpage
\appendix
\clearpage            
\pagenumbering{arabic}  
\setcounter{page}{1}  
{\centering
\Large\bfseries Supplement to `On the Optimality of Functional Sliced Inverse Regression'
\par}
\vspace{1em}

\begin{center}
\begin{tabular}{c}
Rui Chen\textsuperscript{a}, Songtao Tian\textsuperscript{b,*}, Dongming Huang\textsuperscript{c}, 
Qian Lin\textsuperscript{d}, Jun S. Liu\textsuperscript{e,**}
\end{tabular}
\end{center}

\begin{center}
\small\itshape
\renewcommand{\arraystretch}{0.7} 
 \begin{tabular}{>{\centering\arraybackslash}p{0.9\textwidth}}
\textsuperscript{a}Department of Statistics and Data Science, Tsinghua University \\[0.01em]
\textsuperscript{b}Department of Mathematical Sciences, Tsinghua University \\
\textsuperscript{c}Department of Statistics and Data Science, National University of Singapore \\
\textsuperscript{d}Department of Statistics and Data Science, Tsinghua University;\\Beijing Academy of Artificial Intelligence, Beijing, 100084, China \\
\textsuperscript{e}Department of Statistics, Harvard University, 02138
\end{tabular}
\end{center}

\footnotetext[1]{*Co-first author.}
\footnotetext[2]{**Corresponding author.}
\bigskip
Section \ref{sec:appendix:proog of lemma} -  \ref{sec:appendix:assist}  contain the proofs of lemmas, propositions and theorems.
Section \ref{sec:appendix:simulation} contains additional simulation results
of Section \ref{sec:simulationfsir}.
\section{Proof of Lemma \ref{lem:monment condition to sliced stable}}\label{sec:appendix:proog of lemma}

The proof follows the same argument as in the proof the \cite[Theorem $1$ ]{huang2023sliced} and we only need to generalize the multivariate result therein to a functional version. We omit the proof here for simplicity. 

\section{Proof of Proposition \ref{prop:lgvgm}}

Suppose $H$ is any integer greater than the constant $K$ in Assumption~\ref{as: slice stable}.

For any $\bs{u}\in\mb{S}_{\H}$, it holds that
	\begin{equation}\label{eq:gammebeta}
	\begin{aligned}
		\langle \Gamma_{e}(\bs{u}),\bs{u}\rangle &= \int(\E[\langle \boldsymbol{X},\bs{u}\rangle \mid Y=y])^{2} d P_{Y}(y)=\sum_{h:\mathcal{S}_h\in \mathfrak{S}_H(n)}\int_{\mathcal{S}_{h}}(\E[\langle \boldsymbol{X},\bs{u}\rangle \mid Y=y])^{2} d P_{Y}(y) \\
		&=\sum_{h:\mathcal{S}_h\in \mathfrak{S}_H(n)} \bbP(Y\in \mathcal{S}_{h}) \E\left((\E[\langle \boldsymbol{X},\bs{u}\rangle \mid Y=y])^{2} \mid Y\in \mathcal{S}_{h} \right)\\
		&=\sum_{h:\mathcal{S}_h\in \mathfrak{S}_H(n)}\delta_{h}\E[\bs m^2(\bs{u})\mid Y\in \mathcal{S}_{h}]
	\end{aligned}
\end{equation}
	where $\delta_{h}:=\bbP(Y\in \mathcal{S}_{h})$ and $\m(\bs{u}):=\langle\m(Y),\bs{u}\rangle $. 
Furthermore, it holds that 
	\begin{equation}\label{eq:mhbeta}
		\begin{aligned}
		\langle (\overline{\m}_{h}\otimes \overline{\m}_{h}(\bs{u}),\bs{u}\rangle &=\E^2[\m(\bs{u})\mid Y\in \mathcal{S}_{h}]=\E[\m^2(\bs{u})\mid Y\in \mathcal{S}_{h}]-\var(\m(\bs{u})\mid Y\in \mathcal{S}_{h}).
		\end{aligned}
		\end{equation}

For any $\gamma$-partition $\mathfrak{S}_{H}(n):=\{\mathcal{S}_{h}, h=1,..,H\}$, it holds that 
\begin{align*}
&\left|\left\langle\left( \widetilde{\Gamma}_{e}-\Gamma_{e}\right)(\bs{u}),\bs{u}  \right\rangle\right|=\left|\frac{1}{H}\sum_{h:\mathcal{S}_h\in \mathfrak{S}_H(n)}	\langle (\overline{\m}_{h}\otimes \overline{\m}_{h}(\bs{u}),\bs{u}\rangle - \sum_{h:\mathcal{S}_h\in \mathfrak{S}_H(n)}\delta_{h}\E[\bs m^2(\bs{u})\mid Y\in \mathcal{S}_{h}] \right|\\
\leqslant&\left|\sum_{h:\mathcal{S}_h\in \mathfrak{S}_H(n)}(\frac{1}{H}-\delta_{h})\E(\m^2(\bs{u})\mid Y\in \mathcal{S}_{h})\right|+\left|\frac{1}{H}\sum_{h:\mathcal{S}_h\in \mathfrak{S}_H(n)}\var(\m(\bs{u})\mid Y\in \mathcal{S}_{h})\right|\\
\leqslant& \frac{1}{\tau-1}\left|\sum_{h:\mathcal{S}_h\in \mathfrak{S}_H(n)}\delta_{h}\E[\m^2(\bs{u})\mid Y\in \mathcal{S}_{h}]\right|+\left|\frac{1}{H}\sum_{h:\mathcal{S}_h\in \mathfrak{S}_H(n)}\var(\m(\bs{u})\mid Y\in \mathcal{S}_{h})\right| 
\\
\leqslant  &\frac2\tau \left\langle\Gamma_{e}(\bs{u}),\bs{u}\right\rangle+\left|\frac{1}{H}\sum_{h:\mathcal{S}_h\in \mathfrak{S}_H(n)}\var(\m(\bs{u})\mid Y\in \mathcal{S}_{h})\right|\\
\leqslant  &\frac3\tau \left\langle\Gamma_{e}(\bs{u}),\bs{u}\right\rangle,
\end{align*}
where the first inequality is due to \eqref{eq:mhbeta}, the second inequality is because $1-\gamma \leqslant H \delta_{h}$ since $\mathfrak{S}_{H}(n)$ is a $\gamma$-partition, the third is because $\gamma\leqslant\frac1\tau$, and the last is due to the Assumption~\ref{as: slice stable}.

Recall Lemma  \ref{lem:sliced}, which states that there is some $H'$ and $C$ such that for any $H>H'$ and  $n>\frac{4H}{\gamma}+1$, the sample sliced partition is a $\gamma$-partition with probability at least 
		\[1-CH^2\sqrt{n+1}\exp\left(-\frac{\gamma^2(n+1)}{32H^2}\right)\]
		for some absolute constant $C$. Then the proof is completed by choosing $H_0:=\max\{H',K\}$.
\qed
\section{Proof of Lemma  \ref{lem:cvgem}}
By Proposition \ref{prop:lgvgm}, the first statement in Lemma  \ref{lem:cvgem}  is a direct corollary of the second one. Now we begin to prove the second statement.

 Define  an event $\ttE$  as follows:
 \begin{equation}\label{eq: event E}
\ttE:=\{\mathfrak{S}_{H}(n) \text{ is a } \gamma\text{-partition} \}.
 \end{equation}
From the proof of  Proposition \ref{prop:lgvgm}, we know that 
\begin{align*}
\bbP(\ttE)\geqslant 1-CH^2\sqrt{n+1}\exp\left(-\frac{\gamma^2(n+1)}{32H^2}\right).
\end{align*}
and we have the following lemma. 
\begin{lemma}\label{lem:WSSC+ttE}
 Under WSSC, on the event $\ttE$, the following holds: 
\begin{align*}
\left|\left\langle\left( \widetilde{\Gamma}_{e}-\Gamma_{e}\right)(\bs{u}),\bs{u}  \right\rangle\right|
\leqslant \frac3\tau \left\langle\Gamma_{e}(\bs{u}),\bs{u}\right\rangle, \forall\bs{u}\in\mb{\mb{S}}_{\H}.
\end{align*}   
\end{lemma}

 The rest of our effort will be devoted to establishing the following lemma, which directly implies the second statement in Lemma~\ref{lem:cvgem}. 
\begin{lemma}\label{lem:norm control of hatGammae tildeGammae}
Suppose that  Assumptions \ref{as: slice stable} and \ref{as:moment} hold.
For any fixed integer $H>H_{0}$ ($H_0$ is defined in Proposition \ref{prop:lgvgm}) and any sufficiently large  $n>1+4H/\gamma$,  we have 
\begin{align*}
\E\left[1_{\ttE} \norm{\widehat{\Gamma}_{e}-\widetilde{\Gamma}_{e}}^2 \right]\lesssim \frac{H^2}{n}.    
\end{align*}
\end{lemma}

\noindent\textit{Proof of Lemma \ref{lem:norm control of hatGammae tildeGammae}.}
For any $\bs{u}\in\mathbb{S}_\mathcal{H}$, let $\bs{u}:=\sum_{j\geqslant 1}b_{j}\phi_{j}$ with $\sum_{j\geqslant 1}b^{2}_{j}=1$ where $\{\phi_j\}_{j=1}^\infty$ is the eigenfunctions of $\Gamma$. Then we have 
\begin{align*}
\langle (\widehat{\Gamma}_{e}-\widetilde{\Gamma}_{e})\bs{u}, \bs{u}\rangle =& \sum_{ i,j\geqslant 1}b_{i}b_{j}\langle (\widehat{\Gamma}_{e}-\widetilde{\Gamma}_{e})\phi_{i},\phi_{j}\rangle \\
\leqslant & \left(\sum_{ i,j\geqslant 1}b^{2}_{i}b^{2}_{j}\right)^{\tfrac{1}{2}}\left(\sum_{i,j\geqslant 1}\langle (\widehat{\Gamma}_{e}-\widetilde{\Gamma}_{e})\phi_{i},\phi_{j}\rangle ^{2}\right)^{\tfrac{1}{2}}\\
=&\left(\sum_{i,j\geqslant 1}\langle (\widehat{\Gamma}_{e}-\widetilde{\Gamma}_{e})\phi_{i},\phi_{j}\rangle ^{2}\right)^{\tfrac{1}{2}}.
\end{align*}
It shows that $\norm{\widehat{\Gamma}_{e}-\widetilde{\Gamma}_{e}}\leqslant \left(\sum_{i,j\geqslant 1}\langle (\widehat{\Gamma}_{e}-\widetilde{\Gamma}_{e})\phi_{i},\phi_{j}\rangle ^{2}\right)^{\tfrac{1}{2}}$.

Let $\xi_{hj}:=\langle \overline{\bs{m}}_{h},\phi_{j}\rangle $ and $\widehat{\xi}_{hj}:=\langle \overline{\boldsymbol{X}}_{h,\cdot},\phi_{j}\rangle $. 
The operators $\widehat{\Gamma}_{e}$ and $\widetilde{\Gamma}_{e}$ can be written as 
\begin{align*}
\widehat{\Gamma}_{e}=&\frac{1}{H}\sum_{h:\S_{h}\in\mathfrak{S}_{H}(n)}\left(\sum_{j\geqslant 1}\widehat{\xi}_{hj}\phi_{j}\right)\otimes \left(\sum_{j\geqslant 1}\widehat{\xi}_{hj}\phi_{j}\right),  \\ 
 \widetilde{\Gamma}_{e}=&\frac{1}{H}\sum_{h:\S_{h}\in\mathfrak{S}_{H}(n)}\left(\sum_{j\geqslant 1}\xi_{hj}\phi_{j}\right)\otimes \left(\sum_{j\geqslant 1}\xi_{hj}\phi_{j}\right),
\end{align*}
 respectively. Now	we  obtain
	\begin{align*}
		&\sum_{i,j\geqslant 1}\langle (\widehat{\Gamma}_{e}-\widetilde{\Gamma}_{e})\phi_{i},\phi_{j}\rangle ^{2}
		\leqslant\frac{1}{H}\sum_{i,j\geqslant 1}\sum_{h:\S_{h}\in\mathfrak{S}_{H}(n)}\left(\widehat{\xi}_{hi}\widehat{\xi}_{hj}-\xi_{hi}\xi_{hj}   \right)^{2}\\
		= &\frac{1}{H} \sum_{i,j\geqslant 1}\sum_{h:\S_{h}\in\mathfrak{S}_{H}(n)}\left((\widehat{\xi}_{hi}-\xi_{hi})(\widehat{\xi}_{hj}-\xi_{hj})+\xi_{hj}(\widehat{\xi}_{hi}-\xi_{hi})+\xi_{hi}(\widehat{\xi}_{hj}-\xi_{hj})\right)^{2}.
	\end{align*}  

The lemma is proved if we can show that  for any $h$,
\begin{align*}
	A:=	\E\left[1_{\ttE} \sum_{i,j\geqslant 1}\left((\widehat{\xi}_{hi}-\xi_{hi})(\widehat{\xi}_{hj}-\xi_{hj})+\xi_{hj}(\widehat{\xi}_{hi}-\xi_{hi})+\xi_{hi}(\widehat{\xi}_{hj}-\xi_{hj})\right)^{2}\right]\leqslant \frac{C^{\prime}H^{2}}{n}
\end{align*}
for some  positive constant $C^{\prime}$.

Note that 
\begin{align*}
	& \E\left[1_{\ttE}\sum_{i,j\geqslant 1}\left((\widehat{\xi}_{hi}-\xi_{hi})(\widehat{\xi}_{hj}-\xi_{hj})+\xi_{hj}(\widehat{\xi}_{hi}-\xi_{hi})+\xi_{hi}(\widehat{\xi}_{hj}-\xi_{hj})\right)^{2}\right]\\
	\leqslant &3\E\left[1_{\ttE}\sum_{i,j\geqslant 1}(\widehat{\xi}_{hi}-\xi_{hi})^{2}(\widehat{\xi}_{hj}-\xi_{hj})^{2}+\xi^{2}_{hj}(\widehat{\xi}_{hi}-\xi_{hi})^{2}+\xi^{2}_{hi}(\widehat{\xi}_{hj}-\xi_{hj})^{2}\right]\\
	=& 3\E\left[1_{\ttE}\sum_{i,j\geqslant 1}(\widehat{\xi}_{hi}-\xi_{hi})^{2}(\widehat{\xi}_{hj}-\xi_{hj})^{2}\right]+6\E\left[1_{\ttE}\sum_{i,j\geqslant 1}\xi^{2}_{hj}(\widehat{\xi}_{hi}-\xi_{hi})^{2}\right]\\
	:= & 3I+6II.
\end{align*}
We first bound the term $II$. 
Recall that  $\langle \boldsymbol{X},\phi_{j}\rangle=\xi_{j}, \forall j$.
Then  there exists a constant $C_{1}>0$, such that for all $j\geqslant 1$,  we have 
\begin{align}\label{eq:bound expection square}
 1_{\ttE}\xi^{2}_{h j}=1_{\ttE}\E^2[\xi_{j}|Y\in\S_h]
\overset{(a)}{\leqslant}1_{\ttE}\E[\xi_{j}^2|Y\in\S_h]
 \overset{(b)}{\leqslant}1_{\ttE}\E^{1/2}[\xi_{j}^4|Y\in\S_h]\overset{(c)}{\leqslant} C_{1}\sqrt H\E^{1/2}[\xi_{j}^4]\overset{(d)}{\leqslant} C \sqrt H\lambda_{j}
\end{align}
where we have used Jensen inequality for conditional expectation in $(a)$ and $(b)$, Lemma{
\ref{lem:fourth_moment:lemma} in $(c)$, and Assumption \ref{as:moment} in $(d)$.

Assume $\X_{h,j}=\sum_{i=1}^\infty\xi_{h,j,i}\phi_i$, then one has
\begin{align*}
\E[1_{\ttE}(\widehat{\xi}_{hi}-\xi_{hi})^{2}]
=&\E[1_{\ttE}\langle \frac{1}{c}\sum_{j=1}^c \X_{h,j}-\E[\X|Y\in\S_h] ,\phi_i\rangle^2]\\
=&\E[1_{\ttE}(\frac1c\sum_{j=1}^c\xi_{h,j,i}-\E[\xi_i|Y\in\S_h])^2]\\
=&\E[1_{\ttE}(\frac1{c-1}\sum_{j=1}^{c-1}\xi_{h,j,i}-\E[1_{\ttE}\xi_i|Y\in\S_h]+\frac{ \xi_{h,c,i}}{c}-\frac1{c(c-1)}\sum_{j=1}^{c-1}\xi_{h,j,i})^2]\\
\lesssim& \E[1_{\ttE}(\frac1{c-1}\sum_{j=1}^{c-1}\xi_{h,j,i}-\E[1_{\ttE}\xi_i|Y\in \S_h])^2]+\frac1{c^2}\E[1_{\ttE}\xi_{h,c,i}^2]+\frac{1}{c^2(c-1)}\E[1_{\ttE}\sum_{j=1}^{c-1}\xi_{h,j,i}^2]\\
\overset{(a)}{\leqslant}&\frac{1}{c-1}\E[1_{\ttE}\xi_{h,1,i}^2]+\frac1{c^2}\E[1_{\ttE}\xi_{h,1,i}^2]+\frac{1}{c^2}\E[1_{\ttE}\xi_{h,1,i}^2]\\
\lesssim& \frac{1}{c-1}\E[1_{\ttE}\xi_{h,1,i}^2]\overset{(b)}{\leqslant}\frac{1}{c-1}\mb \E^{1/2}[1_{\ttE}\xi_{h,1,i}^4]\\
\overset{(c)}{\lesssim}& \frac{1}{c-1}\sqrt H\E^{1/2}[\xi_{i}^4]\overset{(d)}{\leqslant} \frac{1}{c-1}\sqrt H\lambda_i\leqslant C^{\prime}_{1}H^{3/2}\lambda_{i}/n,
\end{align*}

where we have used Jensen inequality for conditional expectation in $(b)$, Lemma \ref{lem:fourth_moment:lemma} in $(c)$, and Assumption \ref{as:moment} in $(d)$. Furthermore, $(a)$ is based on Lemma~\ref{cor:i.i.d} and the following derivation:

\begin{align*} 
& \E[(\frac1{c-1}\sum_{j=1}^{c-1}\xi_{h,j,i}-\E[\xi_i|Y\in\S_h])^2] \\ 
 &= \E\left[ \E[(\frac1{c-1}\sum_{j=1}^{c-1}\xi_{h,j,i}-\E[\xi_i|Y\in\S_h])^2 \mid \{\S_{h'}\}_{h'=1}^H ] \right]\\
&=\E\left[ \frac{1}{c-1}\var[\xi_{h,1,i}   \mid \{\S_{h'}\}_{h'=1}^H ]  \mid \{\S_{h'}\}_{h'=1}^H \right]\\
& \leqslant  \frac{1}{c-1}\E[\xi_{h,1,i}^2 ].
\end{align*}
Thus 
\begin{equation}\label{eq:termii}
	II\leqslant \sum_{i,j\geqslant 1}\frac{C_{1}H^{2}\lambda_{i}\lambda_{j}}{n}
\end{equation}
for some sufficiently large $C_{1}$. 

Next we handle  the term $I$.  
From \eqref{eq:bound expection square} and Cauchy-Schwarz inequality, one has: $1_{\ttE}\left|\E[\xi_i^3|Y\in\S_h]\right|\leqslant CH^{3/4}\lambda_i^{3/2}$.
By direct calculation of fourth moment  (see, e.g., \cite[Theorem 1]{angelova2012moments}), one has
\begin{align*}
\E[1_{\ttE}(\widehat{\xi}_{hi}-\xi_{hi})^{4}]\leqslant C\E[1_{\ttE}(\frac1{c-1}\sum_{j=1}^{c-1}\xi_{h,j,i}-\E[\xi_i|Y\in\S_h])^4]\leqslant C\frac{H^3\lambda_i^2}{n^2}\leqslant C\frac{H^2\lambda_i^2}{n}.
\end{align*}

We get 
\begin{equation}\label{eq:termi}
	\begin{aligned}
		I=&\E\left[1_{\ttE}\sum_{i,j\geqslant 1}(\widehat{\xi}_{hi}-\xi_{hi})^{2}(\widehat{\xi}_{hj}-\xi_{hj})^{2}\right] \\
		\leqslant& \sum_{i,j\geqslant 1}\E^{\tfrac{1}{2}}\left[1_{\ttE}(\widehat{\xi}_{hi}-\xi_{hi})^{4}\right]\E^{\tfrac{1}{2}}\left[1_{\ttE}(\widehat{\xi}_{hj}-\xi_{hj})^{4}\right]\\
		\leqslant&\sum_{i,j\geqslant 1}\frac{C_{2}H^2\lambda_{i}\lambda_{j}}{n}.
	\end{aligned}
\end{equation}
Since $\Gamma$ has a finite trace, we can
take $C^{\prime}=6\sum_{i,j\geqslant 1}(C_1+C_{2})\lambda_{i}\lambda_{j}<\infty$. Then by Equations \eqref{eq:termii}  and \eqref{eq:termi}, we have 
\begin{align*}
	A\leqslant \frac{C^{\prime}H^{2}}{n}
\end{align*}
as required.
\qed

\section{Proof of Theorem \ref{thm:rootnconsistency}}
We first provide the following lemma.
\begin{lemma}\label{lem:image equality}
 Suppose Assumption \ref{as: slice stable} holds with  $\tau>\frac{6\norm{\Gamma_{e}}}{\lambda^{+}_{\min}(\Gamma_{e})}$. On the event $\ttE$, we have $\mr{Im}(\Gamma_e)=\mr{Im}(\wt\Gamma_e)$  and $\lambda_{\min}^+(\wt\Gamma_e)=\lambda_{d}(\wt\Gamma_e)\geqslant\frac{\lambda_d(\Gamma_e)}{2}$.   
\end{lemma}
\begin{proof}
Recall that $\widetilde{\Gamma}_{e}:=\frac{1}{H}\sum_{h:\mc S_h  \in \mathfrak{S}_{H}(n)}\overline{\bs{m}}_{h}\otimes \overline{\bs{m}}_{h}$ and  $\overline{\bs{m}}_{h} := \mathbb{E}[\bs{m}(Y) \mid Y\in \S_{h}]=\E[\bs X|Y\in\S_h]$. Note that $\overline{\bs{m}}_{h}\in\mc S_e= \mbox{Im}(\Gamma_{e})$    for all $h=1,\dots, H$. Thus, $\mr{Im}(\wt\Gamma_e)\subseteq\mr{Im}(\Gamma_e)$. 

Next we prove by contradiction that $\mr{Im}(\Gamma_e)\subseteq\mr{Im}(\wt\Gamma_e)$. 
Assume that there exists a vector $\u\in\S_{\mc H}$ such that $\Gamma_e\u\neq 0$ and $\wt\Gamma_e\u=0$. Since $\tau>\frac{6\norm{\Gamma_{e}}}{\lambda^{+}_{\min}(\Gamma_{e})}$ and the event $\ttE$ happens, we have 
\begin{align*}
\|\Gamma_e(\u)\| =\|(\Gamma_e-\wt\Gamma_e)(\u)\| \leqslant\frac{3}{\tau}\|\Gamma_e\|\leqslant\frac{\lambda_{\min}^+(\Gamma_e)}{2}, 
\end{align*}
where the first inequality comes from Lemma \ref{lem:WSSC+ttE}.
This is a contradiction because $0<\lambda_{\min}^+(\Gamma_e)\leqslant \|\Gamma_e(\u)\|$. Thus $\mr{Im}(\Gamma_e)=\mr{Im}(\wt\Gamma_e)$ and $\mr{rank}(\wt\Gamma_e)=\mr{rank}(\Gamma_e)=d$. By Lemma \ref{lem:wely ineq operator}, we have 
\begin{align*}
 |\lambda_d(\wt\Gamma_e)-\lambda_d(\Gamma_e)|\leqslant\|\wt\Gamma_e-\Gamma\|\leqslant\frac{3}{\tau}\|\Gamma_e\|\leqslant\frac{\lambda_d(\Gamma_e)}{2}
\end{align*}
where the second inequality comes from Proposition \ref{prop:lgvgm}.
\end{proof}

Under the condition of Theorem \ref{thm:rootnconsistency} with   $\tau>\frac{6\norm{\Gamma_{e}}}{\lambda^{+}_{\min}(\Gamma_{e})}$ and the event $\ttE$, then  
$\mr{Im}(\widetilde{\Gamma}_{e})=\mr{Im}(\Gamma_{e})$ and thus $\mr{rank}(\wt\Gamma_e)=\mr{rank}(\Gamma_e)=d$. 
Let $\{\widehat{\mu}_{i}\}^{d}_{i=1}$ be the $d$ largest eigenvalues of 
$\widehat{\Gamma}_{e}$ with associated eigenfunctions $\{\widehat{v}_{i}\}^{d}_{i=1}$. Define
\begin{equation}\label{eq: def of hatGammaed}
 \widehat{\Gamma}^{d}_{e}:=\sum^{d}_{i=1}\widehat{\mu}_{i} \wh v_{i}\otimes \wh v_{i}.   
\end{equation}
We first prove that
\begin{align*}
	\mb E\left[\norm{\widehat{\Gamma}_{e}^d-\widetilde{\Gamma}_{e}}^2 1_{\ttE}\right]\lesssim \frac{H^{2}}{n}.
\end{align*}
Using Lemma $\ref{lem:wely ineq operator}$, one can get
$
\lambda_i^2\left(\widehat{\Gamma}_{e}\right)\leqslant 2\left(\left\|\widehat{\Gamma}_{e}-\widetilde{\Gamma}_{e}\right\|^2+\lambda_i^2\left( \widetilde{\Gamma}_{e}\right)\right)
$. 
Since $\mr{rank}(\widetilde{\Gamma}_{e})=d$, one can get $\lambda_i(\widetilde{\Gamma}_{e})=0,~i\geqslant d+1$. Thus by Lemma~\ref{lem:norm control of hatGammae tildeGammae}, one has $\E\left[1_{\ttE} \norm{\widehat{\Gamma}_{e}-\widetilde{\Gamma}_{e}}^2 \right]\lesssim \frac{H^2}{n}$, which implies 
\begin{align*}
\mb E\left[\left|\lambda_i\left(\widehat{\Gamma}_{e}\right)\right|^21_{\ttE}\right]\lesssim {\frac{H^{2}}{n}}\quad(i\geqslant d+1).
\end{align*}

Furthermore, 
\begin{align*}
&\mb E\left[\left\|\widehat{\Gamma}_{e}^d- \widetilde{\Gamma}_{e}\right\|^21_{\ttE}\right]\leqslant 2\left(\mb E\left[\left\| \widetilde{\Gamma}_{e}-\widehat{\Gamma}_{e}\right\|^21_{\ttE}+\left\|\widehat{\Gamma}_{e}-\widehat{\Gamma}_{e}^d\right\|^21_{\ttE}\right]\right)\\
=&
2\mb E\left[\left\| \widetilde{\Gamma}_{e}-\widehat{\Gamma}_{e}\right\|^21_{\ttE}+\lambda_{d+1}^2(\wh\Gamma_e)1_{\ttE}\right]\lesssim {\frac{H^{2}}{n}}.
\end{align*}

By Lemma~\ref{lem:image equality}, we have
$\lambda^{+}_{\min}(\widetilde{\Gamma}_{e})\geqslant\lambda_d(\Gamma_e)/2$. 

By Markov inequality, we have 
\begin{align*}
&\bbP(|\lambda_d(\widehat{\Gamma}^{d}_{e})-\lambda_d(\widetilde{\Gamma}_{e})|1_{\ttE}\leqslant\frac{\lambda_d(\widetilde{\Gamma}_{e})}{2})\geqslant
\bbP(|\lambda_d(\widehat{\Gamma}^{d}_{e})-\lambda_d(\widetilde{\Gamma}_{e})|1_{\ttE}\leqslant\frac{\lambda_d(\Gamma_{e})}{4})\\
\geqslant&1-\frac{\bbE[|\lambda_d(\widehat{\Gamma}^{d}_{e})-\lambda_d(\widetilde{\Gamma}_{e})|^21_{\ttE}]}{\lambda_d^2(\Gamma_{e})/16}\\
\geqslant& 1-\frac{\bbE[\|\widehat{\Gamma}^{d}_{e}-\widetilde{\Gamma}_{e}\|^21_{\ttE}]}{\lambda_d^2(\Gamma_{e})/16}
\geqslant 1-C'{\frac{H^2}{ n}}
\end{align*}
for some constant $C'>0$. 
Define $\widetilde\ttF:=\{|\lambda_d(\widehat{\Gamma}^{d}_{e})-\lambda_d(\widetilde{\Gamma}_{e})|1_{\ttE}\leqslant\frac{\lambda_d(\widetilde{\Gamma}_{e})}{2}\}$,  we know under the event $\ttE\cap\wt\ttF$, it hold that $\lambda_d(\widehat{\Gamma}^{d}_{e})\geqslant\lambda_d(\widetilde{\Gamma}_{e})/2\geqslant\lambda_d(\Gamma_e)/4$. Thus we  have $\min\left\{\lambda^{+}_{\min}(\widehat{\Gamma}_{e}^d),\lambda^{+}_{\min}(\widetilde{\Gamma}_{e})\right\}1_{\ttE\cap\wt\ttF}\geqslant \lambda_d(\Gamma_e)/4$.

 Then Applying  $\sin\Theta$ theorem (Lemma \ref{lemma, sin theta of infinite dimension operator}), we have 
 \begin{align*}
 \mb E\left[\norm{ P_{\widehat{\S}_{e}}-P_{\S_{e}}}^2 1_{\ttE\cap\wt\ttF} \right]\lesssim 
 \frac{\mb E\left[\norm{\widehat{\Gamma}_{e}^d-\widetilde{\Gamma}_{e}}^21_{\ttE\cap\wt\ttF}\right]}{\min\left\{\lambda^{+}_{\min}(\widehat{\Gamma}_{e}^d),\lambda^{+}_{\min}(\widetilde{\Gamma}_{e})\right\}^21_{\ttE\cap\wt\ttF}}\lesssim {\frac{H^{2}}{n}}.
\end{align*}
Thus
\begin{align*}
\mb E\left[\norm{ P_{\widehat{\S}_{e}}-
P_{\S_{e}}}^2\right]=&\mb E\left[\norm{ P_{\widehat{\S}_{e}}-P_{\S_{e}}}^21_{\ttE\cap\wt\ttF}\right]+\mb E\left[\norm{ P_{\widehat{\S}_{e}}-P_{\S_{e}}}^21_{(\ttE\cap\wt\ttF)^c}\right]\\
\lesssim &  {\frac{H^{2}}{n}} +  dCH^2\sqrt{n+1}\exp\left(-\frac{\gamma^2(n+1)}{32H^2}\right)\\
\lesssim& {\frac{H^{2}}{n}}.
\end{align*}
\qed
\section{Proof of Theorem \ref{thm:cfsdrs}}\label{sec: proof of upper bound}
By Markov inequality and $\bbE[|\wh\lambda_j-\lambda_j|^2]\lesssim\frac{1}{ n}$(see e.g., Equation (5.26) in  \cite{Hall2007mcflr}), we have 
\begin{align*}
\bbP(|\wh\lambda_j-\lambda_j|\leqslant\frac{\lambda_j}{2})\geqslant 1-\frac{\bbE[|\wh\lambda_j-\lambda_j|^2]}{\lambda_j^2/4}\geqslant 1-\frac{Cj^{2\alpha}}{n}.
\end{align*}
 
Define 
\begin{equation}\label{eq:event F}
 \ttF:=\left\{\frac{\lambda_i}{2}\leqslant\widehat\lambda_i\leqslant\frac{3\lambda_i}{2},\forall i\in[m]\right\}=\left\{|\widehat\lambda_i-\lambda_i|\leqslant\frac{\lambda_i}{2},\forall i\in[m]\right\}.   
\end{equation}
Then we have $\bbP(\ttF)\geqslant 1-\frac{m^{2\alpha+1}}{n}$.

We first need a preparatory theorem.
\subsection{A preparatory theorem}

\begin{theorem}\label{thm:bt1t2}
With the same conditions as in Theorem \ref{thm:cfsdrs}, we have 
\begin{align*}
\mb E\left[\norm{\widehat{\Gamma}^{\dagger}_{m}\widehat{\Gamma}_{e}-\Gamma^{\dagger}_{m}\widetilde{\Gamma}_{e}}1_{\ttE\cap\ttF}\right]\lesssim Hn^{\frac{-(2\beta-1)}{2(\alpha+2\beta)}}.
\end{align*}
\end{theorem}
\begin{proof}
	We can decompose  	 $\widehat{\Gamma}^{\dagger}_{m}\widehat{\Gamma}_{e}-\Gamma^{\dagger}_{m}\widetilde{\Gamma}_{e}$ as follows:
	\begin{equation}\label{eq:dgammagammae}
		\begin{aligned}
			&\widehat{\Gamma}^{\dagger}_{m}\widehat{\Gamma}_{e}-\Gamma^{\dagger}_{m}\widetilde{\Gamma}_{e}\\
			=& \Gamma^{\dagger}_{m}\left(\widehat{\Gamma}_{e}-\widetilde{\Gamma}_{e}\right)+(\widehat{\Gamma}^{\dagger}_{m}-\Gamma^{\dagger}_{m})\widetilde{\Gamma}_{e}+\left(\widehat{\Gamma}^{\dagger}_{m}-\Gamma^{\dagger}_{m}\right)\left(\widehat{\Gamma}_{e}-\widetilde{\Gamma}_{e}\right)\\
			:=& B_{1}+B_{2}+B_{3}.  
		\end{aligned}
	\end{equation}
	Then our theorem is derived directly by the next proposition.
\end{proof}
\begin{proposition}\label{prop:a1a2a3value}
	With the same conditions as in Theorem \ref{thm:cfsdrs},  we have 
	\begin{itemize}
		\item[(1).]
		$\mb E[\norm{B_{1}}^21_{\ttE}]=\mb E\left[\norm{\Gamma^{\dagger}_{m}\left(\widehat{\Gamma}_{e}-\widetilde{\Gamma}_{e}\right)}^21_{\ttE}\right]\lesssim H^2\frac{m^{\alpha+1}}{n}\lesssim H^2n^{\frac{-(2\beta-1)}{\left(\alpha+2\beta\right)}}$;
		\item[(2).] $\mb E[\norm{B_{2}}^21_{\ttE}]=\mb E\left[\norm{\left(\widehat{\Gamma}^{\dagger}_{m}-\Gamma^{\dagger}_{m}\right)\widetilde{\Gamma}_{e}}^21_{\ttE}\right]\lesssim H^2\frac{m^{\alpha+1}}{n}\lesssim H^2n^{\frac{-(2\beta-1)}{\left(\alpha+2\beta\right)}}$;
		\item[(3).] $\mb E[\norm{B_{3}}1_{\ttE\cap\ttF}]=\mb E\left[\norm{\left(\widehat{\Gamma}^{\dagger}_{m}-\Gamma^{\dagger}_{m}\right)\left(\widehat{\Gamma}_{e}-\widetilde{\Gamma}_{e}\right)} 1_{\ttE\cap\ttF}\right]\lesssim  H\frac{m^{(\alpha+1)/2}}{\sqrt n}\lesssim  Hn^{\frac{-(2\beta-1)}{2\left(\alpha+2\beta\right)}}$. 
	\end{itemize}
\end{proposition}
{\noindent \bf \underline{Proof of Proposition \ref{prop:a1a2a3value}-(1)}.~} 
For any $\bs{u}=\sum_{i\geqslant 1}b_{i}\phi_{i}\in\mb S_{\H}$ with $\sum_{i\geqslant 1}b^{2}_{i}=1$, we have 
\begin{align*}
	\norm{\Gamma^{\dagger}_{m}\left(\widehat{\Gamma}_{e}-\widetilde{\Gamma}_{e}\right)\bs{u}}^{2}&=\sum_{j=1}^{\infty}\left\langle\Gamma^{\dagger}_{m}\left(\widehat{\Gamma}_{e}-\widetilde{\Gamma}_{e}\right)\bs{u},\phi_j\right\rangle^2\\
&= \sum_{j=1}^{\infty}\left(\sum_{i=1}^\infty b_i\left\langle\left(\widehat{\Gamma}_{e}-\widetilde{\Gamma}_{e}\right)(\phi_i),\Gamma^{\dagger}_{m}\phi_j\right\rangle\right)^2\\
& \leqslant\sum_{j=1}^{\infty}\left(\sum_{i=1}^\infty b_i^2\sum_{i=1}^\infty\left\langle\left(\widehat{\Gamma}_{e}-\widetilde{\Gamma}_{e}\right)(\phi_i),\Gamma^{\dagger}_{m}\phi_j\right\rangle^2\right)\\
& =\sum_{i,j\geqslant 1}\left\langle\left(\widehat{\Gamma}_{e}-\widetilde{\Gamma}_{e}\right)\phi_{i}, \Gamma^{\dagger}_{m}\phi_{j}\right\rangle^{2}\\
	&= \sum^{m}_{j=1}\frac{1}{\lambda^{2}_{j}}\sum_{i\geqslant 1}\left\langle\left(\widehat{\Gamma}_{e}-\widetilde{\Gamma}_{e}\right)\phi_{i},\phi_{j}\right\rangle^{2}.
\end{align*}
Then we obtain 
 \begin{align*}
\norm{\Gamma^{\dagger}_{m}\left(\widehat{\Gamma}_{e}-\widetilde{\Gamma}_{e}\right)}^{2}\leqslant \sum^{m}_{j=1}\frac{1}{\lambda^{2}_{j}}\sum_{i\geqslant 1}\left\langle\left(\widehat{\Gamma}_{e}-\widetilde{\Gamma}_{e}\right)\phi_{i},\phi_{j}\right\rangle^{2}.
\end{align*}
By the proof of Lemma \ref{lem:norm control of hatGammae tildeGammae}, we see that there exists some constant $C>0$, such that  
\begin{align*}
	\E\left[1_{\ttE}\sum^{m}_{j=1}\frac{1}{\lambda^{2}_{j}}\sum_{i\geqslant 1}\left\langle\left(\widehat{\Gamma}_{e}-\widetilde{\Gamma}_{e}\right)\phi_{i},\phi_{j}\right\rangle^{2}\right]&\leqslant\sum^{m}_{j=1}\frac{1}{\lambda^{2}_{j}}\sum_{i\geqslant 1}\frac{CH^{2}\lambda_{i}\lambda_{j}}{n}\\
	&=CH^{2} \sum^{m}_{j=1}\frac{1}{\lambda_{j}n}\sum_{i\geqslant 1}\lambda_{i}\\
	&\leqslant C'\frac{CH^{2} m^{\alpha+1}}{n},
\end{align*}
where the last inequality comes from $\lambda_j\geqslant Cj^{-\alpha}$ by Assumption \ref{assumption: rate-type condition}.
It implies that  $\mb E[1_{\ttE}\norm{B_{1}}^{2}]\lesssim\frac{H^{2}m^{\alpha+1}}{n}\lesssim H^{2}n^{\frac{-(2\beta-1)}{\alpha+2\beta}}$.

{\noindent \bf \underline{Proof of Proposition \ref{prop:a1a2a3value}-(2)}.~}     We can reformulate  $B_{2}$ to  
\begin{align*}
	B_{2}=\left(\widehat{\Gamma}^{\dagger}_{m}-\Gamma^{\dagger}_{m}\right)\widetilde{\Gamma}_{e}&=\left(\widehat{\Gamma}^{\dagger}_{m}-\Gamma^{\dagger}_{m}\right)\Gamma\Gamma^{-1}\widetilde{\Gamma}_{e}.
\end{align*}

When $\ttE$ happens, it holds that $\mr{Im}(\Gamma_e)=\mr{Im}(\wt\Gamma_e)$. 
For any $\bs{u}\in\H$ with $\norm{\bs{u}}=1$, we can write 
\begin{align*}
	\Gamma^{-1}\widetilde{\Gamma}_{e}(\bs{u})=\sum^{d}_{k=1}c_{k}\eta_{k}
\end{align*}
with $|c_{k}|\leqslant \norm{\Gamma^{-1}\widetilde{\Gamma}_{e}}$ for all $k=1,\dots,d$ and $\{\eta_k\}_{k=1}^d$ are the generalized eigenfunctions  of $\Gamma_e$ associated with eigenvalues $\{\mu_k\}_{k=1}^d$ (i.e.,  $\Gamma_e\eta_k=\mu_k\Gamma\eta_k$). 

We conclude that 
\begin{align*}
\norm{B_{2}}^21_\ttE\leqslant d\sum^{d}_{k=1}\norm{\Gamma^{-1}\widetilde{\Gamma}_{e}}^2	1_\ttE\norm{(\widehat{\Gamma}^{\dagger}_{m}-\Gamma^{\dagger}_{m})\Gamma\eta_{k}}^21_\ttE.
\end{align*}
Now we need two lemmas:
\begin{lemma}
	\label{Lem:temp2}
	\begin{equation}\label{Eqt:temp}
		\mb E\left[\norm{(\widehat{\Gamma}^{\dagger}_{m}-\Gamma^{\dagger}_{m})\Gamma\eta_{k}}^{2}1_\ttE\right]\lesssim\frac{H^{2}m^{\alpha+1}}{n}, \qquad k=1,\dots,d.
\end{equation}\end{lemma}
\begin{proof} 
We first decompose $\left(\widehat{\Gamma}^{\dagger}_{m}-\Gamma^{\dagger}_{m}\right)\Gamma\eta_{k}$ as follows:
\begin{equation}\label{eq:nohggn_temp}
\begin{aligned}
\left(\widehat{\Gamma}^{\dagger}_{m}-\Gamma^{\dagger}_{m}\right)\Gamma\eta_{k}
			&=\widehat{\Gamma}^{\dagger}_{m}(\widehat{\Gamma}+(\Gamma-\widehat{\Gamma}))(\widehat{\eta}_{k}+(\eta_{k}-\widehat{\eta}_{k}))-\Gamma^{\dagger}_{m}\Gamma\eta_{k}\\
   &=\widehat{\Gamma}^{\dagger}_{m}\widehat{\Gamma}\widehat{\eta}_{k}-\Gamma^{\dagger}_{m}\Gamma\eta_{k}+\widehat{\Gamma}^{\dagger}_{m}\widehat{\Gamma}(\eta_{k}-\widehat{\eta}_{k})+\widehat{\Gamma}^{\dagger}_{m}(\Gamma-\widehat{\Gamma})\eta_k.    
\end{aligned}
\end{equation}
Suppose  $\eta_{k}=\sum_{j\geqslant 1}b_{kj}\phi_{j}$ and $\widehat{\eta}_{k}=\sum_{j\geqslant1}b_{kj}\widehat{\phi}_{j}$. 
Let $\eta^{(m)}_{k}=\sum^{m}_{j= 1}b_{kj}\phi_{j}$  and $\widehat{\eta}^{(m)}_{k}=\sum^{m}_{j=1}b_{kj}\widehat{\phi}_{j}$. Recall that we have introduced the notation $\Pi_m:=\sum_{i=1}^m\phi_i\otimes\phi_i$ and $\wh\Pi_m:=\sum_{i=1}^m\wh\phi_i\otimes\wh\phi_i$, then it holds that $\widehat{\Gamma}^{\dagger}_{m}\widehat{\Gamma}\widehat{\eta}_{k}=\wh\Pi_m\widehat{\eta}_{k}=\widehat{\eta}^{(m)}_{k}$ and similarly,  ${\Gamma}^{\dagger}_{m}{\Gamma}{\eta}_{k}=\Pi_m{\eta}_{k}={\eta}^{(m)}_{k}$.
In addition, $\widehat{\Gamma}^{\dagger}_{m}\widehat{\Gamma}\widehat{\eta}_{k}=\widehat{\Gamma}^{\dagger}_{m}\widehat{\Gamma}\widehat{\eta}_{k}^{(m)}$.
 Thus, 
 \begin{align*}
 \widehat{\Gamma}^{\dagger}_{m}\widehat{\Gamma}\widehat{\eta}_{k}-\Gamma^{\dagger}_{m}\Gamma\eta_{k}+\widehat{\Gamma}^{\dagger}_{m}\widehat{\Gamma}(\eta_{k}-\widehat{\eta}_{k})= \widehat{\eta}^{(m)}_{k}-\eta^{(m)}_{k}+  \widehat{\Gamma}^{\dagger}_{m}\widehat{\Gamma}(\eta_{k}-\widehat{\eta}_{k}^{(m)}). 
 \end{align*}
Insert this equality into \eqref{eq:nohggn_temp}, we have 
	\begin{equation}\label{eq:nohggn}
		\begin{aligned}
			\left(\widehat{\Gamma}^{\dagger}_{m}-\Gamma^{\dagger}_{m}\right)\Gamma\eta_{k}
			&=\widehat{\eta}^{(m)}_{k}-\eta^{(m)}_{k}-\widehat{\Gamma}^{\dagger}_{m}\widehat{\Gamma}(\widehat{\eta}^{(m)}_{k}-\eta^{(m)}_{k}+\eta^{(m)}_{k}-\eta_k)+\widehat{\Gamma}^{\dagger}_{m}(\Gamma-\widehat{\Gamma})\eta_{k}\\
   &=(\bs I-\Pi_m)(\widehat{\eta}^{(m)}_{k}-\eta^{(m)}_{k})-\widehat\Pi_{m}(\eta^{(m)}_{k}-\eta_k)+\widehat{\Gamma}^{\dagger}_{m}(\Gamma-\widehat{\Gamma})\eta_{k}
		\end{aligned}
	\end{equation}
 where $\bs I$ is the identity operator.

	We first find a bound for 
	\begin{equation*}
		\|(\bs I-\Pi_m)(\widehat{\eta}^{(m)}_{k}-\eta^{(m)}_{k})\|\leqslant\|\widehat{\eta}^{(m)}_{k}-\eta^{(m)}_{k}\|.
	\end{equation*}
	Note that 
\begin{equation}\label{eq:hatetaetadis}
	\widehat{\eta}^{(m)}_{k}-\eta^{(m)}_{k}=\sum^{m}_{j=1}b_{kj}(\widehat{\phi}_{j}-\phi_{j}).
\end{equation} 
It reduces to analyzing  $\norm{\widehat{\phi}_{j}-\phi_{j}}$.
Note that our predictor $\boldsymbol{X}$ satisfies the assumptions  in \citep{Hall2007mcflr}. 
By Equation (5.22) of \citep{Hall2007mcflr}, we have 
	\begin{equation}\label{eq:hpjpj}
		\E\left[\norm{\widehat{\phi}_{j}-\phi_{j}}^{2}\right]\lesssim j^{2}/n 
	\end{equation}
	uniformly in $1\leqslant j\leqslant m$.
	Substituting it into Equation \eqref{eq:hatetaetadis} and using Cauchy--Schwarz inequality, we obtain 
	\begin{equation*}
		\begin{aligned}	\bbE\left[\norm{\widehat{\eta}^{(m)}_{k}-\eta^{(m)}_{k}}^{2}\right]&=\bbE\left[\norm{\sum^{m}_{j=1}b_{kj}(\widehat{\phi}_{j}-\phi_{j})}^2\right]\\
			& \lesssim\frac{1}{n} \sum^{m}_{j=1}mj^{-2\beta}j^{2} \leqslant\frac mn\sum_{j=1}^\infty j^{2-2\beta}
   \lesssim \frac{m}{n}\lesssim \frac{m^{\alpha+1}}{n} 
		\end{aligned}
	\end{equation*}
by Assumption \ref{assumption: rate-type condition}.	
	Then we have 
	\begin{equation*}
		\bbE\left[\|(\bs I-\Pi_m)(\widehat{\eta}^{(m)}_{k}-\eta^{(m)}_{k})\|^2\right]\lesssim \frac{m^{\alpha+1}}{n}.
	\end{equation*}

Next, we bound $\|\widehat\Pi_{m}(\eta^{(m)}_{k}-\eta_k)\|$ as follows: 
 \begin{align*}
 \left\|\widehat\Pi_{m}(\eta^{(m)}_{k}-\eta_k)\right\|^2=&\left\|\widehat\Pi_{m}\sum_{j=m+1}^\infty b_{kj}\phi_j\right\|^2\leqslant \left\|\sum_{j=m+1}^\infty b_{kj}\phi_j\right\|^2\\
=&\sum_{j=m+1}^\infty b_{kj}^2\leqslant C\sum_{j=m+1}^\infty j^{-2\beta}\asymp m^{1-2\beta}\asymp \frac{m^{\alpha+1}}{n}
.  \end{align*}
by Assumption  \ref{assumption: rate-type condition} and the choice of $m$ in Theorem \ref{thm:cfsdrs}.

Finally, using a similar argument in the proof of Proposition \ref{prop:a1a2a3value} (1), we can easily derive the following inequality:
\begin{align*}
\bbE\left[\norm{\widehat{\Gamma}^{\dagger}_{m}\left(\widehat{\Gamma}-\Gamma\right)}^{2}1_\ttE\right]\lesssim \frac{H^{2} m^{\alpha+1}}{n}
\end{align*}
by noting that $\E\left[\left\|\widehat{\Gamma}-\Gamma\right\|^2\right]=O(n^{-1})$ (see e.g., Equation (5.9) in  \cite{Hall2007mcflr}).
	This proves the desired Equation \eqref{Eqt:temp}
	\end{proof}

\begin{lemma}\label{lem:Claim 1:}
 Suppose Assumption \ref{as:Linearity condition and Coverage condition} holds. 
There exists a constant $C$ that depends only on $\Gamma$ and $\Gamma_{e}$, such that if $\tau>\frac{6\norm{\Gamma_{e}}}{\lambda^{+}_{\min}(\Gamma_{e})}$ and the event $\ttE$ holds, then
	the norm of the operator $\Gamma^{-1}\widetilde{\Gamma}_{e}$ is no greater than $C$. 
\end{lemma}  
\begin{proof}

By Lemma \ref{lem:image equality}, if $\tau>\frac{6\norm{\Gamma_{e}}}{\lambda^{+}_{\min}(\Gamma_{e})}$, then we have $\mr{Im}(\Gamma_e)=\mr{Im}(\wt\Gamma_e)$ on the event $\ttE$. Furthermore, 
\begin{align*}
 \|\Gamma^{-1}\widetilde{\Gamma}_{e}\|1_\ttE \leqslant &
 \|\Gamma^{-1} \mid_{\S_{e}} \| 
 \|\widetilde{\Gamma}_{e}\|1_\ttE
 \leqslant \|\Gamma^{-1}\mid_{\S_{e}} \| \left(\|\widetilde{\Gamma}_{e}-\Gamma_{e}\|+\|\Gamma_{e}\| \right)1_\ttE\\
 \leqslant& \|\Gamma^{-1}\mid_{\S_{e}} \| (\frac{3}{\tau}+1) \|\Gamma_{e}\|
 \leqslant \|\Gamma^{-1}\mid_{\S_{e}} \| (\frac{\lambda_{\min}^+(\Gamma_e)}{2\|\Gamma_e\|}+1) \|\Gamma_{e}\|.
\end{align*}
The last expression is bounded since $\Gamma_e$ is of rank $d$ and 
 under Assumption \ref{as:Linearity condition and Coverage condition}, $\|\Gamma^{-1}\mid_{\S_{e}} \|$ is upper bounded since $\Gamma^{-1}\S_e=\S_{Y|\vX}$.

\end{proof}

Thanks to these two lemmas, we obtain $\bbE\left[\norm{B_{2}}^{2}\right]\lesssim\frac{H^{2}m^{\alpha+1}}{n}$, and finish the proof of Part-(2).

{\noindent \bf \underline{Proof of Proposition \ref{prop:a1a2a3value}-(3)}.~} 
For the term $B_{3}$, we have 
\begin{equation*}
	\widehat{\Gamma}^{\dagger}_{m}-\Gamma^{\dagger}_{m}=\sum^{m}_{j=1}(\widehat{\lambda}^{-1}_{j}-\lambda^{-1}_{j})\widehat{\phi}_{j}\otimes \widehat{\phi}_{j} +\sum^{m}_{j=1}\lambda^{-1}_{j}\left(\widehat{\phi}_{j}\otimes \widehat{\phi}_{j}-\phi_{j}\otimes \phi_{j}\right)=:A_{11}+A_{12}.
\end{equation*}

By Cauchy--Schwarz inequality, we have
\begin{align*}
	\bbE^2\left[\norm{A_{12}\left(\widehat{\Gamma}_{e}-\widetilde{\Gamma}_{e}\right)}1_\ttE\right]\leqslant&
	\bbE\left[\norm{A_{12}}^2\right] 
	\bbE\left[\norm{\left(\widehat{\Gamma}_{e}-\widetilde{\Gamma}_{e}\right)}^21_\ttE\right]  \lesssim \frac{1}{n} \bbE\left[\norm{A_{12}}^{2}\right]\\
 \lesssim&\frac 1n m\sum_{j=1}^m\lambda_j^{-2}\bbE[\|\widehat\phi_j-\phi_j\|^2].
\end{align*}

By Equation (5.22) in \cite{Hall2007mcflr},   we have $\bbE\left[\norm{\widehat{\phi}_{j}-\phi_{j}}^{2}\right]\lesssim\frac{j^{2}}{n}$ holds uniformly for all  $j\leqslant m$. Hence, by Assumption \ref{assumption: rate-type condition}, we have
\begin{equation*}
\frac 1n m\sum_{j=1}^m\lambda_j^{-2}\bbE[\|\widehat\phi_j-\phi_j\|^2]\lesssim\frac{m}{n^2}\sum_{j=1}^m j^{2+2\alpha}\asymp\frac{m^{4+2\alpha}}{n^2}\lesssim\frac{m^{\alpha+1}}{n}.
\end{equation*}
Also, we have 
\begin{equation*}
	A_{11}=\sum^{m}_{j=1}\frac{\widehat{\lambda}_{j}-\lambda_{j}}{\widehat{\lambda}_{j}\lambda_{j}}\left(\widehat{\phi}_{j}\otimes \widehat{\phi}_{j}-\phi_{j}\otimes \phi_{j}\right)+\sum^{m}_{j=1}\frac{\widehat{\lambda}_{j}-\lambda_{j}}{\widehat{\lambda}_{j}\lambda_{j}}\phi_{j}\otimes \phi_{j}=:A_{111}+_{A112}.
\end{equation*}
By  direct calculation, we find that 

\begin{align*}
	\norm{A_{111}\left(\widehat{\Gamma}_{e}-\widetilde{\Gamma}_{e}\right)}
	=&\norm{\sum^{m}_{j=1}\frac{\widehat{\lambda}_{j}-\lambda_{j}}{\widehat{\lambda}_{j}\lambda_{j}}\left(\widehat{\phi}_{j}\otimes \widehat{\phi}_{j}-\phi_{j}\otimes \phi_{j}\right)\left(\widehat{\Gamma}_{e}-\widetilde{\Gamma}_{e}\right)}\\
	\leqslant& C\sum^{m}_{j=1}|\frac{\widehat{\lambda}_{j}-\lambda_{j}}{\widehat{\lambda}_{j}\lambda_{j}}|\|\phi_j-\hat\phi_j\|\|\widehat{\Gamma}_{e}-\widetilde{\Gamma}_{e}\|.
		\end{align*}
  
Recall the definition of $\ttF$ in \eqref{eq:event F}, we have
\begin{align*}
	&\bbE^2\left[\norm{A_{111}\left(\widehat{\Gamma}_{e}-\widetilde{\Gamma}_{e}\right)}1_\ttE1_\ttF\right]
	\leqslant C^2m\sum^{m}_{j=1}\lambda_j^{-2}\bbE^2[\|\phi_j-\hat\phi_j\|\|\widehat{\Gamma}_{e}-\widetilde{\Gamma}_{e}\|]\\
 \leqslant& C^2m\sum^{m}_{j=1}\lambda_j^{-2}
 \bbE[\|\phi_j-\hat\phi_j\|^2]
 \bbE[\|\widehat{\Gamma}_{e}-\widetilde{\Gamma}_{e}\|^21_\ttE]\\
\leqslant &
 \frac{C^2m}{n}\sum^{m}_{j=1}\lambda_j^{-2}\bbE\left[\|\phi_j-\hat\phi_j\|^2\right]\lesssim\frac{m}{n^2}\sum_{j=1}^mj^{2+2\alpha}\lesssim\frac{m^{\alpha+1}}{n}.
		\end{align*} 
For the term $A_{112}$, we have 
	\begin{align*}
		\norm{A_{112}\left(\widehat{\Gamma}_{e}-\widetilde{\Gamma}_{e}\right)}=&\sup_{\bs{u}: \norm{\bs{u}}=1}\norm{\sum^{m}_{j=1}\frac{\widehat{\lambda}_{j}-\lambda_{j}}{\widehat{\lambda}_{j}\lambda_{j}}\phi_{j}\otimes \phi_{j}\left(\widehat{\Gamma}_{e}-\widetilde{\Gamma}_{e}\right)\bs{u}}\\
		=&\sup_{\bs{u}: \norm{\bs{u}}=1}\norm{\sum^{m}_{j=1}\frac{\widehat{\lambda}_{j}-\lambda_{j}}{\widehat{\lambda}_{j}\lambda_{j}}\left\langle \left(\widehat{\Gamma}_{e}-\widetilde{\Gamma}_{e}\right)\bs{u},  \phi_{j}   \right\rangle\phi_{j}	}\\
       =&\sup_{\bs{u}: \norm{\bs{u}}=1}\norm{\sum^{m}_{j=1}\frac{\widehat{\lambda}_{j}-\lambda_{j}}{\widehat{\lambda}_{j}\lambda_{j}}\sum_{i\geqslant 1} u_i\left\langle \left(\widehat{\Gamma}_{e}-\widetilde{\Gamma}_{e}\right)\phi_i,\phi_{j}\right\rangle\phi_{j}	}\quad (\bs{u}=\sum_{i\geqslant1}u_{i}\phi_i)\\
      \leqslant&  \left(\sum^{m}_{j=1}\left(\frac{\widehat{\lambda}_{j}-\lambda_{j}}{\widehat{\lambda}_{j}\lambda_{j}}\right)^{2}\sum_{i\geqslant 1} \left\langle\left(\widehat{\Gamma}_{e}-\widetilde{\Gamma}_{e}\right)\phi_{i},\phi_{j}\right\rangle^{2}\right)^{1/2}            
	\end{align*}
where in the fourth line, we use the  Cauchy-Schwarz inequality and the relation that $\sum_{i\geqslant 1} u_i^2=1$. Then we have
\begin{align*}
&\bbE\left[\norm{A_{112}\left(\widehat{\Gamma}_{e}-\widetilde{\Gamma}_{e}\right)}^2 1_{\ttE\cap\ttF}\right] \\ 
\leqslant & \bbE\left[\left(\sum^{m}_{j=1}\left(\frac{\widehat{\lambda}_{j}-\lambda_{j}}{\widehat{\lambda}_{j}\lambda_{j}}\right)^{2}\sum_{i\geqslant 1}\left\langle\left(\widehat{\Gamma}_{e}-\widetilde{\Gamma}_{e}\right)\phi_{i},\phi_{j}\right\rangle^{2}\right)1_{\ttE\cap\ttF}\right]\\
\leqslant&
\sum^{m}_{j=1}\lambda_j^{-2}\sum_{i\geqslant 1}\bbE\left[1_{\ttE}\left\langle\left(\widehat{\Gamma}_{e}-\widetilde{\Gamma}_{e}\right)\phi_{i},\phi_{j}\right\rangle^{2}\right]
\\
\leqslant&\sum^{m}_{j=1}\lambda_j^{-2}\sum_{i\geqslant 1}\frac{CH^{2}\lambda_{i}\lambda_{j}}{n}
\leqslant CH^{2} \sum^{m}_{j=1}\frac{1}{\lambda_{j}n}\sum_{i\geqslant 1}\lambda_{i}
\lesssim\frac{H^{2} m^{\alpha+1}}{n}.
\end{align*}
where the third inequality follows the same proof of Lemma \ref{lem:cvgem}. Thus, we complete the proof of Proposition \ref{prop:a1a2a3value}-(3). 

\qed
 
\subsection{The proof of Theorem \ref{thm:cfsdrs}}

Let $T:=\Gamma^{-1}\widetilde{\Gamma}_{e}\left(\Gamma^{-1}\widetilde{\Gamma}_{e}\right)^{*}$ and $\widehat{T}_{m}=\widehat{\Gamma}^{\dagger}_{m}\widehat{\Gamma}^{d}_{e}\widehat{\Gamma}^{d}_{e}\widehat{\Gamma}^{\dagger}_{m}$ where $\widehat{\Gamma}^{d}_{e}$ is defined in \eqref{eq: def of hatGammaed}. 
Define 
\begin{equation}\label{eq:event F prim}
 \ttF':=\left\{\frac{\lambda_d(T)}{2}\leqslant\lambda_d(\widehat T_m)\leqslant\frac{3\lambda_d(T)}{2}\right\}=\left\{|\lambda_d(\wh T_m)-\lambda_d(T)|\leqslant\frac{\lambda_d(T)}{2}\right\}.   
\end{equation}

To prove Theorem \ref{thm:cfsdrs}, we only need to prove $\mb{E}\left[\norm{P_{\widehat{\S}_{Y\mid \boldsymbol{X}}}-P_{\S_{Y\mid \boldsymbol{X}}}}^{1/2}1_{\ttE\cap\ttF\cap\ttF'}\right] \lesssim  n^{\frac{-(2\beta-1)}{4(\alpha+2\beta)}}$ and $\bbP(\ttE\cap\ttF\cap\ttF')\xrightarrow{n\to\infty}1$. 

Before we delve into the prove, we introduce some convenient notation. 
For any $\vu\in \H$, define $\vu^*: \H\to\R, \bs v\mapsto \langle \vu,\bs v\rangle$. Then $\vu^*$ is the adjoint operator of $\vu: \R\to \H,\lambda\mapsto \lambda\vu,(\forall\lambda\in\R)$ since   $\vu^*\bs v=\langle \vu^*\bs v,1\rangle=\langle\bs v,\vu\rangle$. 

Similarly, for any $d$ elements in $\H$, say $\bbeta_1,\ldots, \bbeta_d$, we can define $\vB:=(\bbeta_1,\dots,\bbeta_d):\R^d\to L^2[0,1]$ and its adjoint $\vB^{*}$. 

We also define the `truncated central space' 
\begin{align}\label{def: truncated central subspace}\mc S_{Y|\bs X}^{(m)}=\Pi_m \mc S_{Y|\bs X}=\mathrm{span}\{\bs{\beta}_1^{(m)},\dots,\bs{\beta}_d^{(m)}\},\end{align}
 where $\bbeta_{k}^{(m)}:=\Pi_m(\bbeta_k),k\in[d]$.
 For such a truncated central space, we have the following proposition, whose proof is deferred to the end.
\begin{proposition}\label{proposition, truncation error}
Under Assumption $\ref{assumption: rate-type condition}$, if $m$ is sufficiently large, we have
\begin{equation}\label{equation, truncation error}
 \left\|P_{\mathcal S_{Y|\boldsymbol{X}}}-P_{\mathcal S_{Y|\boldsymbol{X}}^{(m)}}\right\|\lesssim m^{-\frac{2\beta-1}{2}}.
\end{equation}
\end{proposition} 

Then we only need to show that 
\begin{align*}
\mb{E}\left[\norm{P_{\widehat{\S}_{Y\mid \boldsymbol{X}}}-P_{\S_{Y\mid \boldsymbol{X}}^{(m)}}}^{1/2}1_{\ttE\cap\ttF\cap\ttF'}\right]\lesssim  n^{\frac{-(2\beta-1)}{4(\alpha+2\beta)}}.
\end{align*}

Recall that $T:=\Gamma^{-1}\widetilde{\Gamma}_{e}\left(\Gamma^{-1}\widetilde{\Gamma}_{e}\right)^{*}$   and  $\widehat{T}_{m}=\widehat{\Gamma}^{\dagger}_{m}\widehat{\Gamma}^{d}_{e}\widehat{\Gamma}^{d}_{e}\widehat{\Gamma}^{\dagger}_{m}$. Let  $T_{m}:=\Pi_{m}T\Pi_{m}$.  We have 
  $P_{\S_{Y\mid \boldsymbol{X}}^{(m)}}=P_{T_m}$ and $T_{m}=\Gamma^{\dagger}_{m}\widetilde{\Gamma}_{e}\widetilde{\Gamma}_{e}\Gamma^{\dagger}_{m}$.   Furthermore, $P_{\widehat{\S}_{Y\mid \boldsymbol{X}}}=P_{\widehat{\Gamma}^{\dagger}_{m}\widehat{\Gamma}^{d}_{e}}=P_{\widehat{T}_{m}}$.
  
  By $\sin \Theta$ theorem (Lemma  \ref{lemma, sin theta of infinite dimension operator}), we have 
\begin{equation}\label{eq:sin theta estimation error}
\begin{aligned}
\norm{P_{\widehat{\S}_{Y\mid \boldsymbol{X}}}-P_{\S_{Y\mid \boldsymbol{X}}^{(m)}}}&=\norm{P_{T_{m}}-P_{\widehat{T}_m}}\leqslant \frac{\pi}{2}\frac{\norm{T_{m}-\widehat{T}_{m}}}{\min\{\lambda_{\min}^{+}(T_{m}),\lambda_{\min}^{+}(\widehat{T}_{m})\}}\\
&\leqslant \frac{\pi}{2}\frac{\left(\norm{\widehat{\Gamma}^{\dagger}_{m}\widehat{\Gamma}^{d}_{e}}+\norm{\Gamma^{\dagger}_{m}\widetilde{\Gamma}_{e}}\right)\norm{\widehat{\Gamma}^{\dagger}_{m}\widehat{\Gamma}^{d}_{e}-\Gamma^{\dagger}_{m}\widetilde{\Gamma}_{e}}}{\min\{\lambda_{\min}^{+}(T_{m}),\lambda_{\min}^{+}(\widehat{T}_{m})\}}.
\end{aligned}
\end{equation}

We first  provide an upper bound on the numerator of the right hand side of  \eqref{eq:sin theta estimation error}.  
Similar to the argument of Lemma \ref{lem:Claim 1:}, we find  that 
\begin{align}\label{eq:upper bound}
\bbE\left[\norm{\Gamma^{\dagger}_{m}\widetilde{\Gamma}_{e}}^21_{\ttE\cap\ttF}\right]\lesssim 1\quad\text{and}\quad \bbE\left[\norm{\widehat{\Gamma}^{\dagger}_{m}\widehat{\Gamma}^{d}_{e}}^21_{\ttE\cap\ttF}\right]\lesssim 1.      
\end{align}

We move on to analyze $\norm{\widehat{\Gamma}^{\dagger}_{m}\widehat{\Gamma}^{d}_{e}-\Gamma^{\dagger}_{m}\widetilde{\Gamma}_{e}}$.

Firstly, for any $\bs{u}\in\widehat{\S}_{e}$ with $\norm{\bs{u}}=1$, we have $\widehat{\Gamma}^{d}_{e}\u=\widehat{\Gamma}_{e}\u$. 
Then by  Theorem \ref{thm:bt1t2}, we have  
\begin{equation}\label{eq:normonse}
\begin{aligned}
&\bbE\left[1_{\ttE\cap\ttF}\norm{\left(\widehat{\Gamma}^{\dagger}_{m}\widehat{\Gamma}^{d}_{e}-\Gamma^{\dagger}_{m}\widetilde{\Gamma}_{e}\right)(\bs{u})}\right]=\bbE\left[1_{\ttE\cap\ttF}\norm{\left(\widehat{\Gamma}^{\dagger}_{m}\widehat{\Gamma}_{e}-\Gamma^{\dagger}_{m}\widetilde{\Gamma}_{e}\right)(\bs{u})}\right]\\
\leqslant&\bbE\left[1_{\ttE\cap\ttF}\norm{\left(\widehat{\Gamma}^{\dagger}_{m}\widehat{\Gamma}_{e}-\Gamma^{\dagger}_{m}\widetilde{\Gamma}_{e}\right)}\right]\lesssim Hn^{\frac{-(2\beta-1)}{2(\alpha+2\beta)}}.
\end{aligned}
\end{equation} 
This shows that $\bbE\left[\norm{\left(\widehat{\Gamma}^{\dagger}_{m}\widehat{\Gamma}^{d}_{e}-\Gamma^{\dagger}_{m}\widetilde{\Gamma}_{e}\right)\mid_{\widehat{\S}_{e}}}1_{\ttE\cap\ttF}\right]\lesssim Hn^{\frac{-(2\beta-1)}{2(\alpha+2\beta)}}.$

Secondly, for any $\bs{u}\in\widehat{\S}_{e}^{\perp}$ with $\norm{\bs{u}}=1$, 
\begin{equation}\label{eq:normonseverp}
\begin{aligned}
   	\norm{\left(\widehat{\Gamma}^{\dagger}_{m}\widehat{\Gamma}^{d}_{e}-\Gamma^{\dagger}_{m}\widetilde{\Gamma}_{e}\right)(\bs{u})}&=	\norm{\Gamma^{\dagger}_{m}\widetilde{\Gamma}_{e}(\bs{u})} 
  =\norm{\Gamma_{m}\Gamma^{-1}\widetilde{\Gamma}_{e}P_{\S_{e}}(\bs{u})}\\
   	&=\norm{\Gamma_{m}\Gamma^{-1}\widetilde{\Gamma}_{e}(P_{\S_{e}}-P_{\widehat{\S}_{e}})\bs{u}}
   	\leqslant \norm{\Gamma_{m}\Gamma^{-1}\widetilde{\Gamma}_{e}(P_{\S_{e}}-P_{\widehat{\S}_{e}})}\\
   	&\leqslant \norm{\Gamma_{m}}\norm{\Gamma^{-1}\widetilde{\Gamma}_{e}}\norm{P_{\S_{e}}-P_{\widehat{\S}_{e}}}.
   	\end{aligned}
\end{equation}
Then by Lemma \ref{lem:Claim 1:} and Theorem \ref{thm:rootnconsistency}, we have
\begin{align*}
\bbE\left[\norm{\left(\widehat{\Gamma}^{\dagger}_{m}\widehat{\Gamma}^{d}_{e}-\Gamma^{\dagger}_{m}\widetilde{\Gamma}_{e}\right)(\bs{u})}1_{\ttE\cap\ttF}\right]\lesssim\sqrt{\frac{H^{2}}{n}}.
\end{align*}
It implies that $\bbE\left[1_{\ttE\cap\ttF}\norm{\left(\widehat{\Gamma}^{\dagger}_{m}\widehat{\Gamma}^{d}_{e}-\Gamma^{\dagger}_{m}\widetilde{\Gamma}_{e}\right)\mid_{\widehat{\S}_{e}^{\perp}}}\right]\lesssim\sqrt{\frac{H^{2}}{n}}$.
Combing Equations \eqref{eq:normonse} and \eqref{eq:normonseverp}, we obtain 
\begin{equation}\label{eq:nomrestiatier}
\bbE\left[1_{\ttE\cap\ttF}\norm{\widehat{\Gamma}^{\dagger}_{m}\widehat{\Gamma}^{d}_{e}-\Gamma^{\dagger}_{m}\widetilde{\Gamma}_{e}}\right]\lesssim\left(Hn^{\frac{-(2\beta-1)}{2(\alpha+2\beta)}}\right).
\end{equation}

Lastly, we provide a lower bound on the denominator of the right hand side of \eqref{eq:sin theta estimation error}. 
If $\tau>\frac{6\norm{\Gamma_{e}}}{\lambda^{+}_{\min}(\Gamma_{e})}$, then $\mr{Im}(\wt\Gamma_e)=\mr{Im}(\Gamma_e)$ on the event $\ttE$ . By Lemma \ref{lem:PimTPimtoT}, one has 
\begin{align}\label{eq:T to Tm}
\|\Gamma^{-1}\wt\Gamma_{e}-\Gamma^{\dagger}_{m}\widetilde{\Gamma}_{e}\|\xrightarrow{m\to\infty}0.
\end{align}
By Lemma \ref{lem:wely ineq operator}, one has
$\sigma^{+}_{\min}(\Gamma^{\dagger}_{m}\widetilde{\Gamma}_{e})\geqslant \frac{\sigma^{+}_{\min}(\Gamma^{-1}\Gamma_{e})}{2}$  for sufficiently large $m$,  where $\sigma^{+}_{\min}$ denotes the  infimum of the positive singular values.

By Markov inequality, we have 
\begin{equation}\label{eq: EFFprimec prop}
 \begin{aligned}
&\bbP(1_{\ttE\cap\ttF}|\lambda_d(\widehat T_m)-\lambda_d(T)|\geqslant\frac{\lambda_d(T)}{2})\leqslant \frac{\bbE[1_{\ttE\cap\ttF}|\lambda_d(\widehat T_m)-\lambda_d(T)|]}{\lambda_d(T)/2}
\leqslant \frac{\bbE[1_{\ttE\cap\ttF}\|\widehat T_m-T\|]}{\lambda_d(T)/2}\\
\lesssim& 
\frac{\bbE[1_{\ttE\cap\ttF}(\|\wh\Gamma^{\dagger}_{m}\widehat{\Gamma}_{e}^d-\Gamma^{-1}\wt\Gamma_e\|^2+2\|\wh\Gamma^{\dagger}_{m}\widehat{\Gamma}_{e}^d-\Gamma^{-1}\wt\Gamma_e\|\|\Gamma^{-1}\wt\Gamma_e\|)]}{\lambda_d(T)/2}\\
\lesssim&\frac{\bbE[1_{\ttE\cap\ttF}\|\wh\Gamma^{\dagger}_{m}\widehat{\Gamma}_{e}^d-\Gamma^{-1}\wt\Gamma_e\|]}{\lambda_d(T)/2}\xrightarrow{m\to\infty}0.
\end{aligned}   
\end{equation}
where the second inequality comes from Lemma \ref{lem:wely ineq operator}, the fourth comes from 
Lemma \ref{lem:Claim 1:} and the fifth comes from
\eqref{eq:nomrestiatier}-\eqref{eq:T to Tm}.

Recall that 
\begin{equation*}
 \ttF':=\left\{\frac{\lambda_d(T)}{2}\leqslant\lambda_d(\widehat T_m)\leqslant\frac{3\lambda_d(T)}{2}\right\}=\left\{|\lambda_d(\wh T_m)-\lambda_d(T)|\leqslant\frac{\lambda_d(T)}{2}\right\}.   
\end{equation*}
Thus \eqref{eq: EFFprimec prop} implies
\begin{align}\label{eq:control on Fprime}
\bbP( (\ttF')^c\cap (\ttE\cap \ttF))\xrightarrow{m\to\infty}0.
\end{align}
Furthermore,
\begin{align*}
\bbP(\ttE\cap\ttF\cap\ttF')=
\bbP(\ttE\cap\ttF)-\bbP(\ttE\cap\ttF\cap(\ttF')^c)
   \xrightarrow{m\to\infty}1.
\end{align*}

On the event $\ttE\cap\ttF\cap\ttF'$, it holds that  $\sigma^{+}_{\min}(\widehat{\Gamma}^{\dagger}_{m}\widehat{\Gamma}^{d}_{e})\geqslant \frac{\sigma^{+}_{\min}(\Gamma^{-1}\Gamma_{e})}{4}$, which implies 
\begin{align}\label{eq:min eigenvalues lower bound}
\min\{\lambda_{\min}^{+}(T_{m}),\lambda_{\min}^{+}(\widehat{T}_{m})\}\geqslant C
\end{align}
for some constant $C>0$.

Inserting  \eqref{eq:upper bound}, \eqref{eq:nomrestiatier} and \eqref{eq:min eigenvalues lower bound} into \eqref{eq:sin theta estimation error}, we have
\begin{align*}
&\mb{E}^2\left[\norm{P_{\widehat{\S}_{Y\mid \boldsymbol{X}}}-P_{\S_{Y\mid \boldsymbol{X}}^{(m)}}}^{1/2}1_{\ttE\cap\ttF\cap\ttF'}\right]=\bbE^2\left[\norm{P_{T_{m}}-P_{\widehat{T}_{m}}}^{1/2}1_{\ttE\cap\ttF\cap\ttF'}\right] \\
\lesssim&\frac{\bbE^2\left[\left(\norm{\widehat{\Gamma}^{\dagger}_{m}\widehat{\Gamma}^{d}_{e}}+\norm{\Gamma^{\dagger}_{m}\widetilde{\Gamma}_{e}}\right)^{1/2}\norm{\widehat{\Gamma}^{\dagger}_{m}\widehat{\Gamma}^{d}_{e}-\Gamma^{\dagger}_{m}\widetilde{\Gamma}_{e}}^{1/2}1_{\ttE\cap\ttF\cap\ttF'}\right]}{\min\{\lambda_{\min}^{+}(T_{m}),\lambda_{\min}^{+}(\widehat{T}_{m})\}1_{\ttE\cap\ttF\cap\ttF'}}\\
\lesssim &\bbE\left[\left(\norm{\widehat{\Gamma}^{\dagger}_{m}\widehat{\Gamma}^{d}_{e}}+\norm{\Gamma^{\dagger}_{m}\widetilde{\Gamma}_{e}}\right)1_{\ttE\cap\ttF\cap\ttF'}\right]\bbE\left[\norm{\widehat{\Gamma}^{\dagger}_{m}\widehat{\Gamma}^{d}_{e}-\Gamma^{\dagger}_{m}\widetilde{\Gamma}_{e}}1_{\ttE\cap\ttF\cap\ttF'}\right]
\lesssim Hn^{\frac{-(2\beta-1)}{2(\alpha+2\beta)}}.
\end{align*}
This completes the proof of  Theorem \ref{thm:cfsdrs}.

\begin{proof}[Proof of Proposition~\ref{proposition, truncation error}]
Let 
$\mc{B}:=\sum\limits_{i=1}^d \bbeta_i\otimes\bbeta_i$ and ${\mc{B}^{(m)}}:=\sum\limits_{i=1}^d \bbeta_i^{(m)}\otimes\bbeta_i^{(m)}$. 
Note that $\mr{Im}(\mc{B})=\mr{span}\{\bbeta_1,\dots,\bbeta_d\}=\mathcal S_{Y|\boldsymbol{X}}$. Similarly, $\mr{Im}(\mc{B}^{(m)})=\mr{span}\{\bbeta_1^{(m)},\dots,\bbeta_d^{(m)}\}=\mathcal S_{Y|\boldsymbol{X}}^{(m)}$. Thus $\left\|P_{\mathcal S_{Y|\boldsymbol{X}}}-P_{\mathcal S_{Y|\boldsymbol{X}}^{(m)}}\right\|=\|P_{\mc{B}}-P_{\mc{B}^{(m)}}\|$.

By Lemma \ref{lemma, sin theta of infinite dimension operator}, we have
\begin{align}\label{eq:sin theta for B Bm}
\|P_{\mc{B}}-P_{\mc{B}^{(m)}}\|\leqslant \frac{\pi}{2}\frac{\|{\mc{B}}-{\mc{B}^{(m)}}\|}{\min\{\lambda_{\min}^+({\mc{B}}),\lambda_{\min}^+({\mc{B}^{(m)}})\}}.
\end{align}

Note that ${\mc{B}}-{\mc{B}^{(m)}}$ is self-adjoint, then
\begin{align*}
&\lno{\mc{B}}-{\mc{B}^{(m)}}\rno=\sup_{{\bs{u}}\in\mathbb{S}_{ \mathcal H}}|\langle ({\mc{B}}-{\mc{B}^{(m)}})({\bs{u}}),{\bs{u}}\rangle|=\sup_{{\bs{u}}\in\mathbb{S}_{\mathcal H}}|\langle {\mc{B}}{\bs{u}},{\bs{u}}\rangle-\langle {\mc{B}^{(m)}}{\bs{u}},{\bs{u}}\rangle|\\
&~~=\sup_{{\bs{u}}\in\mathbb{S}_{\mathcal H}}\hspace{-0.9mm}\left|\sum_{i=1}^d\hspace{-0.9mm}\left[\langle{\bs{\beta}}_i,{\bs{u}}\rangle^2-\langle{\bs{\beta}}_i^{(m)},{\bs{u}}\rangle^2\right]\right|=\sup_{{\bs{u}}\in\mathbb{S}_{\mathcal H}}\hspace{-0.9mm}\left| \sum_{i=1}^d\langle{\bs{\beta}}_i-{\bs{\beta}}_i^{(m)},{\bs{u}}\rangle\langle{\bs{\beta}}_i+{\bs{\beta}}_i^{(m)},{\bs{u}}\rangle\right|\\
&~~\leqslant\sup_{{\bs{u}}\in\mathbb{S}_{\mathcal H}}\sum_{i=1}^d\left| \langle{\bs{\beta}}_i-{\bs{\beta}}_i^{(m)},{\bs{u}}\rangle\langle{\bs{\beta}}_i+{\bs{\beta}}_i^{(m)},{\bs{u}}\rangle\right|
\leqslant\sum_{i=1}^d\left\|{\bs{\beta}}_i-{\bs{\beta}}_i^{(m)}\right\|\left\|{\bs{\beta}}_i+{\bs{\beta}}_i^{(m)}\right\|,
\end{align*}
where the first inequality comes from the triangle inequality, and the 
second inequality comes from the Cauchy-Schwarz inequality and $\|{\bs{u}}\|=1$. 
According to Assumption $\ref{assumption: rate-type condition}$, one can get
\begin{align*}
\left\|{\bs{\beta}}_i-{\bs{\beta}}_i^{(m)}\right\|&=\left\|\sum_{j=m+1}^\infty b_{ij}\phi_j\right\|=\sqrt{\sum_{j=m+1}^\infty b_{ij}^2}\lesssim \sqrt{\sum_{j=m+1}^\infty j^{-2\beta}};\\
\left\|{\bs{\beta}}_i+{\bs{\beta}}_i^{(m)}\right\|&\leqslant\|{\bs{\beta}}_i\|+\lno{\bs{\beta}}_i^{(m)}\rno\leqslant2\|{\bs{\beta}}_i\|=2\sqrt{\sum_{j=1}^\infty b_{ij}^2}\lesssim\sqrt{\sum_{j=1}^\infty j^{- 2\beta}}.
\end{align*}
Because $\beta>1/2$, one has
\[\sum\limits_{j=m+1}^\infty \frac{1}{j^{2\beta}}\lesssim \frac{1}{m^{2\beta-1}};\qquad \sum_{j=1}^\infty \frac 1{j^{2\beta}}<\infty.\]
 Thus, one can get
\begin{equation}\label{eq: upper bound of operator norm of A minus B}
\lno{\mc{B}}-{\mc{B}^{(m)}}\rno\lesssim m^{-\frac{2\beta-1}{2}}.
\end{equation}

Then we show that $\min\{\lambda_{\min}^+({\mc{B}}),\lambda_{\min}^+({\mc{B}^{(m)}})\}\geqslant C$ for some constant $C>0$. 
Since $\mr{rank}(\mc{B})=d$, one can get that $\lambda_{\min}^+(\mc{B})=\lambda_{d}(\mc{B})$. It is easy to see $\mr{rank}(\mc{B}^{(m)})\leqslant d$ by $\mc{B}^{(m)}=\Pi_m \mc{B} \Pi_m$, thus one can assume that $\lambda_{\min}^+(\mc{B}^{(m)})=\lambda_j( \mc{B}^{(m)})$ for some $j\leqslant d$.
By Lemma $\ref{lem:wely ineq operator}$
and \eqref{eq: upper bound of operator norm of A minus B}, one has:
$$
|\lambda_j( \mc{B}^{(m)})-\lambda_j\l \mc{B}\r|\leqslant\lno \mc{B}-\mc{B}^{(m)}\rno\lesssim m^{-\frac{2\beta-1}{2}}.
$$
Thus for sufficiently large $m$, one has
\begin{align}
&\lambda_j\l \mc{B}^{(m)}\r\geqslant \frac{\lambda_d\l \mc{B}\r}{2}\Longrightarrow \min\{\lambda_{\min}^+({\mc{B}}),\lambda_{\min}^+({\mc{B}^{(m)}})\}\geqslant \frac{\lambda_d({\mc{B}})}{2}. \label{eq:lower bound lambda min plus B Bm}
\end{align}
Inserting \eqref{eq: upper bound of operator norm of A minus B} and \eqref{eq:lower bound lambda min plus B Bm} into \eqref{eq:sin theta for B Bm} leads to
\begin{align*}
\left\|P_{\mathcal S_{Y|\boldsymbol{X}}}-P_{\mathcal S_{Y|\boldsymbol{X}}^{(m)}}\right\|\lesssim m^{-\frac{2\beta-1}{2}}.
\end{align*}
Thus we complete the proof of Proposition \ref{proposition, truncation error}.
\end{proof}

\section{Proof of Theorem \ref{thm:opmmultiple}}
\subsection{Proof outline}\label{sec:proof outline}
We sketch the proof outline in this subsection and defer the proof details to the subsequent subsections.

We follow the the standard procedure of applying Fano's inequality to obtain the minimax lower bound. 
The following lemma is one version of the generalized Fano method.
\begin{lemma}[\cite{yu1997assouad}]\label{lem:fano}
	Let $N\geqslant 2$ be an integer and $\{\btheta_1,\dots, \btheta_N\}\subset \Theta_{0}$ index a collection of probability measures $\bbP_{\btheta_i}$ on  a measurable space $(\mc X, \mathcal{A})$. Let $\rho$ be a pseudometric on $\Theta_{0}$ and suppose that for all $i \neq j$
	$$
	\rho(\btheta_i, \btheta_j)\geqslant \alpha_N,\quad\text{and}\quad \kl(\bbP_{\btheta_i}, \bbP_{\btheta_j}) \leqslant \beta_N. 
	$$
	Then every $\mathcal{A}$-measurable estimator $\hat{\btheta}$ satisfies 
	$$
	\max_{i} \bbP \left( \rho(\hat{\btheta}, \btheta_{i}) \geqslant \frac{\alpha_N}{2}  \right) \geqslant 1-\frac{\beta_N +\log 2}{\log N} . 
	$$
\end{lemma}
To apply Lemma~\ref{lem:fano}, we need to construct a family of distributions that are separated from each other in the parameter space but close to each other in terms of the KL-divergence. 

Let us first recall the following  Varshamov--Gilbert bound \cite[Lemma~2.9]{Tsybakov:1315296}. 
\begin{lemma}\label{lem:gbound}
    For any $m>8$, there exists a set $\Theta:=\{\btheta^{(0)},\ldots,\btheta^{(N)}\}\subset \{-1,1\}^{m}$, such that 
    \begin{itemize}
        \item[1).] $\btheta^{(0)}=(-1,\ldots,-1)$;
        \item[2).] for any $\btheta,\btheta'\in\Theta$ and $\btheta\neq \btheta^{\prime}$, $\|\btheta-\btheta'\|^2\geqslant m/2$;
        \item[3).] $N\geqslant 2^{m/8}$.
    \end{itemize}
\end{lemma}
Let $\phi_{1}(t)=1, \phi_{j+1}(t)=\sqrt{2}\cos(j\pi t), j\geqslant 1$. For any $\btheta=(\theta_i)_{i\in[m]}\in\Theta$ in Lemma \ref{lem:gbound} and any  $\beta>3/2$, let us define the central space $\S(\btheta):=\mr{span}\{\bbeta^{\btheta}_{1},\ldots,\bbeta^{\btheta}_{d}\}$ as follows:
\begin{align*}
\bbeta^{\btheta}_{i}:=\sum^{(i+1)m}_{k=im+1}\theta_{k-im}k^{-\beta}\phi_{k}+\phi_i
\end{align*}
where $m=\wt C n^{\frac{1}{\alpha+2\beta}} $ for some $\wt C=\wt C(\alpha,\beta)$ to be determined. We assume $d$ is fixed and $m>d$. 
The collection $\{\S(\btheta):\btheta\in \Theta \}$ satisfies a nice property  stated in the following lemma. 

\begin{lemma}\label{lem:packing set}
For any different $\btheta $ and $\btheta^{\prime}$ in $\Theta$, we have
\begin{align*}
\|P_{\S(\btheta)}-P_{\S(\btheta^{\prime})}\|^2\geqslant  4\vartheta n^{-\frac{2\beta-1}{\alpha+2\beta}}
\end{align*}  
for some constant $\vartheta>0$ that  depends on $\alpha$ and $\beta$.
\end{lemma}

We next move on to construct a population corresponds to each $\S(\btheta)$. 
For this purpose, we will make use of a construction in \cite{lin2021optimalitysupplement}, described as follows. 
For any $x\in\R$, let $\phi(x)$ be a smooth function which maps $(-\infty, 0]$ to 0 and $[1, \infty)$ to 1 and has a positive first derivative over $(0,1)$. For any   $\bs z=(z_i)_{i=1}^d\in\R^d$, let $f(\bs z):=\sum_{i \leqslant d} 2^{i-1} \phi\left(z_i /\zeta \right)$ where $\zeta$ is sufficiently small such that for $\vZ=(Z_i)_{i\in[d]}\sim N(0,\bs I_d)$,  the probability $\mathbb{P}\left(\exists i, 0<Z_i<\zeta\right)\leqslant d\zeta/(\sqrt{2\pi}) <2^{-d}$. 
Let $Y=A f\left(\bs Z\right)+\varepsilon$ for some positive constant $A$ and $\varepsilon\sim N(0,1)$.
If $A$ is sufficiently large, the distribution of $(\vZ, Y)$ satisfies the coverage condition  \citep[Lemma 15]{lin2021optimalitysupplement}. 
Note that the joint distribution of $(\vZ, Y)$ is Lebesgue continuous and it is easy to check that  $\mb E[\vZ|Y=y]$ is a continuous function with respect to $y$ using the formula for conditional expectation. Thus by Lemma \ref{lem:monment condition to sliced stable}, we know that $(\vZ,Y)$ satisfied WSSC.

Now we  describe how to construct a distribution $\bbP_{\btheta}$ of $(\vX,Y)$ for any given $\btheta\in \Theta$. 
For some $\alpha$ satisfying 
$\alpha>1$ and  $\frac12\alpha+1<\beta$, let  $\boldsymbol{X}=\sum^{\infty}_{j=1}j^{-\alpha/2}X_{j}\phi_{j}$ such that $X_{j}\overset{iid}{\sim} N(0,1),j\geqslant 1$. Then $\Gamma=\sum_{j=1}^\infty j^{-\alpha}\phi_j\otimes\phi_j$. 
For any $\bbeta_i(i\in[d])$, we construct the joint distribution of $(\vX,Y)$ as follows: 
\begin{align}\label{eq: XY distribution for lower bound}
Y=Af((\vB^*\Gamma\vB)^{-1/2}\vB^*\vX)+\varepsilon,\quad \boldsymbol{X}=\sum^{\infty}_{j=1}j^{-\alpha/2}X_{j}\phi_{j},\quad \epsilon \sim N(0,1), 
\end{align}
where $\vB:=(\bbeta_1,\dots,\bbeta_d):\R^d\to L^2[0,1]$.

For this distribution of $(\vX, Y)$, we can prove that it belongs to the distribution class $\mathfrak{M}\left(\alpha,\beta,\tau\right)$. 
\begin{lemma}\label{lem: XY in distribution class}
For $(\vX,Y)$ constructed in \eqref{eq: XY distribution for lower bound}, it holds that
\begin{itemize}    \item[\textbf{i)}]$\|\Gamma\|\leqslant C$, $\lambda_{\min}(\Gamma |_{\mathcal{S}_{e}})\geqslant c$ and  $c\leqslant\lambda_d(\Gamma_e)\leqslant\lambda_1(\Gamma_e)\leqslant C$ for two positive constants $c$ and $C$ that do not depend on $\btheta$, $n$, and $m$;
    \item [\textbf{ii)}] the central curve $\vm(y) = \mb E[\vX|Y = y]$ is weak sliced stable with respect to $Y$;
    \item[\textbf{iii)}] $(\vX,Y)\in\mathcal{F}(\alpha,\beta,c_{1},c_{2})$. 
\end{itemize}
\end{lemma}

 Furthermore, we have the following upper bound on the pairwise KL-divergence:
\begin{lemma}\label{lem:KL divergence close}
Let $\mb{P}_{\btheta}$ and $\mb{P}_{\btheta'}$ be the joint distributions of $(\vX,Y)$ induced by $\bbeta_i^{\btheta}(i\in[d])$ and $\bbeta_i^{\btheta'}(i\in[d])$ respectively. Then we have 
    \begin{equation*}
	\mr{KL}(\mb{P}_{\btheta}, \mb{P}_{\btheta'})\lesssim m^{-\alpha}\sum^{d}_{i=1}\|\bbeta^{\theta}_{i}-\bbeta^{\theta^{\prime}}_{i}\|^2.
    \end{equation*}
\end{lemma}

Since
\begin{align*}
\sum^{d}_{i=1}\|\bbeta^{\theta}_{i}-\bbeta^{\theta^{\prime}}_{i}\|^2\lesssim  &m^{-2\beta}\|\btheta-\btheta'\|^2,
\end{align*}
we have
\begin{align*}
\mr{KL}(\mb{P}_{\btheta}, \mb{P}_{\btheta'})\leqslant C(d)m^{-(\alpha+2\beta)}\sum^{d}_{i=1}\|\bbeta^{\theta}_{i}-\bbeta^{\theta^{\prime}}_{i}\|^2\leqslant C(d)m^{-(\alpha+2\beta-1)}.
\end{align*}

We can now prove the minimax lower bound using Lemma \ref{lem:fano}: 

\begin{align*}
 &\inf_{\wh{\mathcal{S}}_{Y\mid \boldsymbol{X}}}\sup_{\mc M\in \mathfrak{M}\left(\alpha,\beta,\tau\right)}\bbP_{\mc M}\left(\norm{P_{\wh{\mathcal{S}}_{Y\mid \boldsymbol{X}}}-P_{\mathcal{S}_{Y\mid \boldsymbol{X}}}}^{2}\geqslant \vartheta  n^{-\frac{2\beta-1}{\alpha+2\beta}}\right)\\
 \geqslant&\inf_{\wh\S_{Y|\vX}}  \sup_{\btheta \in \Theta}\bbP_{\btheta}\left(\|P_{\wh\S_{Y|\vX}}-P_{\S(\btheta)}\|^2\geqslant \vartheta n^{\frac{-2\beta+1}{\alpha+2\beta}} \right)\\
	\geqslant& 1-\frac{\max \kl(\bbP_{\btheta}^{n}, \bbP_{\btheta'}^{n}) +\log(2)}{\log(|\Theta|)}  \\
	\geqslant  &1-\frac{ C(d)n m^{-(\alpha+2\beta-1)}+\log2}{\frac m2\log(2)}\\
=&1-\frac{ C(d)\tilde C^{-(\alpha+2\beta)}m+\log(2)}{\frac m2\log(2)}\geqslant 0.9
\end{align*}
where the constant $\vartheta$ comes from Lemma~\ref{lem:packing set} and in the last equation we have chosen $m=\wt C n^{\frac{1}{\alpha+2\beta}}$ for some $\wt C\geqslant \left(\frac{\log(2)}{30C(d)}\right)^{-\frac{1}{\alpha+2\beta}}$. 
This finishes the proof of Theorem \ref{thm:opmmultiple}.

\subsection{Proof of Lemma \ref{lem:packing set}}
\begin{proof}
    Note that  $P_{\S(\btheta)}=\sum_{i=1}^{d}\widetilde{\bbeta}^{\btheta}_{i}\otimes \widetilde{\bbeta}^{\btheta}_{i}$, where $\widetilde{\bbeta}^{\btheta}_{i}$ is the normalization of $\bbeta^{\btheta}_{i}$.  By the definition of $\bbeta^{\btheta}_{i}$, for any $\btheta$ and $\btheta^{\prime}$, $\langle \bbeta^{\btheta}_{i},  \bbeta^{\btheta^{\prime}}_{j} \rangle=0$  if $i\neq j$.  We have 
   
    \begin{align*}
      & \left\|P_{\S(\btheta)}-P_{\S(\btheta^{\prime})}\right\|^2\geqslant\left\|(P_{\S(\btheta)}-P_{\S(\btheta^{\prime})})(\widetilde{\bbeta}^{\btheta}_{1})\right\|^{2}=\left\|\wt{\bbeta}_1^{\btheta}-\langle \wt{\bbeta}_1^{\btheta},\wt{\bbeta}_1^{\btheta'}\rangle\wt{\bbeta}_1^{\btheta'}\right\|^2\\
          =   &\frac{1}{\|\bbeta_1^{\btheta}\|^2}\left\| \sum_{k=m+1}^{2m}(\theta_{k-m}-\langle\wt\bbeta_1^{\btheta},\wt\bbeta_1^{\btheta'}\rangle\theta'_{k-m})k^{-\beta}\phi_k\right\|^2\\
        =&\frac{1}{\|\bbeta_1^{\btheta}\|^2}\left\| \sum_{\substack{k\in\{m+1,\dots,2m\}\\\theta_{k-m}\neq\theta'_{k-m}}}(\theta_{k-m}-\langle \wt\bbeta_1^{\btheta},\wt\bbeta_1^{\btheta'}\rangle\theta'_{k-m})k^{-\beta}\phi_k+\sum_{\substack{k\in\{m+1,\dots,2m\}\\\theta_{k-m}=\theta'_{k-m}}}(\theta_{k-m}-\langle \wt\bbeta_1^{\btheta},\wt\bbeta_1^{\btheta'}\rangle\theta'_{k-m})k^{-\beta}\phi_k\right\|^2\\    
        \geqslant&\frac{1}{\|\bbeta_1^{\btheta}\|^2}\left \| \sum_{\substack{k\in\{m+1,\dots,2m\}\\\theta_{k-m}\neq\theta'_{k-m}}}(\theta_{k-m}-\langle \wt\bbeta_1^{\btheta},\wt\bbeta_1^{\btheta'}\rangle\theta'_{k-m})k^{-\beta}\phi_k\right\|^2\\   
        \geqslant& \frac{1}{\|\bbeta_1^{\btheta}\|^2}\left\| \sum_{\substack{k\in\{m+1,\dots,2m\}\\\theta_{k-m}\neq\theta'_{k-m}}}k^{-\beta}\phi_k\right\|^2=\frac{1}{\|\bbeta_1^{\btheta}\|^2}\sum_{\substack{k\in\{m+1,\dots,2m\}\\\theta_{k-m}\neq\theta'_{k-m}}}k^{-2\beta}.
    \end{align*}
From the property (2) of Lemma \ref{lem:gbound}, we know that there are at least $m/8$ $k$'s satisfying $k\in\{m+1,\dots,2m\},\theta_{k-m}\neq\theta'_{k-m}$. Thus 
\begin{align*}
\frac{1}{\|\bbeta_1^{\btheta}\|^2}\sum_{\substack{k\in\{m+1,\dots,2m\}\\\theta_{k-m}\neq\theta'_{k-m}}}k^{-2\beta}\geqslant 
\frac{1}{(1+\sum_{k=m+1}^{2m}k^{-2\beta})^2}\sum_{k=15m/8}^{2m}k^{-2\beta}\asymp m^{-2\beta+1}.
\end{align*}

\end{proof}
\subsection{Proof of Lemma \ref{lem: XY in distribution class}}
\begin{proof}
\textbf{Proof of \textbf{i)}}~\\
It is easy to check that $\|\Gamma\|=1$. 
Now we give a lower bound of $\lambda_{\min}(\Gamma|_{\mc S_e})$. For any unit function $\vu\subseteq\mr{Im}(\Gamma_e)= \Gamma\mc S_{Y|\vX}$,
let $\vu=\sum_{i=1}^d a_i\Gamma\bbeta_i$, then
\begin{align*}
\Gamma\bbeta_i&=\sum_{k=im+1}^{(i+1)m}\theta_{k-im}k^{-(\alpha+\beta)}\phi_k+i^{-\alpha}\phi_i\\
\vu&=\sum_{i=1}^d a_i\Gamma\bbeta_i=\sum_{i=1}^d (\sum_{k=im+1}^{(i+1)m}a_i\theta_{k-im}k^{-(\alpha+\beta)}\phi_k+a_ii^{-\alpha}\phi_i)
\end{align*}
where $a_i$ satisfies
$$
\sum_{i=1}^d a_i^2 (\sum_{k=im+1}^{(i+1)m}k^{-(2\alpha+2\beta)}+i^{-2\alpha})=1.
$$
This implies that $\sum_{i=1}^d a_i^2\geqslant 1/4$.

Then
$$
\langle\Gamma(\vu),\vu\rangle=\sum_{i=1}^d a_i^2[\sum_{k=im+1}^{(i+1)m}k^{-(3\alpha+2\beta)}+i^{-3\alpha}]\geqslant \sum_{i=1}^da_i^2i^{-3\alpha}\geqslant d^{-3\alpha}\sum_{i=1}^da_i^2\geqslant d^{-3\alpha}/4
$$
which means 
$$
\lambda_{\min}(\Gamma|_{\mc S_e})\geqslant d^{-3\alpha}/4.
$$

Next we show that there exist positive constants $\lambda_{-}$ and $\lambda_{+}$ such that $\lambda_{-}\leqslant\lambda_d(\Gamma_e)\leqslant\lambda_1(\Gamma_e)\leqslant \lambda_{+}$.

On the one hand, since $\vX$ is a Gaussian process, the linearity condition holds and $\mr{Im}(\Gamma_e)\subseteq \Gamma\mc S_{Y|\vX}$ \citep{lian2014Sefsdr}. Thus, $\mr{rank}(\Gamma_e)\leqslant d$. 
On the other hand,  let $\vZ=(\vB^*\Gamma\vB)^{-1/2}\boldsymbol{B}^{*}\vX$ and $\Gamma_{ez}=\mr{Cov}(\bbE[\vZ|Y])$. 
Then by Lemma \ref{lem:cov TX}, since  $\vB^*\Gamma \vB=\mr{diag}\{\sum_{k=im+1}^{(i+1)m}k^{-(\alpha+2\beta)}+i^{-\alpha}\}_{i=1}^d$ is positive definite and   invertible,  we have 
\begin{align*}
\vZ \sim N(0,\bs I_{d})\quad\text{and}\quad\Gamma_{ez}=(\vB^*\Gamma\vB)^{-1/2}\boldsymbol{B}^{*}\Gamma_e \boldsymbol{B}   (\vB^*\Gamma\vB)^{-1/2}.    
\end{align*}
 Thus, $\mr{rank}(\Gamma_e)\geqslant\mr{rank}(\Gamma_{ez})=d$. We conclude that $\mr{rank}(\Gamma_e)=\mr{rank}(\Gamma_{ez})=d$, i.e., the coverage condition holds and  $\mr{Im}(\Gamma_e)= \Gamma\mc S_{Y|\vX}$.

For the unit function $\vu\subseteq\mr{Im}(\Gamma_e)= \Gamma\mc S_{Y|\vX}$, we have $\langle\Gamma_e(\vu),\vu\rangle=\mr{var}(\mb E[\langle\vX,\vu\rangle|Y])$.
Then we know that $\mathbb{E}[\langle\boldsymbol{X},\vu\rangle|\boldsymbol{B}^* \boldsymbol{X}]=\vu^*\Gamma\boldsymbol{B}(\boldsymbol{B}^*\Gamma\boldsymbol{B})^{-1}\boldsymbol{B}^*\boldsymbol{X}$ by Lemma \ref{lem:Gaussian property}.

By the law of total expectation, 
$$\mathbb{E}[\langle\boldsymbol{X},\vu\rangle|Y]=\mathbb{E}[\mathbb{E}[ \langle\boldsymbol{X},\vu\rangle|\boldsymbol{B}^{*} \boldsymbol{X}]|Y]=\mathbb{E}[\vu^*\Gamma\boldsymbol{B}(\boldsymbol{B}^*\Gamma\boldsymbol{B})^{-1}\boldsymbol{B}^*\boldsymbol{X} |Y]=\vu^*\Gamma\boldsymbol{B}(\boldsymbol{B}^*\Gamma\boldsymbol{B})^{-1/2}\mathbb{E}[\boldsymbol{Z} |Y],$$
thus $\langle\Gamma_e(\vu),\vu\rangle=\vu^*\Gamma\boldsymbol{B}(\boldsymbol{B}^*\Gamma\boldsymbol{B})^{-1/2}\Gamma_{ez}(\boldsymbol{B}^*\Gamma\boldsymbol{B})^{-1/2}\boldsymbol{B}^{*}\Gamma \vu$ where $\Gamma_{ez}:=\mathrm{var}(\mathbb{E}[\boldsymbol{Z}|Y])$.

\textbf{Upper bound on $\lambda_{1}(\Gamma_e)$.}
We immediately have 
\begin{align*}
  \|\Gamma_e\|=& \|\Gamma\boldsymbol{B}(\boldsymbol{B}^*\Gamma\boldsymbol{B})^{-1/2}\Gamma_{ez}(\boldsymbol{B}^*\Gamma\boldsymbol{B})^{-1/2}\boldsymbol{B}^{*}\Gamma\| \\
  \leqslant& \|\Gamma^{1/2}\boldsymbol{B}(\boldsymbol{B}^*\Gamma\boldsymbol{B})^{-1/2}\Gamma_{ez}(\boldsymbol{B}^*\Gamma\boldsymbol{B})^{-1/2}\boldsymbol{B}^{*}\Gamma^{1/2}\| \|\Gamma^{1/2}\|^2 \leqslant  \|\Gamma\| \| \Gamma_{ez}\|
\end{align*}
where $\Gamma^{1/2}$ is the unique positive operator such that $\Gamma^{1/2}\Gamma^{1/2}=\Gamma$, i.e., 
 $\Gamma^{1/2}:=\sum_{j=1}^\infty j^{-\alpha/2}\phi_j\otimes\phi_j$.

\textbf{Lower bound on $\lambda_{d}(\Gamma_e)$.} 
By  Lemma \ref{lem:minimax operator} (min-max theorem), we have

\begin{align*}
  \lambda_i(\Gamma_e)=& \lambda_i(\Gamma\boldsymbol{B}(\boldsymbol{B}^*\Gamma\boldsymbol{B})^{-1/2}\Gamma_{ez}(\boldsymbol{B}^*\Gamma\boldsymbol{B})^{-1/2}\boldsymbol{B}^{*}\Gamma)\\
   =& \lambda_i(\Gamma^{1/2}\boldsymbol{B}(\boldsymbol{B}^*\Gamma\boldsymbol{B})^{-1/2}\Gamma_{ez}(\boldsymbol{B}^*\Gamma\boldsymbol{B})^{-1/2}\boldsymbol{B}^{*}\Gamma^{1/2}\Gamma)\\
  \geqslant&\lambda_i(\Gamma^{1/2}\boldsymbol{B}(\boldsymbol{B}^*\Gamma\boldsymbol{B})^{-1/2}\Gamma_{ez}(\boldsymbol{B}^*\Gamma\boldsymbol{B})^{-1/2}\boldsymbol{B}^{*}\Gamma^{1/2})  \lambda_{\min}(\Gamma|_{\mc S_e}) \geqslant  \lambda_{\min}(\Gamma|_{\mc S_e})\lambda_i(\Gamma_{ez}). 
\end{align*}

Since  $\|\Gamma\|$,  $d^{-3\alpha}/4$, and the matrix $\Gamma_{ez}$  do not depend on $n$, we conclude the existence of the constants $\lambda_{-}$ and $\lambda_{+}$.

\textbf{Proof of \textbf{ii)}}~

 Next we show that the central curve $\vm(y) = \mb E[\vX|Y = y]$ is weak sliced stable with respect to $Y$. The WSSC for $\mathbb{E}[\vZ|Y=y]$ implies WSSC for $\mathbb{E}[\vX|Y=y]$ since
\begin{align*}
&\frac{1}{H}\sum^{H-1}_{h=0}\var\left(\langle\bs{u},\mb E[\vX|Y]\rangle\mid a_h\leqslant Y\leqslant a_{h+1}\right)\\
=&
\frac{1}{H}\sum^{H-1}_{h=0}\var\left(\langle\bs{u},\Gamma\vB(\vB^*\Gamma\vB)^{-1/2}\mathbb{E}[\vZ |Y]\rangle\mid a_h\leqslant Y\leqslant a_{h+1}\right)\\
=&\frac{1}{H}\sum^{H-1}_{h=0}\var\left(\langle(\vB^*\Gamma\vB)^{-1/2}\vB^*\Gamma\bs{u},\mathbb{E}[\vZ |Y]\rangle\mid a_h\leqslant Y\leqslant a_{h+1}\right)\\
\leqslant &\frac1\tau\var\left(\langle(\vB^*\Gamma\vB)^{-1/2}\vB^*\Gamma\bs{u},\mb E[\vZ|Y]\rangle\right)\\=&\frac1\tau\var\left(\langle\bs{u},\mb E[\vX|Y]\rangle\right)\quad(\forall\bs{u}\in\mathbb{S}_{\mathcal{H}})
\end{align*}
where the inequality comes from the WSSC of $\mathbb{E}[\vZ|Y=y]$ and  the fact that $(\vB^*\Gamma\vB)^{-1/2}\vB^*\Gamma\bs{u}\in\R^d$.

\textbf{Proof of \textbf{iii)}}~
Since $\mb E[X_i^4]=3$, we know that Assumption \ref{as:moment} holds by taking $c_1$ to be $3$. 
Since $\alpha>1,\frac12\alpha+1<\beta,\lambda_j=j^{-\alpha}$ and $|b_{ij}|\lesssim j^{-\beta}$ by the definition of $\bbeta^{\btheta}_{i}$, we know that Assumption \ref{assumption: rate-type condition} holds. Combining  these two results with \textbf{i)} and \textbf{ii)}, we know that \textbf{iii)} holds.
\end{proof}

\subsection{Proof of Lemma \ref{lem:KL divergence close}}
\begin{proof}
For simplicity of notation, we define $\vB:=(\bbeta_1^{\btheta},\dots,\bbeta_d^{\btheta}):\R^d\to L^2[0,1]$ and $\vB':=(\bbeta_1^{\btheta'},\dots,\bbeta_d^{\btheta'}):\R^d\to L^2[0,1]$. Let $\mb E_{\btheta}$ denotes the expectation with respective to $\mb P_{\vB}$ and $\phi_d$ be the density function for $N(0,\bs I_d)$. Then we have
\begin{align*}
\text{KL}(\mathbb{P}_{\btheta},\mathbb{P}_{\btheta'})&\leqslant\mathrm{KL}(\mathbb{P}_{\btheta},\mathbb{P}_{\btheta'})+\mathbb{E}_{(\vX,Y)\sim\mathbb{P}_{\btheta}}\left(\mathrm{KL}(\mathbb{P}_{\btheta}(\bs{Z}\mid \vX,Y),\mathbb{P}_{\btheta'}(\bs{Z}\mid \vX,Y)\right)  \\
&=\operatorname{KL}(\mathbb{P}_{\btheta}(\boldsymbol{X},\boldsymbol{Z},Y),\mathbb{P}_{\btheta'}(\boldsymbol{X},\boldsymbol{Z},Y)) \\
&=\mathbb{E}_{\btheta}[\log\left(\frac{\phi_d(\boldsymbol{Z}-\boldsymbol{B}^{*}\boldsymbol{X})}{\phi_d(\boldsymbol{Z}-\vB^{'*}\boldsymbol{X})}\right)] \\
&=\mathbb{E}_{\btheta}\left(-\frac{1}{2}\|\bs Z-\vB^{*}\vX\|^2+\frac{1}{2}\|\vZ-\vB^{'*}\vX\|^2\right) \\
&=\frac12\mb E[\|(\vB-\vB')^*\vX\|^2     ]\\
 & =\frac12\mb E[\sum_{i=1}^d \langle\bbeta^{\theta}_{i}-\bbeta^{\theta^{\prime}}_i,\vX\rangle^2]\\
&\lesssim m^{-\alpha}\sum^{d}_{i=1}\|\bbeta^{\theta}_{i}-\bbeta^{\theta^{\prime}}_{i}\|^2.
\end{align*}
\end{proof}

\section{Assisting Lemmas}\label{sec:appendix:assist}

\begin{lemma}[Minimax theorem]\label{lem:minimax operator}
Assume that $A$ is a positive semi-definite and compact  operator with its eigenvalues $\{\wt\lambda_i\}$ ordered as $\wt\lambda_1\geqslant\dots\geqslant \wt\lambda_n\geqslant\dots\geqslant 0$, then
$$
\wt\lambda_n=\inf_{E_{n-1}}\sup_{x\in E_{n-1}^\perp,\|x\|=1}\langle Ax,x\rangle
$$
where $E_{n-1}$ with dimension $n-1$ is a closed linear subspace of an Hilbert space $\wt{\mc H}$.
\end{lemma}
It is a classic result in standard functional analysis textbook.

\begin{lemma}\label{lem:Gaussian property}
Assume $\vX=\sum_{i=1}^\infty a_iX_i\phi_i\in\mc H,\vu=\sum_{i=1}^{N}b_i\phi_i\in\mc H,\bbeta_j=\sum_{i=1}^{N}c_{ij}\phi_i\in\mc H,\quad\forall j\in[d]$, then we have
$$
\mathbb{E}[\langle\boldsymbol{X},\vu\rangle|\boldsymbol{B}^* \boldsymbol{X}]=\vu^*\Gamma\boldsymbol{B}(\boldsymbol{B}^*\Gamma\boldsymbol{B})^{-1}\boldsymbol{B}^*\boldsymbol{X}
$$
where $\Gamma=\mb E[\vX\otimes \vX], \vB:=(\bbeta_1,\dots,\bbeta_d):\R^d\to L^2[0,1]$.
\end{lemma}
\begin{proof}
Define $\vX'=\{a_iX_i\}_{i=1}^{N}\in\R^{N},\vu'=\{b_i\}_{i=1}^{N},\bbeta_j'=\{c_{ij}\}_{i=1}^{N},\vB'=(\bbeta_1',\dots,\bbeta_d')\in\R^{p\times d}$.
Then using results from multivariate normal distribution, we have
$$
\mathbb{E}[\langle\boldsymbol{X},\vu\rangle|\boldsymbol{B}^* \boldsymbol{X}]=\bbE[\langle \vX',\vu'\rangle|\vB^{'\top}\vX']=\vu^{'\top}\bs\Sigma \vB'(\vB^{'\top}\bs\Sigma \vB')^{-1}\vB^{'\top} \vX'
$$
where $\Sigma$ is the covariance matrix of $\vX'$.

We then complete the proof by using the following relationships:  
$$
\vu^{'\top}\bs\Sigma \vB'=\vu^*\Gamma\boldsymbol{B},(\vB^{'\top}\bs\Sigma \vB')^{-1}=(\boldsymbol{B}^*\Gamma\boldsymbol{B})^{-1},\vB^{'\top} \vX'=\boldsymbol{B}^*\boldsymbol{X}.
$$
\end{proof}

\begin{lemma}\label{lem:cov TX}
If $T$ is an operator defined on $\mc H_1\to\mc H_2$ where $\mc H_i,i=1,2$ is a Hilbert space. $\vX\in\mc H_1$ is a random element satisfying $\mb E[\vX]=0$ . Then we have $\mr{var}(T\vX)=T\mr{var}(\vX)T^*$.
\end{lemma}
\begin{proof}
For any $\u_1,\u_2\in\mc H_2$, we have
\begin{align*}
&\left\langle  T\mr{var}(\vX)T^*\u_1,\u_2  \right\rangle=\left\langle  T\mb E[\vX\otimes\vX]T^*\u_1,\u_2  \right\rangle
=\left\langle  \mb E[\vX\otimes\vX]T^*\u_1,T^*\u_2  \right\rangle    
\end{align*}
since $\mb E[\vX]=0$. By the definition of convariance operator and expectation, we have 
\begin{align*}
\left\langle  \mb E[\vX\otimes\vX]T^*\u_1,T^*\u_2  \right\rangle=&\left\langle  \mb E[\left\langle\vX,  T^*\u_1 \right\rangle       \vX            ],T^*\u_2  \right\rangle
=\mb E[  \left\langle\vX,  T^*\u_1 \right\rangle      \left\langle \vX            ,T^*\u_2  \right\rangle].
\end{align*}
Similarly, we have
\begin{align*}
 \left\langle  \mr{var}(T\vX)\u_1,\u_2  \right\rangle=\left\langle  \mb E[T\vX\otimes T\vX]\u_1,\u_2  \right\rangle=\mb E[  \left\langle T\vX,  \u_1 \right\rangle      \left\langle T\vX            ,\u_2  \right\rangle].\\    
\end{align*}
Then the proof is completed by noticing the following
\begin{align*}
\mb E[  \left\langle T\vX,  \u_1 \right\rangle      \left\langle T\vX            ,\u_2  \right\rangle]=\mb E[  \left\langle\vX,  T^*\u_1 \right\rangle      \left\langle \vX            ,T^*\u_2  \right\rangle].
\end{align*}
\end{proof}
\begin{lemma}\label{lem:PimTPimtoT}If $T$ is of finite rank, then we have $\lim\limits_{m\to \infty}\|\Pi_m T-T\|=0$.
\end{lemma}
\begin{proof}By the triangle inequality and compatibility of operator norm, one has
\begin{align*}
\|\Pi_m T-T\|\leqslant\|(\Pi_m-I)T\|
\end{align*}
where $I=\sum\limits_{i=1}^\infty\phi_i\otimes\phi_i$ for $\{\phi_i\}_{i\in\mb{Z}_{\geqslant 1}}$ being an orthonormal basis of $\mc H$. 

Since $T$ is of finite rank, let us assume that $\{e_i\}_{i=1}^k$ is an orthonormal basis of $\mathrm{Im}(T)$ where $k=\mr{rank}(T)$. For any $\beta\in\mathcal{H}$ such that $\|\beta\|=1$, one has $\|T\beta\|\leqslant\|T\|\|\beta\|=\|T\|$, so one can assume that $T\beta\in\mathrm{Im}(T)$ admits the following expansion under basis $\{e_i\}_{i=1}^k$:
\[T\beta=\sum_{i=1}^k b_ie_i,\quad \sum_{i=1}^k b^2_i\leqslant\|T\|^2<\infty.\]
Thus
\[\|(I-\Pi_m)T\beta\|=\left\|\sum_{i=1}^k(I-\Pi_m) b_ie_i\right\|\leqslant\sum_{i=1}^k |b_i|\cdot\|(I-\Pi_m) e_i\|.\]
Clearly, $\|(\Pi_m-I)\alpha\|~(\forall\alpha\in\H)$ tends to $0$ as $m\to\infty$ since 
\[(I-\Pi_m)\alpha=\left(\sum_{i={m+1}}^\infty\phi_i\otimes\phi_i\right)\left(\sum\limits_{i=1}^\infty c_i\phi_i\right)=\sum_{i=m+1}^\infty c_i\phi_i\xrightarrow{m\to\infty} 0\]
where we have assumed that $\alpha=\sum\limits_{i=1}^\infty c_i\phi_i$ .

Thus $\forall\varepsilon>0$, there exists some $N_i>0$ such that $\forall m> N_i$ one has $\|(\Pi_m-I)e_i\|<\varepsilon$, $(\forall i=1,...,k)$. Let $N=\max\{N_1,\cdots,N_k\}$, then $\forall m>N$ one has
\[\|(I-\Pi_m)T\beta\|\leqslant\sum_{i=1}^k |b_i|\cdot\|(I-\Pi_m) e_i\|\leqslant\sum_{i=1}^k |b_i|\varepsilon\leqslant k\varepsilon\|T\|,\]
which means that $\forall m>N$, one has
\begin{align*}
\|(\Pi_m-I)T\|&=\sup_{\|\beta\|=1}\|(\Pi_m-I)T\beta\|\leqslant k\varepsilon\|T\|. 
\end{align*}
Thus $\lim\limits_{m\to\infty}\|(\Pi_m-I)T\|=0$. 
 Then the proof of Lemma $\ref{lem:PimTPimtoT}$ is completed.
\end{proof}
\subsection{Properties of sliced partition}
\begin{lemma}[Corollary 1 in \cite{lin2018supplement}]\label{cor:i.i.d}
In the slicing inverse regression contexts, recall that $S_{h}$ denotes the h-th interval $(y_{h-1,c},y_{h,c}]$ for $2\leqslant h \leqslant H-1$ and $S_{1}=(-\infty, y_{1,c}]$, $S_{H}=(y_{H-1,c},\infty)$. We have that
$x_{h,i}, i=1,\cdots,c-1$ can be treated as $c-1$  random samples of $x\Big|(y \in S_{h})$ for $h=1,...,H-1$ and $x_{H,1},...,x_{H,c}$ can be treated as $c$ random samples of $x\Big|(y\in S_{H})$. 
\end{lemma}
\begin{lemma}[Lemma $11$ in \cite{lin2018supplement}]\label{lem:sliced}
		For any  sufficiently large $H,c$ and  $n>\frac{4H}{\gamma}+1$, the sliced partition $\mathfrak{S}_{H}(n)$ is a $\gamma$-partition with probability at least 
		\[1-CH^2\sqrt{n+1}\exp\left(-\frac{\gamma^2(n+1)}{32H^2}\right)\]
		for some absolute constant $C$.
\end{lemma}
\begin{lemma}[Lemma $10$ in \cite{lin2018supplement}]\label{lem:fourth_moment:lemma} 
Suppose that $(x,y)$ are defined over $\sigma$-finite space $\mathcal{X}\times\mathcal{Y}$ and $g$ is a non-negative function such that $\E[g(x)]$ exists.  For any fixed positive constants $C_{1}<1<C_{2}$, there exists a constant $C$ which only depends on $C_{1}, C_{2}$ such that for any partition 
$\mathbb{R}=\bigcup_{h=1}^{H}S_{h} $ where $S_{h}$ are intervals  satisfying
\begin{align}\label{partition:condition}
  \frac{C_{1}}{H} \leqslant \P(y \in S_{h})  \leqslant \frac{C_{2}}{H}, \forall h,
\end{align}
we have
\[
\sup_h\E(g(x)\Big|{y\in S'_{h}}) \leqslant  CH\E[g(x)].
\]
\end{lemma}
\subsection{Sin Theta Theorem}\label{ap:Sin Theta theorem}
\begin{lemma}[Proposition 2.3 in \cite{Seelmann2014NotesOT}]\label{lemma, sin theta of infinite dimension operator}
Let $B$ be a self-adjoint operator on a separable Hilbert space $\widetilde{\mathcal{H}}$, and let ${V}\in\mathcal{L}(\widetilde{\mathcal{H}})$ be another self-adjoint operator where $\mathcal{L}\left(\widetilde{\mc H}\right)$ stands for the space of bounded linear operators from a Hilbert space $\widetilde{\mc H}$ to $\widetilde{\mc H}$. 
Write the spectra of $B$ and $B+V$ as \[\mathrm{spec}( B)=\sigma\cup\Sigma\quad\text{and}\quad \mathrm{spec}( B+ V)=\omega\cup\Omega
\]
with $\sigma\cap\Sigma=\varnothing=\omega\cap\Omega$, and suppose that there is $\widehat d>0$ such that
\[\mathrm{dist}(\sigma,\Omega)\geqslant \widehat d\quad\text{and}\quad\mathrm{dist}(\Sigma,\omega)\geqslant \wh d\]
where $\mathrm dist(\sigma,\Sigma):=\min\{|a-b|:a\in\sigma,b\in\Omega\}$. 
Then, it holds that
\[\|P_{{B}}(\sigma)-P_{{B}+{V}}(\omega)\|\leqslant\frac\pi2\frac{\| V\|}{\wh d}\]
where $P_{ B}(\sigma)$ denotes the spectral projection for $ B$ associated with $\sigma$, i.e., 
\[P_{B}(\sigma):=\frac{1}{2\pi\mathrm{i}}\oint_{\gamma}\frac{\mathrm{d}z}{z-B},\]
where $\gamma$ is a contour on $\mathbb{C}$ that encloses $\sigma$ but no other elements of $\mathrm{spec}( B)$.
\end{lemma}

\begin{remark}[Spectral projection]
We note that, 
if further $ B$ is compact, 
the spectral projection $P_{B}(\sigma)$ coincide with the projection operator onto the closure of the space spanned by the eigenfunctions associated with the eigenvalues in $\sigma$.

Specifically, if $B$ is compact, by the spectral decomposition theorem one has
\[B=\sum_{i=1}^\infty\mu_ie_i\otimes e_i\quad\text{and}\quad(z- B)^{-1}=\sum_{i=1}^\infty(z-\mu_i)^{-1}e_i\otimes e_i,\]
where $\mr{spec}(B):=\{\mu_i\}_{i=1}^\infty$ satisfies $|\mu_i|\xrightarrow{i\to\infty} 0$.
Then $\forall v\in \mathcal{H}$, it holds that
\begin{align*}P_{B}(\sigma)v&=\frac{1}{2\pi\mathrm{i}}\oint_{\gamma}({z-B})^{-1}v~{\mathrm{d}z}=\frac{1}{2\pi\mathrm{i}}\oint_{\gamma}\sum_{i=1}^\infty(z-\mu_i)^{- 1}\langle e_i,v\rangle e_i~{\mathrm{d}z}\\
&=\sum_{i=1}^\infty\left[\left(\frac{1}{2\pi\mathrm{i}}\oint_{\gamma}(z-\mu_i)^{-1}~{\mathrm{d}z}\right)\langle e_i,v\rangle e_i\right]=\sum_{i\in\{i:\mu_i\in\sigma\}}\langle e_i,v\rangle e_i.
\end{align*}

In particular, if $\sigma=\mr{spec}(B)\backslash\{0\}$, then $P_{B}(\sigma)$ is the projection operator onto the $\overline{\mathrm{Im}}(B)$.
\end{remark}

\subsection{Wely inequality for a self-adjoint and compact operator}\label{ap:Wely inequality for self-adjoint and compact operators}
\begin{lemma}\label{lem:wely ineq operator}
Let $M$ and $N$ be two self-adjoint, positive semi-definite and compact operators defined on a Hilbert space $\wt{\mc H}$ with their respective eigenvalues $\{\mu_i\},\{\nu_i\}$ ordered as follows
\begin{align*}
M:\mu_1\geqslant\dots\geqslant \mu_n\geqslant\dots\geqslant 0\quad\text{and}\quad
N:\nu_1\geqslant\dots\geqslant \nu_n\geqslant\dots\geqslant 0.
\end{align*}
Then the following inequalities hold: $|\mu_k-\nu_k|\leqslant\|M-N\|$, $ k\geqslant1$.
\end{lemma}

\section{Additional Simulation Results}\label{sec:appendix:simulation}


This section contains additional simulation results
of Section \ref{sec:simulationfsir}. Specifically, in Section \ref{sec:appe:explanation}, we give an  
explanation of why
 $\bs X$ in Model I to III  is equivalent to a construction that satisfies the assumption that $\Gamma$ is non-singular. In Section \ref{sec:guidelines H}, we provide guidelines for the choice of $H$ in practice. Sections \ref{sec:app:optimal m} and \ref{sec:app: synthetic data} contain
additional  simulation results  of Sections \ref{sec:experiment optimal m} and \ref{sec: synthetic data} respectively when 
(i): $\varepsilon\sim N(0,1),H=10$; (ii):$\varepsilon\sim N(0,0.25),H=10$; (iii):$\varepsilon\sim N(0,2),H=10$.
Lastly, we compare FSIR  with PCA on selecting optimal $m$ in Section \ref{sec:compare PCA}.

\subsection{The explanation on difference between  simulation models and the theoretical assumptions}
\label{sec:appe:explanation}
Here, we explain why
 the construction of $\bs X$ in simulation models do not contradict the assumption that $\Gamma$ is non-singular.

Assume $\vX=\sum_{i=1}^\infty a_iX_i\phi_i\in\H \ \text{and } \bbeta_j=\sum_{i=1}^{N}c_{ij}\phi_i\in\H, \ 
\forall j\in[d]$. Then, we have
$$ \begin{aligned}
\E[\vX|f(\vB^*\vX,\epsilon)]
=&\E[\sum_{i=1}^N a_iX_i\phi_i|f(\vB^*\vX,\epsilon)]+\E[\sum_{i=N+1}^\infty a_iX_i\phi_i|f(\vB^*\vX,\epsilon)]\\
=&\E[\sum_{i=1}^N a_iX_i\phi_i|f(\vB^*\vX,\epsilon)]+\E[\sum_{i=N+1}^\infty a_iX_i\phi_i]\\
=&\E[\sum_{i=1}^N a_iX_i\phi_i|f(\vB^*\vX,\epsilon)].
\end{aligned}
$$
Thus, in terms of $\Gamma_e$, 
the truncation on $\bbeta_j$ can be transferred to the truncation on $\vX$. 
Note that the SIR estimate only involves $\Gamma_m$ and  $\Gamma_e$.
This suggests that when $\vX$ has infinitely many terms, the SIR estimate remains the same before and after we do truncation on $\vX$ as long as $m$ is smaller than $N$. 
Therefore, we directly simulate the truncated version of $\vX$ in Model I to III. 

\subsection{Guidelines on selecting $H$}\label{sec:guidelines H}
In practical scenarios, the selection of the slices number $H$ can be influential. 
Our simulation studies suggest that selecting $H \geq \ln(n)$ is often sufficient to achieve desirable numerical results and this selection meets the theoretical requirement that $H > H_0$
since $\ln(n)$ will eventually exceed the constant $H_0$ defined in Theorem 2.
Furthermore, we recommend selecting $H$ within the range $[10, 35]$ to accommodate finite-sample scenarios in practice. 

We conducted a series of experiments to substantiate the practicality of this guideline. 
Specifically, we simulated the FSIR process as described in Section 4.3 and  determined the minimum average subspace estimation error across $100$ repetitions for various values of $m$ in  the set $\{2, 3, \ldots, 13, 14, 20, 30, 40\}$. 
We first set $n=20,000$. 
The value of $H$ is initiated at a baseline value $\ln(n)=\ln(20000) \approx 10$ and increased in increments of 5. 
The results are presented in Table \ref{tab:different H}, which illustrates that the subspace estimation error is relatively insensitive to variations in $H$ as long as $H \geq \ln(n)$. For comparative purposes, at an exceptionally low $H$ value, such as $H=2$, the errors for the three models are recorded at 0.077, 0.291, and 0.02, respectively, which are significantly higher than those obtained for $H \geq 10$. 
This comparison highlights the necessity of adhering to the guideline of $H \geq \ln(n)$ for robust model performances. 
We then expanded the dataset sizes to $n=50,000$ and $n=200,000$ respectively and conducted the experiments with $H$ chosen in the same way as before; the results are consistent with earlier findings. These empirical results substantiate the practicality of our proposed guideline.

\begin{table}[H]
\begin{center}
\renewcommand\arraystretch{1.5}
\begin{tabular}{|c|c|c|c|c|c|c|c|}
\hline
$ \text{FSIR-OT}$ & $H$  & $2$  & $10$ & $15$ & $20$ & $25$  & $30$   \\
\hline
\hline
\multirow{3}{*}{$n=20000$}
& Model I & \textbf{0.077} & \textbf{0.067} & \textbf{0.06} & \textbf{0.066} & \textbf{0.066} & \textbf{0.065} \\
& Model II & \textbf{0.291} & \textbf{0.024} & \textbf{0.03} & \textbf{0.026} & \textbf{0.028} & \textbf{0.026} \\
& Model III& \textbf{0.02}  & \textbf{0.01} & \textbf{0.01} & \textbf{0.01} & \textbf{0.01} & \textbf{0.01} 
\\[1pt]
\hline
$ $ & $H$  & $2$  & $11$ & $16$ & $21$ & $26$  & $31$   \\
\hline
\multirow{3}{*}{$n=50000$}
& Model I & \textbf{0.060} & \textbf{0.051} & \textbf{0.052} & \textbf{0.053} & \textbf{0.052} & \textbf{0.050} \\
& Model II & \textbf{0.464} & \textbf{0.016} & \textbf{0.015} & \textbf{0.015} & \textbf{0.016} & \textbf{0.016} \\
& Model III& \textbf{0.013}  & \textbf{0.010} & \textbf{0.009} & \textbf{0.009} & \textbf{0.010} & \textbf{0.009} 
\\[1pt]
\hline
$ $ & $H$  & $2$  & $13$ & $18$ & $23$ & $28$  & $33$   \\
\hline
\multirow{3}{*}{$n=200000$}
& Model I & \textbf{0.041} & \textbf{0.036} & \textbf{0.035} & \textbf{0.036} & \textbf{0.036} & \textbf{0.036} \\
& Model II & \textbf{0.338} & \textbf{0.008} & \textbf{0.007} & \textbf{0.008} & \textbf{0.008} & \textbf{0.008} \\
& Model III& \textbf{0.007}  & \textbf{0.004} & \textbf{0.004} & \textbf{0.005} & \textbf{0.004} & \textbf{0.005} 
\\[1pt]
\hline
\end{tabular}
\caption{
The minimum average subspace estimation error over $100$ repeated experiments with respect to different $m$ of FSIR-OT for various models in the case of $\varepsilon\sim N(0,2)$. 
}\label{tab:different H}
\end{center}
\end{table}

\subsection{Additional simulation results of Section \ref{sec:experiment optimal m}}\label{sec:app:optimal m}

The left panels of Figures \ref{fig:optimal m noise1} - \ref{fig:optimal m noise2} are the average subspace estimation error under Model (I) where $n$ ranges in $\{2\times 10^3,2\times 10^4,5\times10^4,2\times 10^5,5\times 10^5,10^6\}$, $m$ ranges in $\{3,4,\dots,25\}$. The optimal value of $m$ (denoted by $m^*$) for each $n$ is marked with a red circle.
The shaded areas represent the standard error bands associated with these estimates (all smaller than  $0.011$). 
The right panel of Figures~\ref{fig:optimal m noise1} - \ref{fig:optimal m noise2} illustrate the linear dependence of  $\log(m^*)$ on $\log(n)$. The solid line characterizes the linear trend of $\log(m^*)$ against $\log(n)$. The dotted lines are their least-squares fittings, with their slopes estimated as $0.183$, $0.2$ and $0.196$ respectively, which are close to the theoretical value of $2/11$. These results are consistent with the theoretically optimal choice of $m$ in FSIR-OT.
\begin{figure}[H]
	\centering
	\begin{minipage}{0.50\textwidth}
		\includegraphics[width=\textwidth]{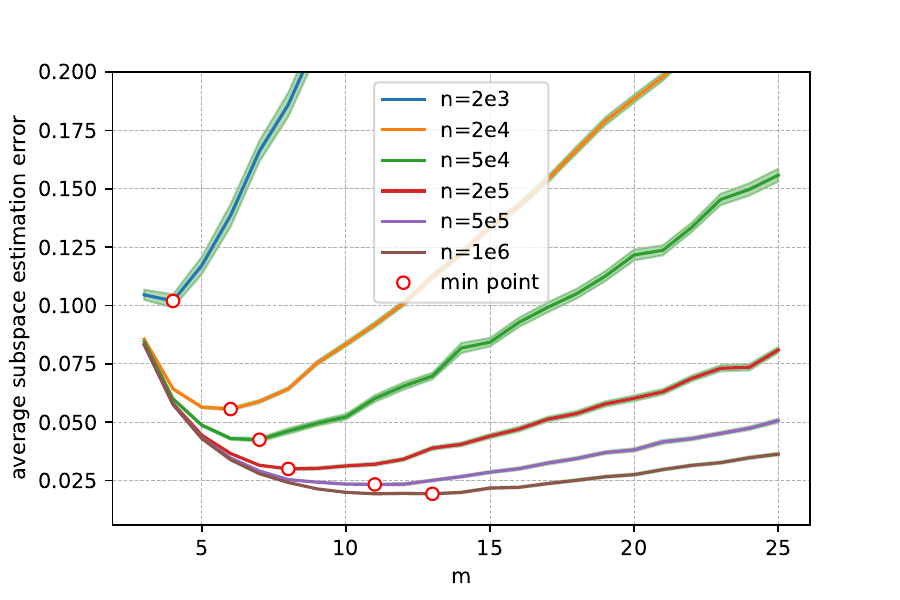}
	\end{minipage}\hfill
	\begin{minipage}{0.50\textwidth}
		\includegraphics[width=\textwidth]{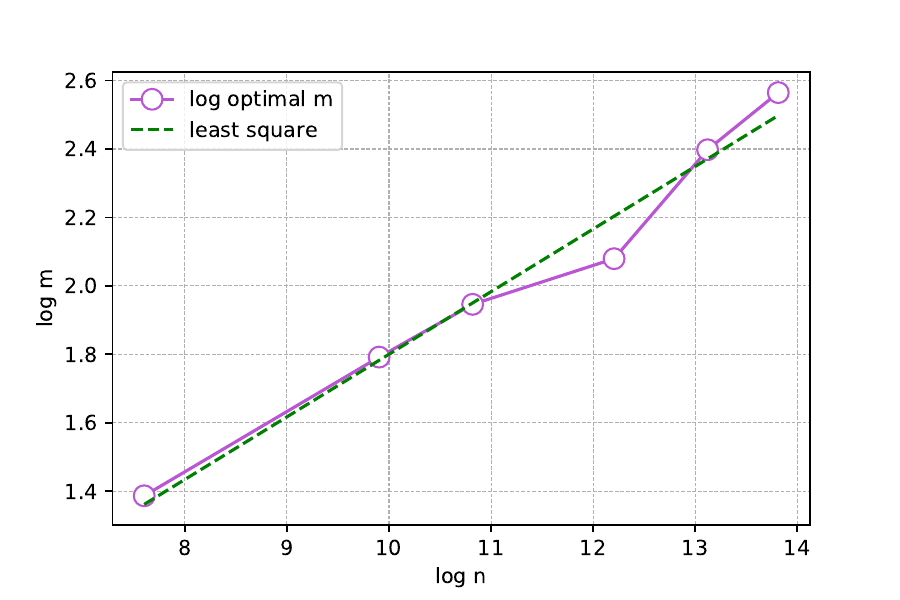}
	\end{minipage}
\caption{
Experiments for the optimal choice of truncation parameter $m$ with $\varepsilon\sim N(0,1)$ and $H=10$. 
Left: average subspace estimation error with increasing $m$ for different $n$. 
Right: linear trend of $\log(m^*)$ against $\log(n)$, with a slope of $0.183$  and $R^2>0.98$.}
\label{fig:optimal m noise1}
\end{figure}
\begin{figure}[H]
	\centering
	\begin{minipage}{0.50\textwidth}
		\includegraphics[width=\textwidth]{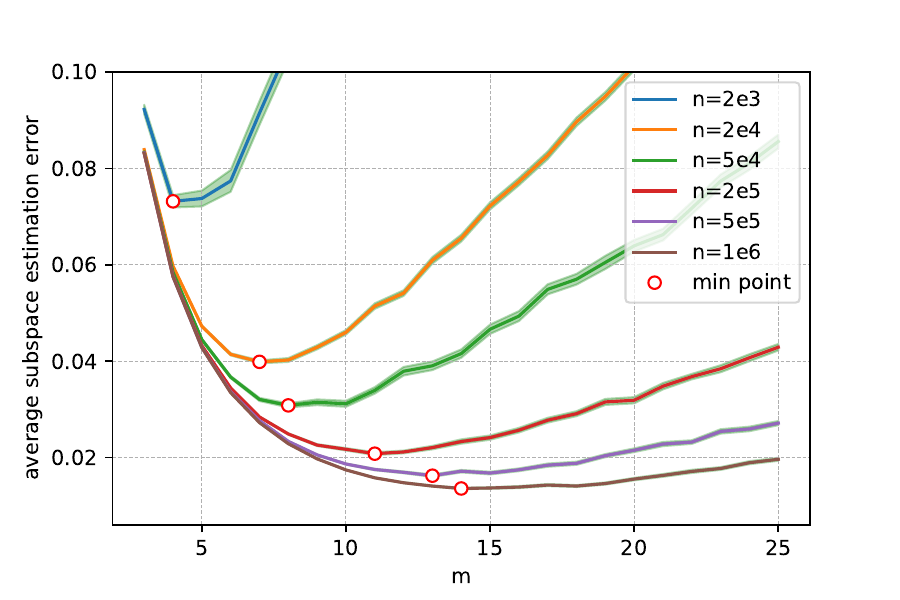}
	\end{minipage}\hfill
	\begin{minipage}{0.50\textwidth}
		\includegraphics[width=\textwidth]{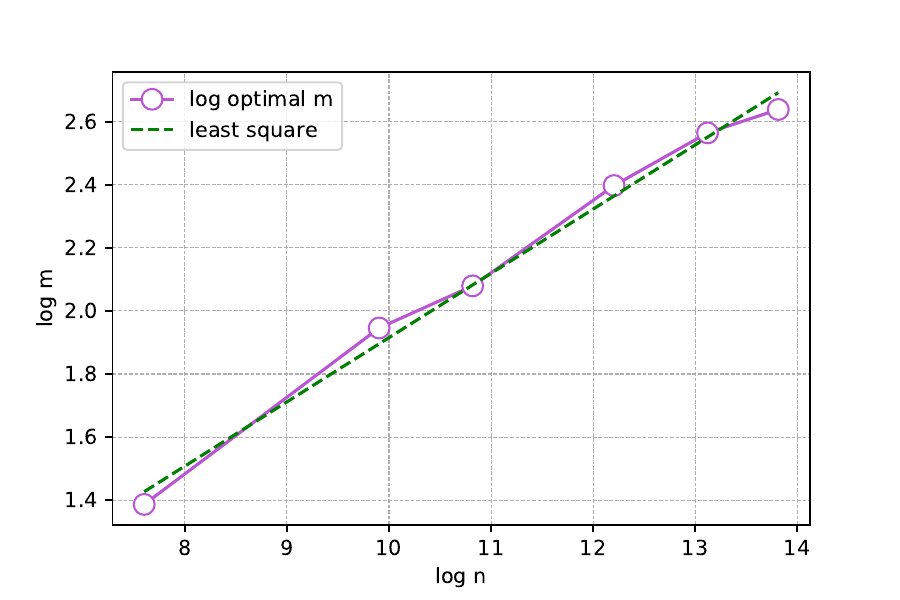}
	\end{minipage}
\caption{Experiments for the optimal choice of truncation parameter $m$ with  $\varepsilon\sim N(0,0.25)$ and $H=10$. Left: average subspace estimation error with increasing $m$ for different $n$. 
Right: linear trend of $\log(m^*)$ against $\log(n)$, with a slope of $0.2$  and $R^2>0.99$.}
\label{fig:optimal m noise0.5}
\end{figure}

\begin{figure}[H]
	\centering
	\begin{minipage}{0.50\textwidth}
		\includegraphics[width=\textwidth]{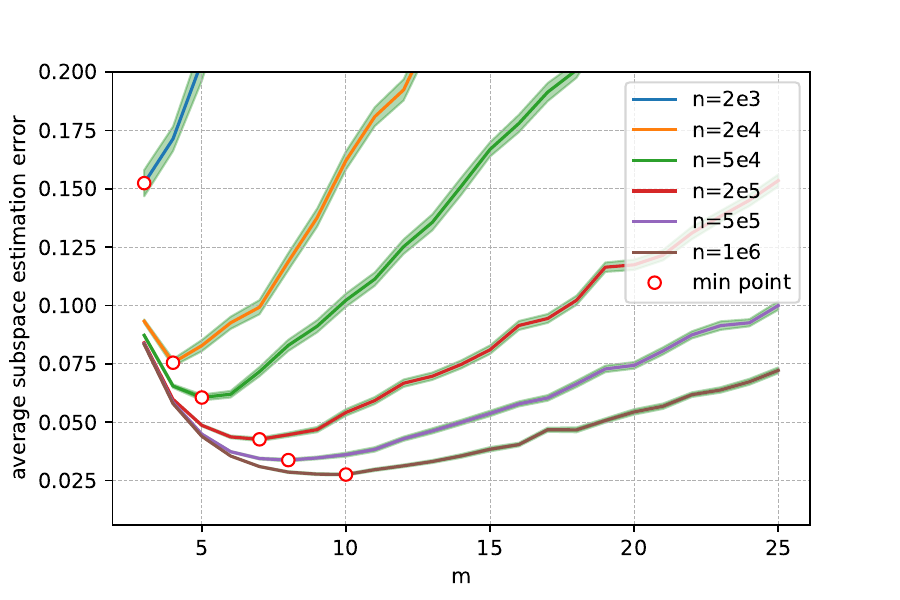}
	\end{minipage}\hfill
	\begin{minipage}{0.50\textwidth}
		\includegraphics[width=\textwidth]{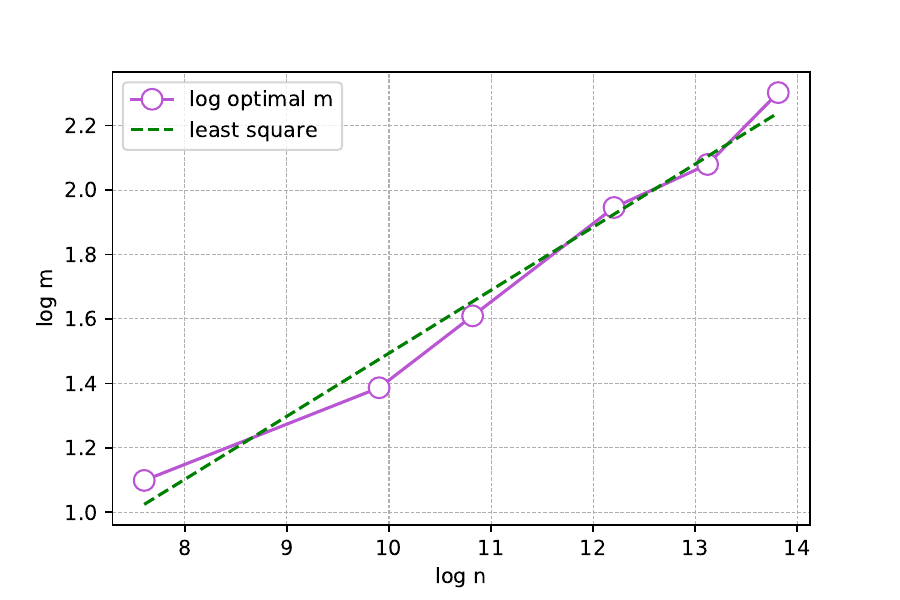}
	\end{minipage}
\caption{Experiments for the optimal choice of truncation parameter $m$ with  $\varepsilon\sim N(0,2)$ and $H=10$.
 Left: average subspace estimation error with increasing $m$ for different $n$. 
Right: linear trend of $\log(m^*)$ against $\log(n)$, with a slope of $0.196$  and $R^2>0.98$.}
\label{fig:optimal m noise2}
\end{figure}

\subsection{Additional simulation results of Section \ref{sec: synthetic data}}\label{sec:app: synthetic data}

For each model of Models (I) - (III), we calculate the average subspace estimation error of FSIR-OT and RFSIR 
based on $100$ replications, where $n=20000$, the truncation parameter of FSIR-OT $m$ ranges in $\{2,3,\dots,13,14,20,30,40\}$, and the regularization parameter in RFSIR $\rho$ ranges in $0.01\times \{1,2,\cdots,9,10,15,20,25,30,40,\cdots,140,150\}$. Detailed results are presented in Figures \ref{fig:error 3models noise1} - \ref{fig:error 3models noise2}, 
where we mark the minimal error in each model with red `$\times$'.
The shaded areas represent the corresponding standard errors, all of which are less than $0.01$. When $\varepsilon\sim N(0,1), H=10$, 
for FSIR-OT, the  minimal errors for $\mc M_1$, $\mc M_2$, and $\mc M_3$ are  $0.06,0.02$, and $0.01$ 
respectively.  For RFSIR,  the corresponding minimal errors are $0.08,0.06$, and $0.01$.
When $\varepsilon\sim N(0,0.25), H=10$, for FSIR-OT, the  minimal errors for $\mc M_1$, $\mc M_2$, and $\mc M_3$ are  $0.04,0.02$, and $0.01$ 
respectively. For RFSIR,  the corresponding minimal errors are $0.06,0.04$, and $0.01$.
When $\varepsilon\sim N(0,2),H=10$, for FSIR-OT, the  minimal errors for $\mc M_1$, $\mc M_2$, and $\mc M_3$ are  $0.078,0.032$, and $0.017$ 
respectively. For RFSIR,  the corresponding minimal errors are $0.109,0.105$, and $0.014$.
The results here suggest that the performance of FSIR-OT is generally superior to, or at the very least equivalent to, that of the RFSIR. 

\begin{figure}[H]
	\centering
	\begin{minipage}{0.50\textwidth}
		\includegraphics[width=\textwidth]{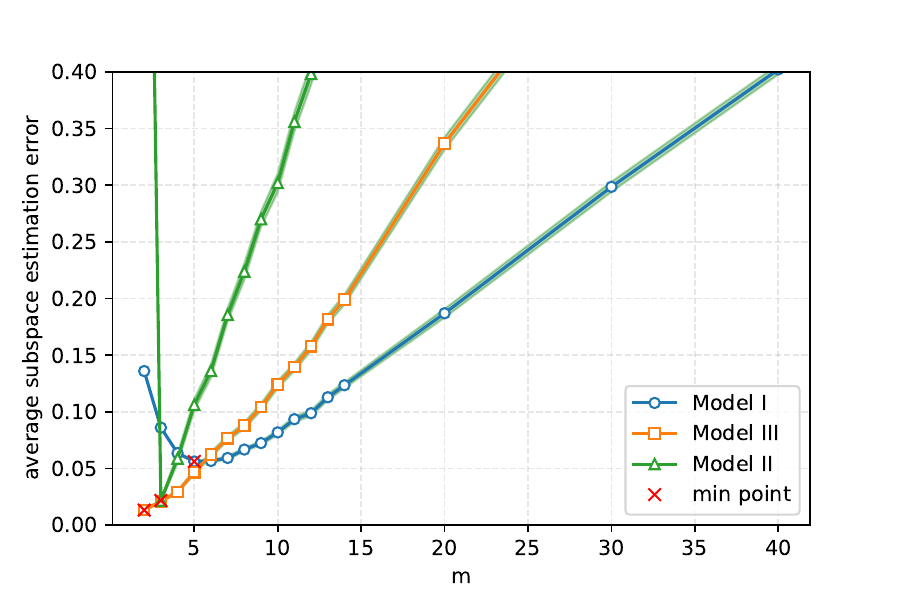}
	\end{minipage}\hfill
	\begin{minipage}{0.50\textwidth}
		\includegraphics[width=\textwidth]{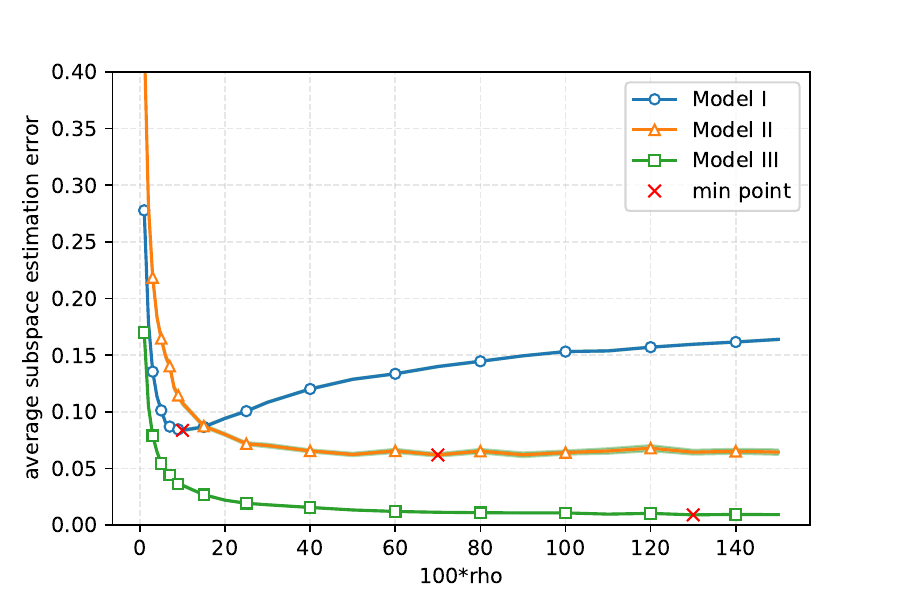}
	\end{minipage}
\caption{
 Average subspace estimation error of FSIR-OT and RFSIR for various models in the case of $\varepsilon\sim N(0,1)$ and $H=10$. The standard errors are all below $0.01$. 
Left: FSIR-OT with different truncation parameter $m$; Right: RFSIR with different values of the regularization parameter $\rho$. 
}
\label{fig:error 3models noise1}
\end{figure}

\begin{figure}[H]
	\centering
	\begin{minipage}{0.50\textwidth}
		\includegraphics[width=\textwidth]{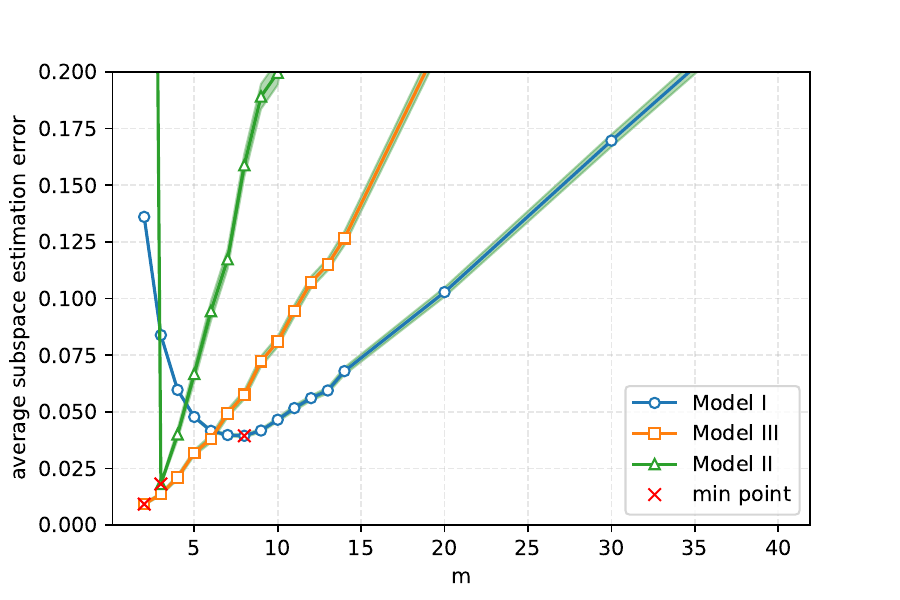}
	\end{minipage}\hfill
	\begin{minipage}{0.50\textwidth}
		\includegraphics[width=\textwidth]{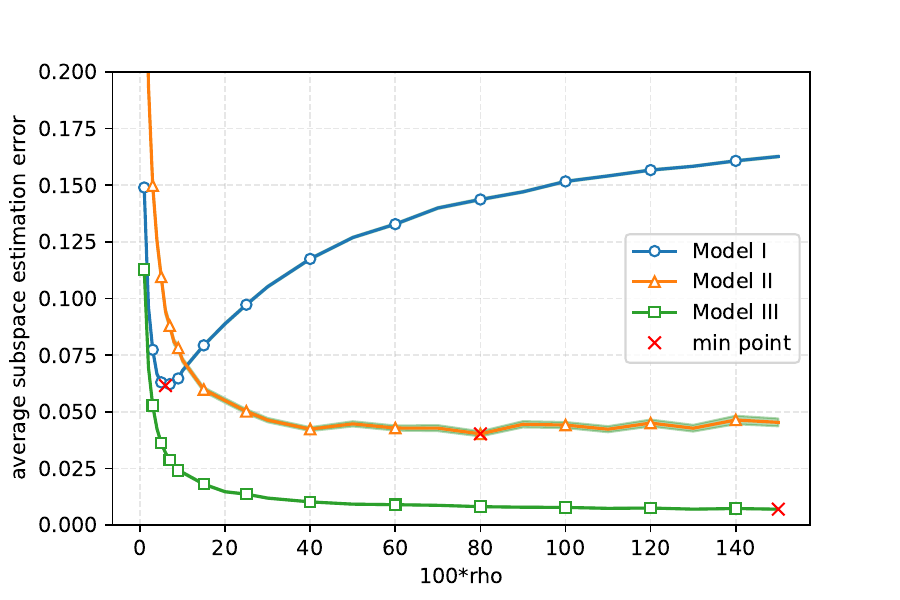}
	\end{minipage}
\caption{
 Average subspace estimation error of FSIR-OT and RFSIR for various models in the case of $\varepsilon\sim N(0,0.25)$ and $H=10$. Average subspace estimation error of FSIR-OT and RFSIR for various models. The standard errors are all below $0.008$. 
Left: FSIR-OT with different truncation parameter $m$; Right: RFSIR with different values of the regularization parameter $\rho$. 
}
\label{fig:error 3models noise0.5}
\end{figure}

\begin{figure}[H]
	\centering
	\begin{minipage}{0.50\textwidth}
		\includegraphics[width=\textwidth]{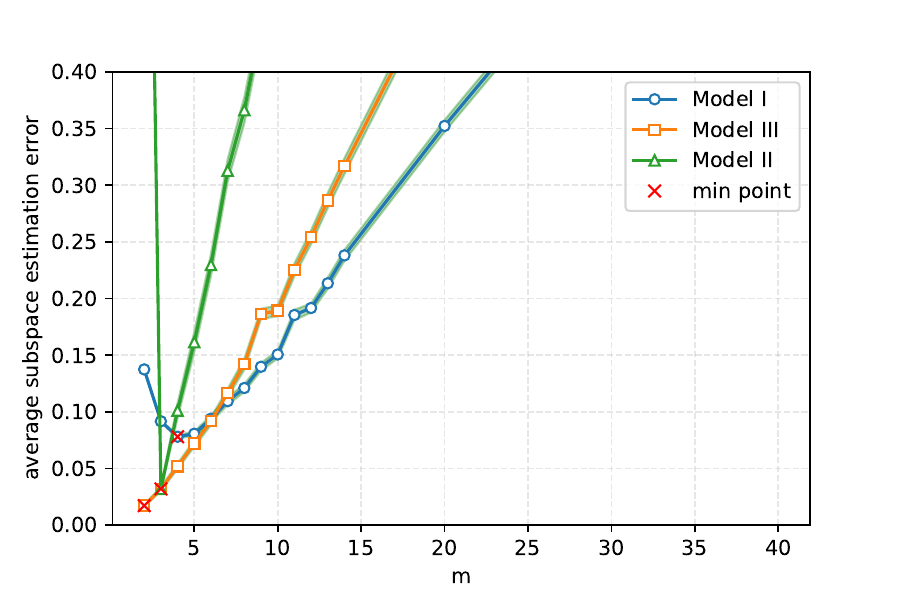}
	\end{minipage}\hfill
	\begin{minipage}{0.50\textwidth}
		\includegraphics[width=\textwidth]{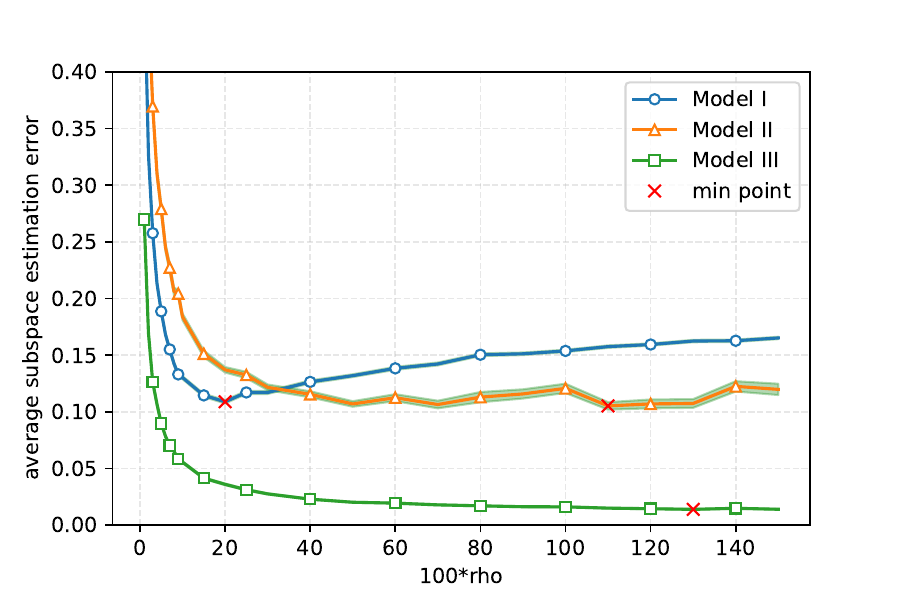}
	\end{minipage}
\caption{
 Average subspace estimation error of FSIR-OT and RFSIR for various models in the case of $\varepsilon\sim N(0,2)$ and $H=10$.
Average subspace estimation error of FSIR-OT and RFSIR for various models. The standard errors are all below $0.012$. 
Left: FSIR-OT with different truncation parameter $m$; Right: RFSIR with different values of the regularization parameter $\rho$. 
}
\label{fig:error 3models noise2}
\end{figure}

\subsection{Comparison with PCA on selecting optimal $m$}\label{sec:compare PCA}
SIR is a supervised learning method, whereas PCA is an unsupervised learning method. Consequently, using PCA to select 
$m$ without any information from the response variable is intuitively incorrect. 

Specifically, we attempt to use PCA to select 
$m$ and highlight the drawbacks of this approach in the following. 
Here we consider the bike sharing data set studied in Section 4.4. 
First, we calculate the proportion of the total eigenvalue sum explained by the first 
$i$ ($i \leq14$) eigenvalues of $\wh\Gamma$. We find that the first eigenvalue alone accounts for over $99.4\%\geq99\%$ of the total, suggesting that choosing 
$m$ as $1$ is a good option, and increasing $m$ further adds little value. Under this parameter selection (with  $d$ only taking $1$), the Gaussian process regression error after SIR dimension reduction is $0.201$, significantly higher than the optimal result $0.188$ in Table \ref{tab:real data error}, corresponding to $d=2$ and $m=6$. Therefore, using the unsupervised learning method PCA does not provide a reasonable way to select $m$.

\end{document}